\documentclass[12pt,diff1_eng]{article} 
\usepackage[english]{babel}
\usepackage{amssymb,amscd,amsmath}     
\makeatletter
\input{DIFF1_ENG.STY}
\makeatother
\usepackage{graphicx}
\usepackage{caption}
\usepackage[linktocpage=true, colorlinks=true, linkcolor=blue, citecolor=blue, urlcolor=blue]{hyperref}

\begin{document}
  
\large

$$
{}
$$

\begin{center}
\Large{\bf A New Modifications of the Milstein and Euler Methods and the Taylor--It\^{o} Expansion
for It\^{o} Stochastic Differential 
Equations with Drift and Diffusion Coefficients Having Non-Differentiable or 
Singular Time-Dependent Multipliers}
\end{center}

$$
{}
$$

\begin{center}
{\large {Dmitriy F. Kuznetsov}}\\
Peter the Great Saint-Petersburg Polytechnic University\\
\hspace{-2mm}e-mail: sde\_\hspace{0.5mm}kuznetsov@inbox.ru
\end{center}

\begin{center}
{\large {Mikhail D. Kuznetsov}}\\
Saint-Petersburg
Federal Research Center of the Russian Academy of Sciences\\
(SPC RAS)\\
International Digital Forensic Center
\end{center}

$$
{}
$$

\noindent
{\bf Abstract.}  The article presents a new modification of the Milstein method for systems 
of It\^{o} stochastic differential equations with 
drift and diffusion coefficients having differentiable
and non-differentiable or 
singular power-law time-dependent components
(the latter depending only on 
the time variable). By iteratively applying It\^{o}'s formula to the mentioned 
differentiable components of the drift and diffusion coefficients, a new 
modification of the Milstein method is derived in the paper. A key feature 
of this modification of the Milstein method is the presence of time-dependent  
weight functions in the iterated It\^{o} stochastic integrals of the second and first
multiplicity (the specified weight functions are identically equal to unity in the 
classical Milstein method). To approximate the specified iterated It\^{o} stochastic 
integrals (with respect to the components of a multidimensional Wiener process) that 
involve time-dependent weight functions, the paper proposes using the method of 
generalized multiple Fourier series, introduced by the first author in 2006.
The article also provides a modification of the Euler method and 
the Taylor--It\^{o} expansion for It\^{o} SDEs
with non-differentiable or singular drift and diffusion coefficients.
The convergence rates of the proposed modifications of the Milstein and Euler
methods are established.

\renewcommand{\baselinestretch}{1.3}

{\normalsize 

\linespread{1.0}

\tableofcontents

\linespread{1.0}

}





\section{Introduction}

Let $(\Omega,$ ${\rm F},$ ${\sf P})$ be a complete probability space and let 
$\{{\rm F}_t, t\in[0,T]\}$ be a nondecreasing right-continuous 
family of $\sigma$-algebras of ${\rm F}.$
Let ${\bf w}_t$ be an ${\rm F}_t$-measurable for all $t\in[0,T]$ standard vector Wiener process with
independent components ${\bf w}_t^{(i)},$ $i=1,\ldots,m.$
Consider an It\^{o} stochastic differential equation (SDE) in the integral form
\begin{equation}
\label{1.5.2}
{\bf x}_t={\bf x}_0+\int\limits_0^t {\bf a}({\bf x}_{s},s)ds+
\sum\limits_{j=1}^m
\int\limits_0^t B_j({\bf x}_{s},s)d{\bf w}_{s}^{(j)},\ \ \
{\bf x}_0={\bf x}(0,\omega),\ \ \ \omega\in\Omega,
\end{equation}
where $t\in [0, T],$ ${\bf x}(t,\omega)\stackrel{\sf def}{=}{\bf x}_t: 
[0,T]\times \Omega\to {\bf R}^n$ is a strong solution of 
the It\^{o} SDE (\ref{1.5.2}),
the nonrandom functions ${\bf a}({\bf x},t),
B_j({\bf x},t)$ $(j=1,\ldots,m):$ 
${\bf R}^n\times[0, T]\to{\bf R}^n$
guarantee the existence and uniqueness up to stochastic equivalence 
of a strong solution
of (\ref{1.5.2}) \cite{1982}, the second integral on 
the right-hand side of (\ref{1.5.2}) is 
the It\^{o} stochastic integral, 
the random vector variable ${\bf x}_{0}\in {\bf R}^n$ is ${\rm F}_0$-measurable, 
independent of ${\bf w}_t-{\bf w}_0$ $(t>0)$ and satisfies the condition
${\sf M}\bigl\{\left\Vert{\bf x}_0\right\Vert^2\bigr\}<\infty$ 
(where here and below $\left\Vert\cdot\right\Vert$ 
denotes the Euclidean norm of a vector and ${\sf M}$ is an expectation operator).

It\^{o} SDEs have found widespread application as mathematical models of dynamic 
systems in various scientific fields \cite{2}-\cite{15}.
This explains the need to develop numerical methods for It\^{o} SDEs.
It should be noted that even to simulate the exact solution of It\^{o} SDE 
(such solutions, as a rule, are known only in rare cases), one often has 
to apply specific numerical methods for It\^{o} SDEs.

The classical theory of numerical integration of It\^{o} SDEs assumes that 
their drift and diffusion coefficients are sufficiently smooth functions
\cite{2}, \cite{3}, \cite{6}, \cite{15}.
Among classical numerical methods for It\^{o} SDEs, the so-called Milstein 
method \cite{15a}, \cite{15b} is of great interest, as it converges strongly with order 1.0 
under standard conditions \cite{2}, \cite{15a}, \cite{15b}.
This is due to the fact that the Euler--Maruyama method (the simplest numerical 
method for It\^{o} SDEs) converges strongly with order 0.5, which turns out to be 
insufficient for solving a number of practical problems.

There are a number of papers dedicated to the Milstein method that relax 
the classical conditions 
that ensure its convergence \cite{16}-\cite{22a}.

The It\^{o} SDEs with singular coefficients and numerical methods for such equations 
also attract significant attention of researchers
\cite{22b}-\cite{22ffff} (see also the literature review in \cite{22ffff}).

We specifically note the recent work \cite{22ffff} (2026), which, in particular,  investigates 
a polygonal-type Euler--Maruyama method for the numerical solution of It\^{o} SDEs with 
a singular drift coefficient that has a singularity at zero. 
More precisely, in the numerical experiments in \cite{22ffff}, the drift coefficient has 
a singular multiplier (depending only on 
the time variable) that tends to infinity at the initial time moment $t=0$.
The numerical method considered in \cite{22ffff} uses integral smoothing of the 
singularity instead of evaluating it at grid points, thereby avoiding division by zero.

Note that in our article, not only the drift coefficient is singular, but also 
the diffusion coefficients, which, just like in \cite{22ffff}, are integrally smoothed in 
the process of the numerical solution. At the same time, we propose a 
new methodology allowing, with the help of the selective application
of the It\^{o} formula, to obtain 
numerical methods of higher orders of strong convergence for It\^{o} SDEs with 
singular coefficients than the numerical method from \cite{22ffff}.

This article is devoted to a 
new modification of the Milstein method for non-commutative systems 
of It\^{o} SDEs (the case of a multidimensional Wiener process) with 
drift and diffusion coefficients having differentiable
and non-differentiable or 
singular power-law time-dependent components.
Namely, the above coefficients are products of a differentiable component (depending on the 
state variable and the time variable) and a non-differentiable 
or singular power-law one (depending only on the time variable).
By iteratively applying It\^{o}'s formula to the mentioned 
differentiable components, a new 
modification of the Milstein method is derived in the paper. 

A key feature 
of the above modification of the Milstein method is the presence of time-dependent  
weight functions in the iterated It\^{o} stochastic integrals of the second and first
multiplicity (the specified weight functions are identically equal to unity in the 
classical Milstein method \cite{2}, \cite{15a}, \cite{15b}). 
Thus, the problem of the non-smoothness or singularity of the It\^{o} SDE coefficients has shifted 
to the weight functions of the iterated It\^{o} stochastic integrals.
To approximate these iterated It\^{o} stochastic integrals (with respect to the 
components of a multidimensional Wiener process) involving time-dependent weight 
functions, we use the method of generalized multiple Fourier series, introduced by 
the first author in 2006 \cite{23}-\cite{54} (also see \cite{15}, \cite{55}).

The article also provides a modification of the Euler method and 
the Taylor--It\^{o} expansion for It\^{o} SDEs
with non-differentiable or singular drift and diffusion coefficients.
The convergence rates of the proposed modifications of the Milstein and Euler
methods are established.

\section{The Milstein Method}

In this section, we consider the explicit Milstein method 
for the It\^{o} SDE (\ref{1.5.2}) \cite{2}, \cite{15a}, \cite{15b}.

Let ${\bf x}_t,\ t\in[0,T]$ be a strong solution 
of the It\^{o} SDE (\ref{1.5.2}). 
Consider a partition $\{\tau_p\}_{p=0}^N$ 
of the interval $[0, T]$ such that
\begin{equation}
\label{ccc2}
0=\tau_0<\tau_1<\ldots <\tau_N=T,\ \ \
\Delta_N=
\hbox{\vtop{\offinterlineskip\halign{
\hfil#\hfil\cr
{\rm max}\cr
$\stackrel{}{{}_{0\le j\le N-1}}$\cr
}} }\left|\tau_{j+1}-\tau_j\right|.
\end{equation}

Let ${\bf y}_{\tau_j},$\ $j=0, 1,\ldots,N$ denote a discrete 
approximation of the stochastic process ${\bf x}_t,\ t\in[0,T]$
corresponding to the maximum discretization step $\Delta_N$
(instead of ${\bf y}_{\tau_j}$, it would be more correct to write 
${\bf y}_{\tau_j}^{\Delta_N},$ however, we omit this superscript 
to simplify the notation).

Let us define the strong convergence of a discrete approximation
and the rate of strong convergence.

{\bf Definition 1}\ \cite{2}. {\it We say that 
the discrete approximation {\rm (}numerical me\-thod{\rm )}
${\bf y}_{\tau_j},$ $j=0, 1,\ldots,N,$
corresponding to the maximum discretization step $\Delta_N,$
converges strongly to the stochastic process ${\bf x}_t,$ $t\in[0,T]$
at time $T$ if
$\lim\limits_{\Delta_N\downarrow 0}
{\sf M}\{\left\Vert{\bf x}_{T}-{\bf y}_{T}\right\Vert\}=0.$
}

{\bf Definition 2}\ \cite{2}.
{\it We say that the discrete approximation 
{\rm (}numerical me\-thod{\rm )}
${\bf y}_{\tau_j},$ $j=0, 1,\ldots,N,$
corresponding to the maximum discretization step $\Delta_N,$
converges strongly with order $\gamma>0$ at time
$T$ to the stochastic process ${\bf x}_t,$ $t\in[0,T]$
if there exist a constant $C\in(0,\infty),$ which is 
independent of $\Delta_N,$ and a number $\delta>0$ such that
\begin{equation}
\label{u1.1}
\varepsilon\stackrel{\sf def}{=}{\sf M}\{\left\Vert{\bf x}_{T}-{\bf y}_{T}\right\Vert\}\le
C(\Delta_N)^{\gamma}
\end{equation}
for all $\Delta_N\in(0, \delta).$}

Milstein method 
for the It\^{o} SDE (\ref{1.5.2}) can be 
written in the following form
\cite{2}, \cite{15a}, \cite{15b}
$$
{\bf y}_{p+1}={\bf y}_p+\sum_{j=1}^m B_{j}({\bf y}_p,\tau_p)
I_{\tau_{p+1},\tau_p}^{(j)}+(\tau_{p+1}-\tau_p){\bf a}({\bf y}_p,\tau_p)+
$$
\begin{equation}
\label{4.10}
+\sum_{l,j=1}^{m}\sum^{n}_{k=1} 
B^{(kl)}({\bf y}_p,\tau_p)
{\partial B_{j} \over \partial {\bf x} ^ {(k)}}({\bf y}_p,\tau_p)
I_{\tau_{p+1},\tau_p}^{(l j)},
\end{equation}
where $\tau_p=p\Delta$ $(p=0,1,\ldots,N),$
$\Delta=T /N$ $(N>1)$ is a constant (for simplicity) step of integration,
${\bf y}_{\tau_p}\stackrel{\sf def}{=}{\bf y}_p,$ 
$B^{(kl)}({\bf x},t)$ is the $kl$th element of the matrix function $B({\bf x},t):$
${\bf R}^n \times [0, T]\to {\bf R}^{n\times m}$ with columns $B_j({\bf x},t)$
$(j=1,\ldots,m),$ 
and 
$I_{\tau_{p+1},\tau_p}^{(j)},$ 
$I_{\tau_{p+1},\tau_p}^{(l j)}$ are iterated It\^{o} stochastic integrals of the form
\begin{equation}
\label{int1}
I_{\tau_{p+1},\tau_p}^{(j)}=\int\limits_{\tau_p}^{\tau_{p+1}}d{\bf w}_{s}^{(j)},\ \ \
I_{\tau_{p+1},\tau_p}^{(l j)}=\int\limits_{\tau_p}^{\tau_{p+1}}\int\limits_{\tau_p}^{s}
d{\bf w}_{\tau}^{(l)}d{\bf w}_{s}^{(j)},\ \ \ l,j=1,\ldots,m.
\end{equation}

It is known that the Milstein method (\ref{4.10}) converges strongly
with order 1.0 under standard conditions \cite{2}, \cite{15a}, \cite{15b}.

Note that if the so-called commutativity condition of the form
\begin{equation}
\label{int2}
\sum^{n}_{k=1} 
B^{(kl)}({\bf x},t)
{\partial B_{j} \over \partial {\bf x} ^ {(k)}}({\bf x},t)=
\sum^{n}_{k=1} 
B^{(kj)}({\bf x},t)
{\partial B_{l} \over \partial {\bf x} ^ {(k)}}({\bf x},t)\ \ \ (\hbox{for all}\ l,j=1,\ldots,m)
\end{equation}
is not satisfied, the iterated It\^{o} stochastic integrals (\ref{int1})
cannot be simulated numerically using simple formulas \cite{2}, \cite{15},
\cite{23}-\cite{26a}. This issue will be discussed in detail below.

Further, throughout this paper, we will assume that the commutativity 
condition (\ref{int2}) does not hold, unless otherwise specified.

\section{A New Modification of the Milstein Method for It\^{o} SDEs
with Drift and Diffusion Coefficients Having Non-Differentiable 
Time-Dependent Multipliers}

Consider the following special case of the It\^{o} SDE (\ref{1.5.2})
\begin{equation}
\label{1.5.2s}
{\bf x}_t={\bf x}_0+\int\limits_0^t \sigma_0(s){\bf a}({\bf x}_{s},s)ds+
\sum\limits_{j=1}^m
\int\limits_0^t \sigma_j(s) B_j({\bf x}_{s},s)d{\bf w}_{s}^{(j)},\ \ \
{\bf x}_0={\bf x}(0,\omega),
\end{equation}
where $t\in [0, T],$ $\sigma_0(t),$ $\sigma_j(t)$ $(j=1,\ldots,m):$ $[0, T]\to{\bf R}$ are measurable
functions that satisfy the following condition
\begin{equation}
\label{300fff}
\left|\sigma_0(t)\right|+\sum\limits_{j=1}^m \left|\sigma_j(t)\right|\le K
\end{equation}
for all $t\in [0, T],$ where the constant $K$ is independent of $t,$
i.e., the functions $\sigma_0(t),$ $\sigma_j(t)$ $(j=1,\ldots,m)$ are allowed to be non-differentiable.
The nonrandom functions ${\bf a}({\bf x},t),$
$B_j({\bf x},t)$ $(j=1,\ldots,m):$ 
${\bf R}^n\times[0, T]\to{\bf R}^n$
are measurable with respect to the collection of 
variables $({\bf x},t)$, satisfy the global Lipschitz condition, and the linear growth condition.
Furthermore, the random vector variable ${\bf x}_{0}\in {\bf R}^n$ is ${\rm F}_0$-measurable, 
independent of ${\bf w}_t-{\bf w}_0$ $(t>0)$ and satisfies the condition
${\sf M}\bigl\{\left\Vert{\bf x}_0\right\Vert^2\bigr\}<\infty.$
The other notations 
appearing in (\ref{1.5.2s}) are the same as in (\ref{1.5.2}). 

The above conditions 
guarantee the existence and uniqueness up to stochastic equivalence 
of a strong solution
of (\ref{1.5.2s}) \cite{1982}. 
Moreover, all moments ${\sf M}\bigl\{\left\Vert {\bf x}_t\right\Vert^{2r}\bigr\}$ $(r\ge 1)$
of the strong solution ${\bf x}_t$ to the It\^{o} SDE
(\ref{1.5.2s}) are finite and bounded by a continuous function of $t$ \cite{1982}.

Furthermore, let us assume that the functions ${\bf a}({\bf x},t)$ and
$B_j({\bf x},t)$ $(j=1,\ldots,m)$ are sufficiently smooth with respect 
to the variables ${\bf x}$ and $t.$ At the same time,
the functions 
$\sigma_0(t)$ and $\sigma_j(t)$ $(j=1,\ldots,m)$
are allowed to be non-differentiable and satisfy the 
condition (\ref{300fff}).

Let us rewrite the system of It\^{o} SDEs (\ref{1.5.2s})
in component-wise form
\begin{equation}
\label{1.5.2ss}
{\bf x}_t^{(i)}={\bf x}_0^{(i)}+\int\limits_0^t \sigma_0(s){\bf a}^{(i)}({\bf x}_{s},s)ds+
\sum\limits_{j=1}^m
\int\limits_0^t \sigma_j(s) B^{(ij)}({\bf x}_{s},s)d{\bf w}_{s}^{(j)},
\end{equation}
where ${\bf x}^{(i)}$ is the $i$th component of
${\bf x}$, ${\bf a}^{(i)}({\bf x},t)$ is the $i$th component
of ${\bf a}({\bf x},t)$, $B^{(ij)}({\bf x},t)$
is the $ij$th element of $B({\bf x}, t)$, and $i=1,\ldots,n.$

Using (\ref{1.5.2ss}), we obtain w.~p.~1
\begin{equation}
\label{1.5.2sss}
{\bf x}_t^{(i)}={\bf x}_{\tau}^{(i)}+\int\limits_{\tau}^t \sigma_0(s){\bf a}^{(i)}({\bf x}_{s},s)ds+
\sum\limits_{j=1}^m
\int\limits_{\tau}^t \sigma_j(s) B^{(ij)}({\bf x}_{s},s)d{\bf w}_{s}^{(j)},
\end{equation}
where $t>\tau$ and $i=1,\ldots,n.$

Applying It\^{o}'s formula, we have for $s>\tau$ w.~p.~1 
\begin{equation}
\label{5}
{\bf a}^{(i)}({\bf x}_{s},s)={\bf a}^{(i)}({\bf x}_{\tau},\tau)+
\int\limits_{\tau}^s L{\bf a}^{(i)}({\bf x}_{\theta},\theta)d\theta+
\sum\limits_{l=1}^m
\int\limits_{\tau}^s G_0^{(l)}{\bf a}^{(i)}({\bf x}_{\theta},\theta)d{\bf w}_{\theta}^{(l)},
\end{equation}
\begin{equation}
\label{6}
B^{(ij)}({\bf x}_{s},s)=B^{(ij)}({\bf x}_{\tau},\tau)+
\int\limits_{\tau}^s LB^{(ij)}({\bf x}_{\theta},\theta)d\theta+
\sum\limits_{l=1}^m
\int\limits_{\tau}^s G_0^{(l)}B^{(ij)}({\bf x}_{\theta},\theta)d{\bf w}_{\theta}^{(l)},
\end{equation}
where $i=1,\ldots,m,$ 
the differential operators $L$ and $G_0^{(l)}$ $(l=1,\ldots,m)$ have the following form
\begin{equation}
\label{operl}
L= {\partial \over \partial t}
+ \sigma_0(t) \sum^{n}_{l=1}a^{(l)}({\bf x},t) 
{\partial \over  \partial  {\bf  x}^{(l)}} +
{1\over 2} \sum^{m}_{j=1} \sigma_j^2(t)\sum^{n}_{l,k=1}
B^{(lj)}({\bf x}, t) B^ {(kj)}({\bf x}, t) 
{\partial^{2} \over\partial{\bf x}^{(l)}\partial{\bf x}^{(k)}},
\end{equation}
$$
G_0^{(l)}=\sigma_l(t)G^{(l)},\ \ \ 
G^{(l)}=\sum^{n}_{k=1} B^{(kl)}({\bf x}, t)
{\partial \over \partial {\bf x}^{(k)}},\ \ \
l=1,\ldots,m.
$$

Combining (\ref{1.5.2sss})--(\ref{6}), we get w.~p.~1
$$
{\bf x}_t^{(i)}={\bf x}_{\tau}^{(i)}+{\bf a}^{(i)}({\bf x}_{\tau},\tau)\int\limits_{\tau}^t \sigma_0(s)ds+
\sum\limits_{j=1}^m
B^{(ij)}({\bf x}_{\tau},\tau)\int\limits_{\tau}^t \sigma_j(s)d{\bf w}_{s}^{(j)}+
$$
\begin{equation}
\label{6s}
+\sum\limits_{l,j=1}^m
\int\limits_{\tau}^t \sigma_j(s)\int\limits_{\tau}^s
\sigma_l(\theta)G^{(l)}B^{(ij)}({\bf x}_{\theta},\theta)d{\bf w}_{\theta}^{(l)}d{\bf w}_{s}^{(j)}+
R_{t,\tau}^{(i)},
\end{equation}
where
$$
R_{t,\tau}^{(i)}=
\int\limits_{\tau}^t \sigma_0(s)\int\limits_{\tau}^s
L{\bf a}^{(i)}({\bf x}_{\theta},\theta)d\theta ds+
$$
$$
+
\sum\limits_{l=1}^m
\int\limits_{\tau}^t \sigma_0(s)\int\limits_{\tau}^s
G_0^{(l)}{\bf a}^{(i)}({\bf x}_{\theta},\theta)d{\bf w}_{\theta}^{(l)}ds+
$$
\begin{equation}
\label{6s111}
+
\sum\limits_{j=1}^m
\int\limits_{\tau}^t \sigma_j(s)\int\limits_{\tau}^s
LB^{(ij)}({\bf x}_{\theta},\theta)d\theta d{\bf w}_{s}^{(j)}.
\end{equation}

Omitting the remainder term $R_{t,\tau}^{(i)}$ and 
replacing $G^{(l)}B^{(ij)}({\bf x}_{\theta},\theta)$ with
$G^{(l)}B^{(ij)}({\bf x}_{\tau},\tau)$ in (\ref{6s}),
we obtain the following new modification of the Milstein method 
in component-wise form ($\tau=\tau_p,$ $t=\tau_{p+1}$)
$$
{\bf y}_{p+1}^{(i)}={\bf y}_{p}^{(i)}+
{\bf a}^{(i)}({\bf y}_{p},\tau_p) I^{\sigma_0}_{\tau_{p+1},\tau_p}+
\sum\limits_{j=1}^m
B^{(ij)}({\bf y}_{p},\tau_p)I^{\sigma_j(j)}_{\tau_{p+1},\tau_p}
+
$$
\begin{equation}
\label{10}
+\sum\limits_{l,j=1}^m 
\sum^{n}_{k=1} 
B^{(kl)}({\bf y}_p,\tau_p)
{\partial B^{(ij)} \over \partial {\bf x} ^ {(k)}}({\bf y}_p,\tau_p)
I^{\sigma_l\sigma_j(lj)}_{\tau_{p+1},\tau_p},
\end{equation}
where
\begin{equation}
\label{11a}
I^{\sigma_0}_{\tau_{p+1},\tau_p}=\int\limits_{\tau_p}^{\tau_{p+1}} \sigma_0(s)ds,
\end{equation}
\begin{equation}
\label{11b}
I^{\sigma_j(j)}_{\tau_{p+1},\tau_p}=
\int\limits_{\tau_p}^{\tau_{p+1}} \sigma_j(s)d{\bf w}_{s}^{(j)},
\end{equation}
\begin{equation}
\label{11c}
I^{\sigma_l\sigma_j(lj)}_{\tau_{p+1},\tau_p}
=
\int\limits_{\tau_p}^{\tau_{p+1}} \sigma_j(s)\int\limits_{\tau_p}^s
\sigma_l(\theta)d{\bf w}_{\theta}^{(l)}d{\bf w}_{s}^{(j)},
\end{equation}
$\tau_p=p\Delta$ $(p=0,1,\ldots,N),$
$\Delta=T /N$ $(N>1)$ is a constant (for simplicity) step of integration,
${\bf y}_{\tau_p}\stackrel{\sf def}{=}{\bf y}_p,$ 
$B^{(kl)}({\bf x},t)$ is the $kl$th element of the matrix function $B({\bf x},t):$
${\bf R}^n \times [0, T]\to {\bf R}^{n\times m}$ with columns $B_j({\bf x},t)$
$(j=1,\ldots,m).$

The numerical scheme (\ref{10}) can be written in vector form as follows
$$
{\bf y}_{p+1}={\bf y}_{p}+
{\bf a}({\bf y}_{p},\tau_p)I^{\sigma_0}_{\tau_{p+1},\tau_p}+
\sum\limits_{j=1}^m
B_j({\bf y}_{p},\tau_p)I^{\sigma_j(j)}_{\tau_{p+1},\tau_p}
+
$$
\begin{equation}
\label{11}
+\sum\limits_{l,j=1}^m 
\sum^{n}_{k=1} 
B^{(kl)}({\bf y}_p,\tau_p)
{\partial B_j \over \partial {\bf x} ^ {(k)}}({\bf y}_p,\tau_p)
I^{\sigma_l\sigma_j(lj)}_{\tau_{p+1},\tau_p}.
\end{equation}

{\bf Remark~1.}\ {\it Obviously$,$ under the standard conditions {\rm \cite{2}}
on the functions ${\bf a}({\bf x},t),$ 
$B_j({\bf x},t),$ 
$L{\bf a}^{(i)}({\bf x},t),$ $G^{(l)}{\bf a}^{(i)}({\bf x},t),$
$LB^{(ij)}({\bf x}_t,t),$ $G^{(l)}B^{(ij)}({\bf x},t)$ 
as well as under the condition {\rm (\ref{300fff})}
on the functions $\sigma_0(t),$ $\sigma_j(t)$ $(j=1,\ldots,m)$, 
it can be shown in a standard way {\rm \cite{2}} that the 
numerical methods {\rm (\ref{10})} and {\rm (\ref{11})} have a strong order {\rm 1.0} of convergence.}

Assume that the following commutativity condition 
\begin{equation}
\label{2000}
\sum^{n}_{k=1} 
B^{(kl)}({\bf x},t)
{\partial B_{j} \over \partial {\bf x} ^ {(k)}}({\bf x},t)=
\sum^{n}_{k=1} 
B^{(kj)}({\bf x},t)
{\partial B_{l} \over \partial {\bf x} ^ {(k)}}({\bf x},t)\ \ \ (\hbox{for all}\ l,j=1,\ldots,m)
\end{equation}
holds. Then the numerical method (\ref{11}) can be written in the form
$$
{\bf y}_{p+1}={\bf y}_{p}+
{\bf a}({\bf y}_{p},\tau_p)I^{\sigma_0}_{\tau_{p+1},\tau_p}+
\sum\limits_{j=1}^m
B_j({\bf y}_{p},\tau_p)I^{\sigma_j(j)}_{\tau_{p+1},\tau_p}
+
$$
$$
+\frac{1}{2}\sum\limits_{l,j=1}^m 
\sum^{n}_{k=1} 
B^{(kl)}({\bf y}_p,\tau_p)
{\partial B_j \over \partial {\bf x} ^ {(k)}}({\bf y}_p,\tau_p)\times
$$
\begin{equation}
\label{200}
\times
\left( I^{\sigma_j(j)}_{\tau_{p+1},\tau_p}\cdot I^{\sigma_l(l)}_{\tau_{p+1},\tau_p}-
{\bf 1}_{\{l=j\}}
\int\limits_{\tau_p}^{\tau_{p+1}}\sigma_l^2(s)ds\right),
\end{equation}

\noindent
where ${\bf 1}_A$ denotes the indicator of the set $A.$

Note that the classical Milstein method for the It\^{o} SDE
(\ref{1.5.2s}) takes the form
$$
{\bf y}_{p+1}={\bf y}_{p}+
\sigma_0(\tau_p){\bf a}({\bf y}_{p},\tau_p)(\tau_{p+1}-\tau_p)+
\sum\limits_{j=1}^m
\sigma_j(\tau_p) B_j({\bf y}_{p},\tau_p)I^{(j)}_{\tau_{p+1},\tau_p}
+
$$
\begin{equation}
\label{11vv}
+\sum\limits_{l,j=1}^m \sigma_l(\tau_p)\sigma_j(\tau_p)
\sum^{n}_{k=1} 
B^{(kl)}({\bf y}_p,\tau_p)
{\partial B_j \over \partial {\bf x} ^ {(k)}}({\bf y}_p,\tau_p)
I^{(lj)}_{\tau_{p+1},\tau_p},
\end{equation}
where the iterated It\^{o} stochastic integrals 
$I^{(j)}_{\tau_{p+1},\tau_p},$ $I^{(lj)}_{\tau_{p+1},\tau_p}$ are defined by 
(\ref{int1}).

Under the commutativity condition (\ref{2000}), the Milstein method
(\ref{11vv}) can be written as
$$
{\bf y}_{p+1}={\bf y}_{p}+
\sigma_0(\tau_p){\bf a}({\bf y}_{p},\tau_p)(\tau_{p+1}-\tau_p)+
\sum\limits_{j=1}^m
\sigma_j(\tau_p)B_j({\bf y}_{p},\tau_p)I^{(j)}_{\tau_{p+1},\tau_p}
+
$$
$$
+\frac{1}{2}\sum\limits_{l,j=1}^m \sigma_l(\tau_p)\sigma_j(\tau_p)
\sum^{n}_{k=1} 
B^{(kl)}({\bf y}_p,\tau_p)
{\partial B_j \over \partial {\bf x} ^ {(k)}}({\bf y}_p,\tau_p)\times
$$
\begin{equation}
\label{11vvv}
\times
\biggl( I^{(j)}_{\tau_{p+1},\tau_p}\cdot I^{(l)}_{\tau_{p+1},\tau_p}-
{\bf 1}_{\{l=j\}}
(\tau_{p+1}-\tau_p)\biggr).
\end{equation}

\section{Application of the New Modification (\ref{11}) of the Milstein Method 
to Numerical Integration of It\^{o} SDEs
with Drift and Diffusion Coefficients Having 
Singular Power-Law Time-Dependent Multipliers}

First, let us consider the following non-commutative linear It\^{o} SDE
\begin{equation}
\label{1.5.2svv}
x_t=x_0+\int\limits_0^t \sigma_0(s)\left(a(s)x_{s}+c(s)\right)ds+
\sum\limits_{j=1}^m
\int\limits_0^t \sigma_j(s) \left(b_j(s)x_{s}+f_j(s)\right)d{\bf w}_{s}^{(j)},
\end{equation}

\vspace{1mm}
\noindent
where $x_0=x(0,\omega),$
$t\in [0, T],$ $\sigma_0(t)=t^{\alpha},$
$\sigma_j(t)=t^{\beta_j},$ $\alpha\in (-1,0)$, $\beta_j\in (-1/2,0)$ $(j=1,\ldots,m)$,
the nonrandom functions $a(s),$ $c(s),$ $b_j(s),$ $f_j(s):$ $[0, T]\to{\bf R}$ $(j=1,\ldots,m)$
are continuous, and their first-order derivatives are uniformly bounded.

From the above conditions on the functions $\sigma_0(t),$ $\sigma_j(t)$ $(j=1,\ldots,m)$, we 
have
\begin{equation}
\label{300}
\sigma_0(t)\in L_1([0, T],{\bf R}),\ \ \ \sigma_j(t)\in L_2([0, T],{\bf R}),
\end{equation}

\vspace{1mm}
\noindent
i.e.,
$$
\int\limits_0^T \left|\sigma_0(s)\right|ds+\sum\limits_{j=1}^m \int\limits_0^T \sigma_j^2(s)ds<\infty.
$$

\vspace{1mm}

According to Theorem~1.1 from \cite{Jan}, under the above conditions, the It\^{o} SDE (\ref{1.5.2svv})
has a unique strong solution.

We can repeat the derivation of the numerical method (\ref{10}) for 
the It\^{o} SDE (\ref{1.5.2svv}). Then, all integrals on the right-hand 
sides of (\ref{5}), (\ref{6}), (\ref{6s}), (\ref{6s111}) 
will exist in the mean-square sense by virtue of condition (\ref{300}),
the linearity
with respect to $x_t$ of the functions on the right-hand sides of these equalities,
and povided that
${\sf M}\left\{x_t^2\right\}\le K,$ where $x_t$ is a strong solution 
to the It\^{o} SDE (\ref{1.5.2svv}) and $K$ is a constant that does not depend on $t$
(see the part of this section concerning the It\^{o} SDE (\ref{1.5.2svd})
for details).

Let us verify that the condition ${\sf M}\left\{x_t^2\right\}\le K$ 
holds for the It\^{o} SDE (\ref{1.5.2svv}).
By applying It\^{o}'s formula, it is easy to obtain the following system
of differential equations for ${\sf M}\left\{x_t\right\}\stackrel{\sf def}{=}m_t$
and ${\sf M}\left\{x_t^2\right\}\stackrel{\sf def}{=}\Gamma_t$
\cite{2}

\newpage
\noindent
\begin{equation}
\label{syst1}
\left\{\begin{matrix}
\displaystyle{\frac{dm_t}{dt}}=\sigma_0(t)\left(a(t)m_t+c(t)\right)\cr \cr
\displaystyle{\frac{d\Gamma_t}{dt}}=A(t)\Gamma_t+B(t)
\end{matrix}\right.,
\end{equation}

\noindent
where
$$
A(t)=2\sigma_0(t)a(t)+\sum\limits_{j=1}^m \sigma_j^2(t)b_j^2(t),
$$
$$
B(t)=2\left(\sigma_0(t)c(t)+\sum\limits_{j=1}^m \sigma_j^2(t)b_j(t)f_j(t)\right)m_t+
\sum\limits_{j=1}^m \sigma_j^2(t)f_j^2(t).
$$

The solution to the system (\ref{syst1}) has the following form
\begin{equation}
\label{syst2}
m_t=\Phi_0(t)\left(m_0+\int\limits_0^t \frac{\sigma_0(s)c(s)}{\Phi_0(s)}ds\right),
\end{equation}
\begin{equation}
\label{syst3}
\Gamma_t=\Phi_1(t)\left(\Gamma_0+\int\limits_0^t \frac{B(s)}{\Phi_1(s)}ds\right),
\end{equation}
where
\begin{equation}
\label{syst4}
\Phi_0(t)=\exp \left(\int\limits_0^t \sigma_0(s)a(s)ds\right), \ \ \ 
\Phi_1(t)=\exp \left(\int\limits_0^t A(s)ds\right).
\end{equation}

From (\ref{syst2})--(\ref{syst4}), the properties of the functions 
$\sigma_0(t),$
$\sigma_j(t),$
$a(s),$ $c(s),$ $b_j(s),$ $f_j(s)$ $(j=1,\ldots,m),$
and the 
properties of the Lebesgue integral, we conclude that 
the function (\ref{syst3}) is continuous on $[0, T]$ and therefore bounded, i.e., 
${\sf M}\left\{x_t^2\right\}\le K,$ where $K$ is a constant that does not depend on $t.$

The convergence rate of the numerical method (\ref{10}) applied to the It\^{o} SDE (\ref{1.5.2svv}) 
and It\^{o} SDEs (\ref{1.5.2svd}), (\ref{1600t}) (see below)
will be considered in Sect.~12.

It is obvious that the problem of numerical integration of the It\^{o} SDE (\ref{1.5.2s}) or
(\ref{1.5.2svv})
with non-differentiable or singular time-dependent drift and diffusion coefficients (due to the functions
$\sigma_0(t)$ and $\sigma_j(t)$ $(j=1,\ldots,m)$) 
has shifted, in numerical methods (\ref{10}) and (\ref{11}), to the weight functions 
$\sigma_0(t),$ $\sigma_j(t)$ $(j=1,\ldots,m)$ appearing 
in the iterated It\^{o} stochastic integrals (\ref{11b}), (\ref{11c}) and the 
integral (\ref{11a}).

Even in the case when $\sigma_j(t)\equiv 1$ $(j=1,\ldots,m)$, 
numerical simulation of the iterated It\^{o} stochastic 
integrals (\ref{11b}), (\ref{11c}) is a difficult problem \cite{2}, \cite{15a}, \cite{15b}.
Nevertheless, in the cases when $\sigma_j(t)\not\equiv 1$ or $\sigma_j(t)\equiv 1$ $(j=1,\ldots,m)$, 
the method of generalized multiple Fourier series 
(proposed by the first author in 2006 \cite{23}-\cite{26a} (also see \cite{15}, \cite{27})) 
can be successfully applied to approximate 
the iterated It\^{o} stochastic integrals (\ref{11b}), (\ref{11c}). 
The next section will be devoted to this issue.

Due to the singularity of the functions 
$\sigma_0(t)=t^{\alpha},$
$\sigma_j(t)=t^{\beta_j}$ $(\alpha\in (-1,0)$, $\beta_j\in (-1/2,0)$ $(j=1,\ldots,m))$
at $t=0,$ it is reasonable 
to expect that the classical Milstein method applied to the It\^{o} SDE (\ref{1.5.2svv}) 
(see (\ref{11vv})) will diverge (it will be simply impossible to execute 
the first step of this method due to division by zero if the integration 
interval is $[0, T]$), whereas the new modification of the Milstein 
method (\ref{11}) will still converge. 
However, one can set $\tau_0=\varepsilon>0$ $(\varepsilon<\Delta)$
when performing the first step by the classical Milstein method (\ref{11vv}), 
or use other approaches that help to perform the first integration step
by the classical Milstein method (\ref{11vv}).
We will verify these hypotheses using numerical experiments in Sect. 13 and 14.

Now, let us consider a nonlinear system of It\^{o} SDEs of the form
\begin{equation}
\label{1.5.2svd}
{\bf x}_t={\bf x}_0+\int\limits_0^t \sigma_0(s){\bf a}({\bf x}_{s})ds+
\sum\limits_{j=1}^m
\int\limits_0^t \sigma_j(s) B_j({\bf x}_{s})d{\bf w}_{s}^{(j)},\ \ \
{\bf x}_0={\bf x}(0,\omega),
\end{equation}
where $t\in [0, T],$ $\sigma_0(t)=t^{\alpha},$ $\sigma_j(t)=t^{\beta_j},$ 
$\alpha\in (-\left|h_0\right|, 0),$ $\beta_j\in (-\left|h_1\right|, 0)$
(the real constants $h_0,$ $h_1$ will be defined below), $j=1,\ldots,m,$
the nonrandom functions ${\bf a}({\bf x}),$
$B_j({\bf x})$ $(j=1,\ldots,m):$ 
${\bf R}^n\to{\bf R}^n$
are measurable, satisfy the global Lipschitz condition, and the linear growth condition
(further, 
${\bf x}^{(i)}$ is the $i$th component of
${\bf x}$, ${\bf a}^{(i)}({\bf x})$ is the $i$th component
of ${\bf a}({\bf x})$, and
$B^{(ij)}({\bf x})$ is the $ij$th element of the matrix $B({\bf x})$
with columns $B_j({\bf x})$).
The other notations 
appearing in (\ref{1.5.2svd}) are the same as in (\ref{1.5.2}).

Note that if $\left|h_0\right|\le 1$ and $\left|h_1\right|\le 1/2,$
then (under the above conditions) according to Theorem 1.1 \cite{Jan}, the It\^{o} SDE
(\ref{1.5.2svd}) has a unique strong solution.

We can repeat the derivation of the numerical method (\ref{10}) for 
the It\^{o} SDE (\ref{1.5.2svd}). However, all the integrals on the right-hand 
sides of (\ref{5}), (\ref{6}), (\ref{6s}), (\ref{6s111})
must exist in the mean-square sense, i.e., integrals of the second moments
of the integrand processes in (\ref{5}), (\ref{6}), (\ref{6s}), (\ref{6s111})  must exist.

On the other hand, depending on the growth conditions for
the partial derivatives of orders 1 and 2 for the functions 
${\bf a}^{(i)}({\bf x}),$ $B^{(ij)}({\bf x})$ $(i=1,\ldots,n,\ j=1,\ldots,m)$, 
ensuring the existence of the integrals 
on the right-hand sides of (\ref{5}), (\ref{6}), (\ref{6s}), (\ref{6s111}) 
in the mean-square sense will 
require the existence of certain moments ${\sf M}\left\{\left\Vert{\bf x}_t\right\Vert^r\right\}$ 
of the stochastic process ${\bf x}_t$ and the boundedness of these moments, for example, by a constant
that does not depend on $t$ or by a continuous function of $t$.

In \cite{Jan} (Theorem~1.7), sufficient conditions for the existence of
${\sf M}\left\{\left\Vert{\bf x}_t\right\Vert^r\right\}$ $(r>2)$ 
and the boundedness of these moments by a continuos function of $t$
are provided. However, the conditions 
of Theorem~1.7 \cite{Jan} automatically require that the condition 
\begin{equation}
\label{12321}
\sigma_0^r(t)\in L_1([0, T],{\bf R}),\ \ \ \sigma_j^{2r}(t)\in L_1([0, T],{\bf R})\ \ \
(j=1,\ldots,m)
\end{equation}
holds, where $r>2.$

First, we restrict ourselves to such growth conditions on the aforementioned 
partial derivatives of the functions ${\bf a}^{(i)}({\bf x}),$
$B^{(ij)}({\bf x}),$ $i=1,\ldots,n,$ $j=1,\ldots,m$ (see below), which lead to the requirement of
the existence of only the second moment ${\sf M}\bigl\{\left\Vert{\bf x}_t\right\Vert^2\bigr\}$ 
of ${\bf x}_t$. Thereby, we 
can take $r=2+\varepsilon$ in the condition (\ref{12321}) 
on the functions 
$\sigma_0(t),$ $\sigma_j(t),$ $j=1,\ldots,m$ ($\varepsilon$ is an arbitrary
small positive real number).

The aforementioned conditions on the partial derivatives, together with the linear
growth condition, take the form
\begin{equation}
\label{qw11}
\left|{\bf a}^{(i)}({\bf x})\right|+\left|B^{(ij)}({\bf x})\right|\le K\left(1+\left\Vert{\bf x}
\right\Vert\right),
\end{equation}
\begin{equation}
\label{qw12}
\left|\frac{\partial {\bf a}^{(i)}}{\partial {\bf x}^{(l)}}({\bf x})\right|+
\left|\frac{\partial B^{(ij)}}{\partial {\bf x}^{(l)}}({\bf x})\right|\le K,
\end{equation}
\begin{equation}
\label{qw13}
\left|\frac{\partial^2 {\bf a}^{(i)}}{\partial {\bf x}^{(l)}
\partial {\bf x}^{(k)}}({\bf x})\right|+
\left|\frac{\partial^2 B^{(ij)}}{\partial {\bf x}^{(l)}
\partial {\bf x}^{(k)}}({\bf x})\right|
\le \frac{K}{\left(1+\left\Vert{\bf x}\right\Vert\right)},
\end{equation}

\vspace{2mm}
\noindent
where the constant $K$ is independent of ${\bf x},$
and $i,l,k=1,\ldots,n,$\ $j=1,\ldots,m.$

As examples of functions satisfying the conditions (\ref{qw11})--(\ref{qw13}),
one can take $\bigl(1+\left\Vert{\bf x}\right\Vert^2\bigr)^{1/2},$
$\bigl(1+\left\Vert{\bf x}\right\Vert^4\bigr)^{1/4}.$

Let us consider the stochastic integral 

\begin{equation}
\label{rr100}
\int\limits_{\tau}^t \sigma_j(s)\int\limits_{\tau}^s
\sigma_0(\theta) \sum^{n}_{l=1}a^{(l)}({\bf x}_{\theta}) 
{\partial B^{(ij)} \over  \partial  {\bf  x}^{(l)}}({\bf x}_{\theta})
d\theta d{\bf w}_s^{(j)}\stackrel{\sf def}{=}I,
\end{equation}
which appears when considering the stochastic integral 
$$
\int\limits_{\tau}^t \sigma_j(s)\int\limits_{\tau}^s
LB^{(ij)}({\bf x}_{\theta})d\theta d{\bf w}_s^{(j)}
$$
(this integral is included in the formula (\ref{6s111})).

Using the stochastic integral (\ref{rr100})
as an example, 
we will show that 
the requirement for the existence of the integrals on 
the right-hand sides of (\ref{5}), (\ref{6}), (\ref{6s}), (\ref{6s111})
in the mean-square sense (under the conditions 
(\ref{qw11})--(\ref{qw13}))
does not 
require stronger conditions on the functions $\sigma_0(t),$ $\sigma_j(t)$ $(j=1,\ldots,m)$
than the condition (\ref{12321}) for $r=2+\varepsilon$
($\varepsilon$ is an arbitrary
small positive real number).

Using the It\^{o} isometry, the Cauchy--Bunyakovsky inequality, the Fubini Theorem,
conditions (\ref{qw11})--(\ref{qw13}), and condition (\ref{12321}) for $r=2+\varepsilon$, we obtain
$$
{\sf M}\left\{I^2\right\}=
\int\limits_{\tau}^t \sigma_j^2(s){\sf M}\left\{\left(\int\limits_{\tau}^s
\sigma_0(\theta) \sum^{n}_{l=1}a^{(l)}({\bf x}_{\theta}) 
{\partial B^{(ij)} \over  \partial  {\bf  x}^{(l)}}({\bf x}_{\theta})
d\theta\right)^2\right\} ds\le
$$
$$
\le
\int\limits_{\tau}^t \sigma_j^2(s){\sf M}\left\{\left(\int\limits_{\tau}^s
\left|\sigma_0(\theta)\right| \sum^{n}_{l=1}\left|a^{(l)}({\bf x}_{\theta}) 
{\partial B^{(ij)} \over  \partial  {\bf  x}^{(l)}}({\bf x}_{\theta})\right|
d\theta\right)^2\right\} ds\le
$$
$$
\le \int\limits_{\tau}^t \sigma_j^2(s)
\int\limits_{\tau}^s
\left|\sigma_0(\theta)\right|d\theta 
\int\limits_{\tau}^s \left|\sigma_0(\theta)\right| {\sf M}\left\{\left(
\sum^{n}_{l=1}\left|a^{(l)}({\bf x}_{\theta})
{\partial B^{(ij)}\over  \partial  {\bf  x}^{(l)}}({\bf x}_{\theta})  \right|
\right)^2\right\}d\theta  ds\le
$$
$$
\le C_1\int\limits_{\tau}^t \sigma_j^2(s)
\int\limits_{\tau}^s
\left|\sigma_0(\theta)\right|d\theta 
\int\limits_{\tau}^s \left|\sigma_0(\theta)\right| {\sf M}\left\{
\left(1+\left\Vert {\bf x}_{\theta}\right\Vert\right)^2\right\}d\theta  ds\le
$$
$$
\le C_2\int\limits_{\tau}^t \sigma_j^2(s)
\int\limits_{\tau}^s
\left|\sigma_0(\theta)\right|d\theta 
\int\limits_{\tau}^s \left|\sigma_0(\theta)\right| d\theta  ds\le
$$
\begin{equation}
\label{gg111}
\le C_2\int\limits_{0}^T \sigma_j^2(s) ds
\left(\int\limits_{0}^T
\left|\sigma_0(\theta)\right|d\theta\right)^2 <\infty,
\end{equation}
where $C_1,$ $C_2$ are constants, and
$\varepsilon$ is an arbitrary
small positive real number.

We see that the right-hand side of (\ref{gg111}) exists under the following condition
\begin{equation}
\label{12321s}
\sigma_0(t)\in L_1([0, T],{\bf R}),\ \ \ \sigma_j(t)\in L_2([0, T],{\bf R})\ \ \
(j=1,\ldots,m).
\end{equation}

However, in obtaining (\ref{gg111}), we used the stronger 
condition (\ref{12321}) for $r=2+\varepsilon$.
Combining the conditions (\ref{12321s}) and (\ref{12321}) ($r=2+\varepsilon$), we find
\begin{equation}
\label{12321dd}
\sigma_0^{2+\varepsilon}(t)\in L_1([0, T],{\bf R}),\ \ \ \sigma_j^{4+2\varepsilon}
(t)\in L_1([0, T],{\bf R})\ \ \
(j=1,\ldots,m),
\end{equation}
where $\varepsilon$ is an arbitrary
small positive real number.

Using the reasoning presented above for the stochastic integral (\ref{rr100}), we can obtain 
condition (\ref{12321dd}) for the remaining integrals included in the right-hand 
sides of (\ref{5}), (\ref{6}), (\ref{6s}), (\ref{6s111}) as well.

Since $\sigma_0(t)=t^{\alpha},$ $\sigma_j(t)=t^{\beta_j},$ 
$\alpha\in (-\left|h_0\right|, 0),$ $\beta_j\in (-\left|h_1\right|, 0),$
then from (\ref{12321dd}) it follows that $\alpha\in (-1/(2+\varepsilon), 0),$
$\beta_j\in (-1/(4+2\varepsilon), 0),$ i.e. $\left|h_0\right|=1/(2+\varepsilon),$
$\left|h_1\right|=1/(4+2\varepsilon).$ In the more general case, we obtain
$\sigma_0^{2+\varepsilon}(t)\in L_1([0, T],{\bf R}),$ 
$\sigma_j^{4+2\varepsilon}(t)\in L_1([0, T],{\bf R})$.

Suppose now that the functions ${\bf a}^{(i)}({\bf x}),$
$B^{(ij)}({\bf x})$ $(i=1,\ldots,n,\ j=1,\ldots,m)$
are bounded, and 
the partial derivatives of orders 1 and 2 for the functions 
${\bf a}^{(i)}({\bf x}),$ $B^{(ij)}({\bf x})$ 
are bounded.
As examples of functions satisfying to the above conditions,
one can take $\exp\bigl(-\left\Vert{\bf x}\right\Vert^2\bigr),$
$\bigl(1+\left\Vert{\bf x}\right\Vert^2\bigr)^{-1/2},$
$\left\Vert{\bf x}\right\Vert^2/\bigl(1+\left\Vert{\bf x}\right\Vert^2\bigr).$

In this case, in obtaining (\ref{gg111}), we use the estimate
\begin{equation}
\label{gg111ddd}
\left|a^{(l)}({\bf x}_{\theta})
{\partial B^{(ij)}\over  \partial  {\bf  x}^{(l)}}({\bf x}_{\theta})  \right|\le C,
\end{equation}
where the constant $C$ is independent of $\theta$. Therefore, in this case, the right-hand side
of (\ref{gg111}) will exist under the condition (\ref{12321s}) (without any
additional conditions).

Using the above reasoning for the stochastic integral (\ref{rr100}), we can obtain 
condition (\ref{12321s}) for the remaining integrals included in the right-hand 
sides of (\ref{5}), (\ref{6}), (\ref{6s}), (\ref{6s111}) as well.

Since $\sigma_0(t)=t^{\alpha},$ $\sigma_j(t)=t^{\beta_j},$ 
$\alpha\in (-\left|h_0\right|, 0),$ $\beta_j\in (-\left|h_1\right|, 0),$
then from (\ref{12321s}) it follows that $\alpha\in (-1, 0),$
$\beta_j\in (-1/2, 0),$ i.e. $\left|h_0\right|=1,$
$\left|h_1\right|=1/2$ for this case.
Generally, we get $\sigma_0(t)\in L_1([0, T],{\bf R}),$ 
$\sigma_j^{2}(t)\in L_1([0, T],{\bf R})$.

Obviously, other variants of the
growth conditions on the mentioned 
partial derivatives of the functions ${\bf a}^{(i)}({\bf x}),$
$B^{(ij)}({\bf x})$ $(i=1,\ldots,n,\ j=1,\ldots,m)$ are also possible.

For example, if these functions and the specified partial derivatives 
satisfy the linear growth condition (\ref{qw11}), then for the existence 
of stochastic integrals on 
the right-hand sides of (\ref{5}), (\ref{6}), (\ref{6s}), (\ref{6s111})
in the mean-square sense,
the existence of the  
moment ${\sf M}\bigl\{\left\Vert{\bf x}_t\right\Vert^{6}\bigr\}$ of the process 
${\bf x}_t$
will be required, as well as the boundedness of this moment, for example, 
by a continuous function of $t$.

Further, since $\sigma_0(t)=t^{\alpha},$ $\sigma_j(t)=t^{\beta_j},$ 
$\alpha\in (-\left|h_0\right|, 0),$ $\beta_j\in (-\left|h_1\right|, 0),$
then from (\ref{12321dd}) ($r=6)$ it follows that $\alpha\in (-1/6, 0),$
$\beta_j\in (-1/12, 0),$ i.e. $\left|h_0\right|=1/6,$
$\left|h_1\right|=1/12.$ More generally, we obtain
$\sigma_0^{6}(t)\in L_1([0, T],{\bf R}),$ $\sigma_j^{12}(t)\in L_1([0, T],{\bf R})$.
Thus, the singularities in this case
become significantly weaker than in previous cases.

It should be noted that conditions (\ref{nomer1})--(\ref{nomer4})
of Theorem~4 (see. Sect.~12) for the autonomous case 
(${\bf a}({\bf x},t)\equiv {\bf a}({\bf x}),$
$B_{j}({\bf x},t)\equiv B_{j}({\bf x})$) are satisfied
if conditions (\ref{qw11})--(\ref{qw13}) hold.

\section{Method of Generalized Multiple Fourier Series for the Approximation of
Iterated It\^{o} Stochastic Integrals}

Let us consider 
the following iterated It\^{o} 
stochastic integrals
\begin{equation}
\label{ito}
J[\psi^{(k)}]_{\bar T,t}^{(i_1\ldots i_k)}=\int\limits_t^{\bar T}\psi_k(t_k) \ldots \int\limits_t^{t_{2}}
\psi_1(t_1) d{\bf w}_{t_1}^{(i_1)}\ldots
d{\bf w}_{t_k}^{(i_k)},
\end{equation}
where $i_1,\ldots,i_k=0, 1,\ldots,m,$
$\psi_1(\tau),\ldots,\psi_k(\tau)\in L_2([t, \bar T],{\bf R}),$
${\bf w}_{\tau}^{(i)}$
($i=1,\ldots,m$) are independent standard Wiener processes, and
${\bf w}_{\tau}^{(0)}=\tau.$

Let us consider the approach to expansion of the iterated
It\^{o} stochastic integrals (\ref{ito}) 
\cite{23}-\cite{26a} (the so-called
method of generalized
multiple Fourier series). 
The method of generalized multiple Fourier series was introduced by 
the first author in 2006 \cite{23} and developed in \cite{23a}-\cite{26a}.

Let us
consider the unordered
set $\{1, 2, \ldots, k\}$ 
and separate it into two parts:
the first part consists of $r$ unordered 
pairs (sequence order of these pairs is also unimportant) and the 
second one consists of the 
remaining $k-2r$ numbers.
So, we have
\begin{equation}
\label{leto5007}
(\{
\underbrace{\{g_1, g_2\}, \ldots, 
\{g_{2r-1}, g_{2r}\}}_{\small{\hbox{part 1}}}
\},
\{\underbrace{q_1, \ldots, q_{k-2r}}_{\small{\hbox{part 2}}}
\}),
\end{equation}

\noindent
where 
$$
\{g_1, g_2, \ldots, 
g_{2r-1}, g_{2r}, q_1, \ldots, q_{k-2r}\}=\{1, 2, \ldots, k\},
$$

\noindent
braces   
mean an unordered 
set, and pa\-ren\-the\-ses mean an ordered set.

We will say that (\ref{leto5007}) is a partition 
and consider the sum with respect to all possible
partitions
\begin{equation}
\label{leto5008}
\sum_{\stackrel{(\{\{g_1, g_2\}, \ldots, 
\{g_{2r-1}, g_{2r}\}\}, \{q_1, \ldots, q_{k-2r}\})}
{{}_{\{g_1, g_2, \ldots, 
g_{2r-1}, g_{2r}, q_1, \ldots, q_{k-2r}\}=\{1, 2, \ldots, k\}}}}
a_{g_1 g_2, \ldots, 
g_{2r-1} g_{2r}, q_1 \ldots q_{k-2r}},
\end{equation}

\noindent
where $a_{g_1 g_2, \ldots, 
g_{2r-1} g_{2r}, q_1 \ldots q_{k-2r}}\in {\bf R}$.

Define the following function on the hypercube $[t, \bar T]^k$
\begin{equation}
\label{ppp}
K(t_1,\ldots,t_k)=
\left\{\begin{matrix}
\psi_1(t_1)\ldots \psi_k(t_k),\ &t_1<\ldots<t_k\cr\cr
0,\ &\hbox{\rm otherwise}
\end{matrix}
\right.,
\end{equation}

\noindent
where $\psi_1(\tau),\ldots,\psi_k(\tau)\in L_2([t, \bar T],{\bf R}),$
$t_1,\ldots,t_k\in [t, \bar T]$ $(k\ge 2)$ and 
$K(t_1)\equiv\psi_1(t_1)$ for $t_1\in[t, \bar T].$

Suppose that $\{\phi_j(x)\}_{j=0}^{\infty}$
is an arbitrary complete orthonormal system of functions in 
the space $L_2([t, \bar T],{\bf R})$.

Consider the Fourier coefficient corresponding to the kernel (\ref{ppp})
$$
C_{j_k\ldots j_1}=\int\limits_{[t,\bar T]^k}
K(t_1,\ldots,t_k)\prod_{l=1}^{k}\phi_{j_l}(t_l)dt_1\ldots dt_k=
$$
\begin{equation}
\label{ppppa}
=\int\limits_t^{\bar T}\psi_k(t_k)\phi_{j_k}(t_k)\ldots
\int\limits_t^{t_2}
\psi_1(t_1)\phi_{j_1}(t_1)
dt_1\ldots dt_k.
\end{equation}

{\bf Theorem 1}\ \cite{23a}-\cite{26a} (also see \cite{23} (2006)).\
{\it Suppose that
$\psi_1(\tau),\ldots,\psi_k(\tau)\in L_2([t, \bar T],{\bf R})$ and
$\{\phi_j(x)\}_{j=0}^{\infty}$ is an arbitrary complete orthonormal system  
of functions in the space $L_2([t,\bar T],{\bf R}).$
Then the following expansion
$$
J[\psi^{(k)}]_{\bar T,t}^{(i_1\ldots i_k)}=
\hbox{\vtop{\offinterlineskip\halign{
\hfil#\hfil\cr
{\rm l.i.m.}\cr
$\stackrel{}{{}_{p_1,\ldots,p_k\to \infty}}$\cr
}} }
\sum\limits_{j_1=0}^{p_1}\ldots
\sum\limits_{j_k=0}^{p_k}
C_{j_k\ldots j_1}\Biggl(
\prod_{l=1}^k\zeta_{j_l}^{(i_l)}+\sum\limits_{r=1}^{[k/2]}
(-1)^r \times
\Biggr.
$$
\begin{equation}
\label{razzar1}
\times
\sum_{\stackrel{(\{\{g_1, g_2\}, \ldots, 
\{g_{2r-1}, g_{2r}\}\}, \{q_1, \ldots, q_{k-2r}\})}
{{}_{\{g_1, g_2, \ldots, 
g_{2r-1}, g_{2r}, q_1, \ldots, q_{k-2r}\}=\{1, 2, \ldots, k\}}}}
\prod\limits_{s=1}^r
{\bf 1}_{\{i_{g_{{}_{2s-1}}}=~i_{g_{{}_{2s}}}\ne 0\}}
\Biggl.{\bf 1}_{\{j_{g_{{}_{2s-1}}}=~j_{g_{{}_{2s}}}\}}
\prod_{l=1}^{k-2r}\zeta_{j_{q_l}}^{(i_{q_l})}\Biggr)
\end{equation}

\noindent
con\-verg\-ing in the mean-square sense is valid$,$
where the iterated It\^{o} stochastic integral
$J[\psi^{(k)}]_{\bar T,t}^{(i_1\ldots i_k)}$ is defined by {\rm (\ref{ito}),}
$i_1,\ldots,i_k=0, 1,\ldots,m,$ $[x]$ is an integer part of a real number $x,$
${\bf 1}_A$ denotes the indicator of the set $A,$
$$
\prod\limits_{\emptyset}
\stackrel{\sf def}{=}1,\ \ \ 
\sum\limits_{\emptyset}
\stackrel{\sf def}{=}0,
$$
${\rm l.i.m.}$ is a limit in the mean-square sense$,$
\begin{equation}
\label{rr23}
\zeta_{j}^{(i)}=
\int\limits_t^{\bar T} \phi_{j}(s) d{\bf w}_s^{(i)}\ \ \ (i=0, 1,\ldots,m;\ j=0,1,2,\ldots)
\end{equation} 
are independent standard Gaussian random variables
for various
$i$ or $j$ {\rm(}in the case when $i\ne 0${\rm),}
$C_{j_k\ldots j_1}$ is the Fourier coefficient {\rm(\ref{ppppa})}.}

Consider particular cases of (\ref{razzar1}) for $k=1,2,3$
\begin{equation}
\label{a1}
J[\psi^{(1)}]_{\bar T,t}^{(i_1)}
=\hbox{\vtop{\offinterlineskip\halign{
\hfil#\hfil\cr
{\rm l.i.m.}\cr
$\stackrel{}{{}_{p_1\to \infty}}$\cr
}} }\sum_{j_1=0}^{p_1}
C_{j_1}\zeta_{j_1}^{(i_1)},
\end{equation}
\begin{equation}
\label{a2}
J[\psi^{(2)}]_{\bar T,t}^{(i_1 i_2)}
=\hbox{\vtop{\offinterlineskip\halign{
\hfil#\hfil\cr
{\rm l.i.m.}\cr
$\stackrel{}{{}_{p_1,p_2\to \infty}}$\cr
}} }\sum_{j_1=0}^{p_1}\sum_{j_2=0}^{p_2}
C_{j_2j_1}\Biggl(\zeta_{j_1}^{(i_1)}\zeta_{j_2}^{(i_2)}
-{\bf 1}_{\{i_1=i_2\ne 0\}}
{\bf 1}_{\{j_1=j_2\}}\Biggr),
\end{equation}

\vspace{3mm}
$$
J[\psi^{(3)}]_{\bar T,t}^{(i_1 i_2 i_3)}=
\hbox{\vtop{\offinterlineskip\halign{
\hfil#\hfil\cr
{\rm l.i.m.}\cr
$\stackrel{}{{}_{p_1,p_2,p_3\to \infty}}$\cr
}} }\sum_{j_1=0}^{p_1}\sum_{j_2=0}^{p_2}\sum_{j_3=0}^{p_3}
C_{j_3j_2j_1}\Biggl(
\zeta_{j_1}^{(i_1)}\zeta_{j_2}^{(i_2)}\zeta_{j_3}^{(i_3)}
-\Biggr.
$$
\begin{equation}
\label{a3}
\Biggl.
-
{\bf 1}_{\{i_1=i_2\ne 0\}}
{\bf 1}_{\{j_1=j_2\}}
\zeta_{j_3}^{(i_3)}
-{\bf 1}_{\{i_2=i_3\ne 0\}}
{\bf 1}_{\{j_2=j_3\}}
\zeta_{j_1}^{(i_1)}-
{\bf 1}_{\{i_1=i_3\ne 0\}}
{\bf 1}_{\{j_1=j_3\}}
\zeta_{j_2}^{(i_2)}\Biggr).
\end{equation}

\section{Formulas for Numerical Simulation of Iterated It\^{o}
Stochastic Integrals (\ref{11b}), (\ref{11c}) with Mean-Square Approximation Error Control}

Using (\ref{a1}) and (\ref{a2}), we obtain the following formulas
for approximating the iterated It\^{o} stochastic integrals
(\ref{11b}), (\ref{11c})
\begin{equation}
\label{100}
I^{\sigma_j(j)q_1}_{\tau_{p+1},\tau_p}=\sum_{k=0}^{q_1}
C_{k}^{\sigma_j}\zeta_{k}^{(j)},
\end{equation}
\begin{equation}
\label{101}
I^{\sigma_l\sigma_j(lj)q_2}_{\tau_{p+1},\tau_p}
=\sum_{k,r=0}^{q_2}
C_{r k}^{\sigma_l\sigma_j}\Biggl(\zeta_{k}^{(l)}\zeta_{r}^{(j)}
-{\bf 1}_{\{l=j\}}
{\bf 1}_{\{k=r\}}\Biggr),
\end{equation}
where $q_1, q_2\in {\bf N},$ $l,j=1,\ldots,m,$ 
\begin{equation}
\label{102}
C_{k}^{\sigma_j}=\int\limits_{\tau_p}^{\tau_{p+1}} \sigma_j(s)\phi_{k}(s)ds,
\end{equation}
\begin{equation}
\label{103}
C_{r k}^{\sigma_l\sigma_j}=\int\limits_{\tau_p}^{\tau_{p+1}}
\sigma_j(s) \phi_{r}(s)
\int\limits_{\tau_p}^{s}
\sigma_l(\theta)
\phi_{k}(\theta)d\theta ds,
\end{equation}
another notations are the same as in Theorem~1.

According to the results obtained in \cite{15} (Sect.~5.2, 5.9), 
\cite{25} (Sect.~1.2, 1.12), \cite{26a} (also see \cite{25a} (2018),
Sect.~6.2), we have
\begin{equation}
\label{104}
{\sf M}\left\{\left(I^{\sigma_j(j)}_{\tau_{p+1},\tau_p}-
I^{\sigma_j(j)q_1}_{\tau_{p+1},\tau_p}\right)^2\right\}=
\int\limits_{\tau_p}^{\tau_{p+1}}\sigma_j^2(s)ds-
\sum_{k=0}^{q_1}\left(C_{k}^{\sigma_j}\right)^2,
\end{equation}
$$
{\sf M}\left\{\left(
I^{\sigma_l\sigma_j(lj)}_{\tau_{p+1},\tau_p}-
I^{\sigma_l\sigma_j(lj)q_2}_{\tau_{p+1},\tau_p}
\right)^2\right\}=
$$
\begin{equation}
\label{105}
=
\int\limits_{\tau_p}^{\tau_{p+1}}\sigma_j^2(s)\int\limits_{\tau_p}^{s}
\sigma_l^2(\theta)d\theta ds-
\sum_{k,r=0}^{q_2}\left(C_{r k}^{\sigma_l\sigma_j}\right)^2
-{\bf 1}_{\{l=j\}}\sum_{k,r=0}^{q_2}C_{r k}^{\sigma_l\sigma_j}C_{k r}^{\sigma_l\sigma_j},
\end{equation}
or
\begin{equation}
\label{106}
{\sf M}\left\{\left(
I^{\sigma_l\sigma_j(lj)}_{\tau_{p+1},\tau_p}-
I^{\sigma_l\sigma_j(lj)q_2}_{\tau_{p+1},\tau_p}
\right)^2\right\}\le 2!
\left(\int\limits_{\tau_p}^{\tau_{p+1}}\sigma_j^2(s)\int\limits_{\tau_p}^{s}
\sigma_l^2(\theta)d\theta ds-
\sum_{k,r=0}^{q_2}\left(C_{r k}^{\sigma_l\sigma_j}\right)^2\right),
\end{equation}
where $q_1, q_2\in {\bf N},$ $l,j=1,\ldots,m.$ 

As we noted above, 
under the conditions on the functions 
${\bf a}({\bf x},t)$ and
$B_j({\bf x},t)$ $(j=1,\ldots,m)$
formulated at 
the beginning of Sect.~3 as well as under the condition
(\ref{300fff}) on the functions $\sigma_0(t),$ $\sigma_j(t)$ $(j=1,\ldots,m),$ 
it can be shown in a standard way \cite{2} that the                       
numerical methods (\ref{10}) and (\ref{11}) have a strong order 1.0 of convergence.

To implement the numerical methods (\ref{10}) or (\ref{11}), it is necessary
to replace the iterated It\^{o} stochastic integrals (\ref{11b}), (\ref{11c})
with their mean-square approximations (\ref{100}), (\ref{101}).
In this case, the numbers $q_1$ and $q_2$ must be chosen
from the following conditions (according to the classical theory
of the numerical integration of It\^{o} SDEs using the Milstein method \cite{2})
\begin{equation}
\label{107}
{\sf M}\left\{\left(I^{\sigma_j(j)}_{\tau_{p+1},\tau_p}-
I^{\sigma_j(j)q_1}_{\tau_{p+1},\tau_p}\right)^2\right\}\le C\Delta^{3},
\end{equation}
\begin{equation}
\label{108}
{\sf M}\left\{\left(
I^{\sigma_l\sigma_j(lj)}_{\tau_{p+1},\tau_p}-
I^{\sigma_l\sigma_j(lj)q_2}_{\tau_{p+1},\tau_p}
\right)^2\right\}\le C\Delta^{3},
\end{equation}

\vspace{1mm}
\noindent
where constant $C$ is independent of $\Delta$ 
($\Delta= T /N$ $(N>1)$ is a constant (for simplicity) step of integration,
$\tau_p=p\Delta$ $(p=0,1,\ldots,N)$).

Using formulas (\ref{100})--(\ref{108}), one can perform numerical simulation of the iterated It\^{o}
stochastic integrals (\ref{11b}), (\ref{11c}) with mean-square approximation error control.

Calculation of the Fourier coefficients (\ref{102}) and (\ref{103}) is a separate problem, 
and in the general case, these coefficients can only be calculated approximately.
However, if $\sigma_0(s)=s^{\alpha},$
$\sigma_j(s)=s^{\beta_j},$ $\alpha, \beta_j\in {\bf R}$ $(j=1,\ldots,m)$
and $\{\phi_j(s)\}_{j=0}^{\infty}$ is a complete orthonormal system of 
Legendre polynomials in the space $L_2([\tau_p,\tau_{p+1}],{\bf R}),$ 
then the Fourier coefficients (\ref{102}) and (\ref{103}) 
can be calculated using exact formulas.

Note that condition (\ref{300}) is satisfied, for example, when $\sigma_0(s)=s^{\alpha},$
$\sigma_j(s)=s^{\beta_j},$ $\alpha\in (-1,1]$, $\beta_j\in (-1/2,1]$ $(j=1,\ldots,m).$

\section{Exact Calculation of the Fourier--Legendre Co\-ef\-fi\-ci\-ents (\ref{102}), (\ref{103})}

Suppose that $\sigma_0(s)=s^{\alpha},$
$\sigma_j(s)=s^{\beta_j}$ ($\alpha\in (-1,1]$, $\beta_j\in (-1/2,1]$ $(j=1,\ldots,m)$)
and $\{\phi_j(s)\}_{j=0}^{\infty}$ is a complete orthonormal system of 
Legendre polynomials in the space $L_2([t, \bar T],{\bf R}).$

Then, the 
Fourier--Legendre coefficients 
(\ref{102}), (\ref{103}) calculated on the interval $[t, \bar T]$ take the form

$$
C_{k}^{\sigma_j}=\int\limits_{t}^{\bar T} s^{\beta_j}\phi_{k}(s)ds,
$$
$$
C_{r k}^{\sigma_l\sigma_j}=\int\limits_{t}^{\bar T}
s^{\beta_j} \phi_{r}(s)
\int\limits_{t}^{s}
\theta^{\beta_l}
\phi_{k}(\theta)d\theta ds,
$$
where
$$
\phi_k(s)=\sqrt{\frac{2k+1}{\bar T-t}}P_k\biggl(\biggl(
s-\frac{\bar T+t}{2}\biggr)\frac{2}{\bar T-t}\biggr),\ \ \ k=0, 1, 2,\ldots,
$$
where 
\begin{equation}
\label{400}
P_k(x)=\sum_{g=0}^{[k/2]} (-1)^g \frac{(2k-2g)!}{2^k g! (k-g)!(k-2g)!}\ x^{k-2g}
\end{equation}

\noindent
is the Legendre polynomial, where $[x]$ is an integer part of a real number $x.$

Further, we have
$$
C_{k}^{\sigma_j}=\int\limits_{t}^{\bar T} s^{\beta_j}\phi_{k}(s)ds=
$$
\begin{equation}
\label{700}
=
\frac{\sqrt{2k+1}}{2^{\beta_j+1}}(\bar T-t)^{\beta_j+1/2}\int\limits_{-1}^1 (s+b)^{\beta_j} P_k(s)ds,
\end{equation}
$$
C_{r k}^{\sigma_l\sigma_j}=\int\limits_{t}^{\bar T}
s^{\beta_j} \phi_{r}(s)
\int\limits_{t}^{s}
\theta^{\beta_l}
\phi_{k}(\theta)d\theta ds=
$$
\begin{equation}
\label{701}
=\frac{\sqrt{(2k+1)(2r+1)}}{2^{\beta_j+\beta_l+2}}(\bar T-t)^{\beta_j+\beta_l+1}
\int\limits_{-1}^1 (s+b)^{\beta_j} P_r(s)\int\limits_{-1}^s (\theta+b)^{\beta_l} P_k(\theta)d\theta ds,
\end{equation}
$$
x^{k-2g}=((x+b)-b)^{k-2g}=
$$
\begin{equation}
\label{500}
=\sum\limits_{h=0}^{k-2g}\frac{(k-2g)!}{h!(k-2g-h)!}(x+b)^h (-b)^{k-2g-h},
\end{equation}
where
$$
b=\frac{\bar T+t}{\bar T-t}.
$$

Let us substitute (\ref{500}) into (\ref{400})
\begin{equation}
\label{600}
P_k(x)=\sum_{g=0}^{[k/2]} \sum\limits_{h=0}^{k-2g} A_{g,h}(k,b) (x+b)^h,
\end{equation}
where
\begin{equation}
\label{900}
A_{g,h}(k,b)=\frac{(-1)^{k-g-h}(2k-2g)!b^{k-2g-h}}
{2^k g! (k-g)!h!(k-2g-h)!}.
\end{equation}

Using (\ref{600}), we obtain
\begin{equation}
\label{800}
\int\limits_{-1}^1 (s+b)^{\beta_j} P_k(s)ds=
\sum_{g=0}^{[k/2]} \sum\limits_{h=0}^{k-2g} A_{g,h}(k,b)\ \frac{(b+1)^{\beta_j+h+1} 
-(b-1)^{\beta_j+h+1}}{\beta_j+h+1},
\end{equation}
$$
\int\limits_{-1}^1 (s+b)^{\beta_j} P_r(s)\int\limits_{-1}^s (\theta+b)^{\beta_l} P_k(\theta)d\theta ds=
$$
$$
=\sum_{m_1=0}^{[k/2]} \sum\limits_{k_1=0}^{k-2m_1} \sum_{m_2=0}^{[r/2]} \sum\limits_{k_2=0}^{r-2m_2} 
\frac{A_{m_1,k_1}(k,b) A_{m_2,k_2}(r,b)}{\beta_l+k_1+1}\times
$$
$$
\times \left(\frac{(b+1)^{\beta_j+\beta_l+k_1+k_2+2}-(b-1)^{\beta_j+\beta_l+k_1+k_2+2}}
{\beta_j+\beta_l+k_1+k_2+2}-\right.
$$
\begin{equation}
\label{801}
\left.-(b-1)^{\beta_l+k_1+1}\ \frac{(b+1)^{\beta_j+k_2+1}-(b-1)^{\beta_j+k_2+1}}
{\beta_j+k_2+1}\right).
\end{equation}

\vspace{2mm}

Let $t=\tau_p=p\Delta$ $(p=0, 1,\ldots,N-1),$
$\bar T-t =\Delta = T /N$ $(N>1)$ is a constant (for simplicity) step of integration.
Then
$$
b=\frac{\bar T+t}{\bar T-t}=1+\frac{2t}{\bar T-t}=1+2p.
$$

Further, combining (\ref{700}), (\ref{701}), (\ref{800}), (\ref{801}), we finally have
$$
C_{k}^{\sigma_j}=\int\limits_{\tau_p}^{\tau_{p+1}} s^{\beta_j}\phi_{k}(s)ds=
$$
\begin{equation}
\label{1000}
=
\frac{\sqrt{2k+1}}{2^{\beta_j+1}}\Delta^{\beta_j+1/2}
\sum_{g=0}^{[k/2]} \sum\limits_{h=0}^{k-2g} A_{g,h}(k,b)
\frac{(b+1)^{\beta_j+h+1} 
-(b-1)^{\beta_j+h+1}}{\beta_j+h+1},
\end{equation}

\vspace{1mm}
$$
C_{r k}^{\sigma_l\sigma_j}=\int\limits_{\tau_p}^{\tau_{p+1}}
s^{\beta_j} \phi_{r}(s)
\int\limits_{\tau_p}^{s}
\theta^{\beta_l}
\phi_{k}(\theta)d\theta ds=
$$
$$
=\frac{\sqrt{(2k+1)(2r+1)}}{2^{\beta_j+\beta_l+2}}\Delta^{\beta_j+\beta_l+1}
\sum_{m_1=0}^{[k/2]} \sum\limits_{k_1=0}^{k-2m_1} \sum_{m_2=0}^{[r/2]} \sum\limits_{k_2=0}^{r-2m_2} 
\frac{A_{m_1,k_1}(k,b) A_{m_2,k_2}(r,b)}{\beta_l+k_1+1}\times
$$
$$
\times \left(\frac{(b+1)^{\beta_j+\beta_l+k_1+k_2+2}-(b-1)^{\beta_j+\beta_l+k_1+k_2+2}}
{\beta_j+\beta_l+k_1+k_2+2}-\right.
$$
\begin{equation}
\label{1001}
\left.-(b-1)^{\beta_l+k_1+1}\ \frac{(b+1)^{\beta_j+k_2+1}-(b-1)^{\beta_j+k_2+1}}
{\beta_j+k_2+1}\right),
\end{equation}

\vspace{2mm}
\noindent
where $A_{g,h}(k,b)$ is defined by (\ref{900}), $b=1+2p,$ and $p=0,1,\ldots,N-1.$

It should be noted that the Fourier--Legendre coefficients (\ref{1000}) and (\ref{1001})
depend on the integration step number $p.$

In conclusion of this section, we provide formulas 
for approximations 
of iterated It\^{o} stochastic integrals $I_{\tau_{p+1},\tau_p}^{(j)},$ 
$I_{\tau_{p+1},\tau_p}^{(l j)}$ $(l,j=1,\ldots,m)$ of the form (\ref{int1})
obtained using Theorem~1 and 
a complete orthonormal system of Legendre polynomials in the space $L_2([\tau_p,\tau_{p+1}],{\bf R})$
\cite{15}, \cite{23}, \cite{23a}, \cite{25}, \cite{25a} (for the first time, these formulas
were obtained in \cite{old-art-1} (1997) using another approach)
\begin{equation}
\label{15000a}
I_{\tau_{p+1},\tau_p}^{(j)}=\sqrt{\Delta}\zeta_0^{(j)},
\end{equation}
\begin{equation}
\label{15000b}
I_{\tau_{p+1},\tau_p}^{(l j)q}
=\frac{\Delta}{2}\left(\zeta_0^{(l)}\zeta_0^{(j)}+\sum_{i=1}^{q}
\frac{1}{\sqrt{4i^2-1}}\left(
\zeta_{i-1}^{(l)}\zeta_{i}^{(j)}-
\zeta_i^{(l)}\zeta_{i-1}^{(j)}\right) - {\bf 1}_{\{l=j\}}\right),
\end{equation}
where $q\in {\bf N},$ 
$$
\zeta_{i}^{(l)}=\int\limits_{\tau_p}^{\tau_{p+1}}\phi_i(s)d{\bf w}_s^{(l)}\ \ \ 
(l=1,\ldots,m;\ i=0,1,2,\ldots)
$$ 
are independent standard Gaussian random variables
for various
$i$ or $l,$
$\{\phi_i(s)\}_{i=0}^{\infty}$
is a complete orthonormal system of Legendre polynomials 
in the space $L_2([\tau_p,\tau_{p+1}],{\bf R}).$

Moreover \cite{15}, \cite{23}, \cite{23a}, \cite{25}, \cite{25a}, \cite{old-art-1}
\begin{equation}
\label{15000c}
{\sf M}\biggl\{\left(I_{\tau_{p+1},\tau_p}^{(l j)}-
I_{\tau_{p+1},\tau_p}^{(l j)q}
\right)^2\biggr\}
=\frac{\Delta^2}{2}\left(\frac{1}{2}-\sum_{i=1}^q
\frac{1}{4i^2-1}\right)\ \ \ (l\ne j).
\end{equation}

Also note that for $l=j$ the following equality
\begin{equation}
\label{15000d}
I_{\tau_{p+1},\tau_p}^{(j j)}=\frac{1}{2}\left(\left(I_{\tau_{p+1},\tau_p}^{(j)}\right)^2-\Delta\right)=
\frac{1}{2}\Delta\left(\left(\zeta_0^{(j)}\right)^2-1\right)
\end{equation}
holds w.~p.~1.

\section{The Taylor--It\^{o} Expansion for It\^{o} SDEs with Non-Differentiable Time-Dependent
Diffusion Coefficients}

Consider the following special case of the It\^{o} SDE (\ref{1.5.2s})
\begin{equation}
\label{1.5.2sqqq}
{\bf x}_t={\bf x}_0+\int\limits_0^t {\bf a}({\bf x}_{s},s)ds+
\sum\limits_{j=1}^m
\int\limits_0^t \sigma_j(s) B_j({\bf x}_{s},s)d{\bf w}_{s}^{(j)},\ \ \
{\bf x}_0={\bf x}(0,\omega),
\end{equation}

\noindent
where $t\in [0, T],$ the functions $\sigma_j(t): [0, T]\to{\bf R}$ $(j=1,\ldots,m)$ 
are continuous but non-differentiable on $[0, T],$ whereas
$\sigma_j^2(t)\in C^{\infty}([0, T])$ (the space of functions 
with continuous derivatives of all orders on $[0, T]$).
Furthermore, 
the functions ${\bf a}({\bf x},t),$
$B_j({\bf x},t)$ $(j=1,\ldots,m):$ 
${\bf R}^n\times[0, T]\to{\bf R}^n$
are measurable with respect to the collection of 
variables $({\bf x},t)$, satisfy the global Lipschitz condition, and the linear growth condition;
another notations 
appearing in (\ref{1.5.2sqqq}) are the same as in (\ref{1.5.2}).

The above conditions 
guarantee the existence and uniqueness up to stochastic equivalence 
of a strong solution
of (\ref{1.5.2sqqq}) \cite{1982}.

As examples of functions $\sigma_j(t)$ $(j=1,\ldots,m)$ satisfying 
the above conditions,
one can take
$\sqrt{t-t_0},$ $\left|t-t_0\right|$ $(t_0\in [0, T]),$ 
$\left|t\right|\cos(t),$ 
$\left|\sin (at)\right|$ $(a=10\pi/T),$ $\left|\cos (at)\right|$ $(a=19\pi/2T),$
$\sqrt{(t-t_1)\ldots(t-t_k)},$ 
$\left|(t-t_1)\ldots(t-t_k)\right|$ $(t_1,\ldots,t_k\in [0, T]),$ etc.

Let us assume that the functions ${\bf a}({\bf x},t)$ and
$B_j({\bf x},t)$ $(j=1,\ldots,m)$ are sufficiently smooth with respect 
to the variables ${\bf x}$ and $t.$

Consider the following differential operators
$$
L= {\partial \over \partial t}
+ \sum^ {n} _ {i=1} a^{(i)} ({\bf x},  t) 
{\partial \over  \partial  {\bf  x}^{(i)}} +
{1\over 2} \sum^ {m} _ {j=1}\sigma_j^2(t) \sum^ {n} _ {l,i=1}
B^ {(lj)} ({\bf x}, t) B^ {(ij)} ({\bf x}, t) {\partial
^{2} \over\partial{\bf x}^{(l)}\partial{\bf x}^{(i)}},
$$
$$
G_0^{(j)}=\sigma_j(t)G^{(j)},\ \ \ 
G^{(j)}=\sum^{n}_{i=1} B^{(ij)}({\bf x}, t)
{\partial \over \partial {\bf x}^{(i)}},\ \ \
j=1,\ldots,m.
$$

We also introduce the following notations
$$
{M}_k=\biggl\{(\lambda_k,\ldots,\lambda_1):\
\lambda_l=1\ \hbox{or}\ \lambda_l=0,\ l=1,\ldots,k\biggr\},\ \ \ k\ge 1,
$$
\begin{equation}
\label{9191}
J_{(\lambda_{k}\ldots \lambda_1)s,t}^{\sigma_{i_k}\ldots\sigma_{i_1}(i_k\ldots
i_1)}=
\int\limits_t^s\sigma_{i_1}(t_1)\int\limits_t^{t_1}\sigma_{i_2}(t_2)\ldots
\int\limits_t^{t_{k-1}}\sigma_{i_k}(t_k)
d{\bf w}_{t_{k}}^{(i_k)}\ldots
d{\bf w}_{t_2}^{(i_2)}d{\bf w}_{t_1}^{(i_1)},\ \ \ k\ge 1,
\end{equation}

\noindent
where $i_1,\ldots,i_k=0,1,\ldots,m,$ $\lambda_l=1$ or $\lambda_l=0$,
$Q_{\lambda_l}^{(i_l)}\stackrel{\sf def}{=}{L}$
and $i_l=0$ for $\lambda_l=0,$  $Q_{\lambda_l}^{(i_l)}\stackrel{\sf def}{=}G_0^{(i_l)}$
and $i_l=1,\ldots,m$ for $\lambda_l=1,$ 
${\bf w}_{\tau}^{(i)}$ $(i=1,\ldots,m)$ are
independent standard Wiener processes,
${\bf w}_{\tau}^{(0)}=\tau,$ and $\sigma_0(t)\equiv 1.$

Suppose that the function $R({\bf x}, t): {\bf R}^n\times [0, T]\to {\bf R}$ is 
sufficiently smooth with respect 
to the variables ${\bf x}$ and $t.$

Consider the stochastic process $R({\bf x}_t, t),$ where ${\bf x}_t$
is the strong solution to the It\^{o} SDE (\ref{1.5.2sqqq}).

Iteratively 
applying It\^{o}'s formula (as in Sect.~3) to the process 
$R({\bf x}_s, s),$ we obtain the following 
Taylor--It\^{o} expansion
$$
R({\bf x}_s,s)=R({\bf x}_t,t)
+
$$
\begin{equation}
\label{7000}
+\sum_{k=1}^r \sum_{(\lambda_{k},\ldots,\lambda_1)
\in{M}_k}
\sum_{i_1=\lambda_1}^{m\lambda_1}
\ldots 
\sum_{i_k=\lambda_k}^{m\lambda_k}
Q_{\lambda_k}^{(i_k)}\ldots Q_{\lambda_1}^{(i_1)}
R({\bf x}_t,t){J}_{(\lambda_{k}\ldots \lambda_1)s,t}^{\sigma_{i_k}\ldots\sigma_{i_1}(i_k\ldots i_1)}
+\left(D_{r+1}\right)_{s,t}
\end{equation}

\vspace{2mm}
\noindent
w.~p.~1, where $s>t,$ $r\in {\bf N},$ and
$$
\left(D_{r+1}\right)_{s,t}=
$$
$$
=
\sum_{(\lambda_{r+1},\ldots,\lambda_{1})\in{M}_{r+1}}
\sum_{i_1=\lambda_1}^{m\lambda_1}
\ldots 
\sum_{i_{r+1}=\lambda_{r+1}}^{m\lambda_{r+1}}
\int\limits_{t}^{s}\int\limits_{t}^{t_1}\ldots
\int\limits_{t}^{t_r}
Q_{\lambda_{r+1}}^{(i_{r+1})}\ldots
Q_{\lambda_{1}}^{(i_1)}
R({\bf x}_{t_{r+1}}, t_{r+1})\times
$$
$$
\times
d{\bf w}_{t_{r+1}}^{(i_{r+1})}
\ldots
d{\bf w}_{t_{2}}^{(i_2)}d{\bf w}_{t_{1}}^{(i_1)}.
$$

\vspace{2mm}

Equality (\ref{7000}) can be used to construct high-order
strong numerical methods for It\^{o} SDE (\ref{1.5.2sqqq}), while iterated 
It\^{o} stochastic integrals (\ref{9191}) can be approximated 
using Theorem~1.

\section{The Taylor--It\^{o} Expansion for It\^{o} SDEs with 
Drift and Diffusion Coefficients Having 
Singular Time-Dependent Multipliers}

Now consider an the It\^{o} SDE of the form
\begin{equation}
\label{1.5.2sqqqx}
{\bf x}_t={\bf x}_0+\int\limits_0^t \sigma_0(s){\bf a}({\bf x}_{s})ds+
\sum\limits_{j=1}^m
\int\limits_0^t \sigma_j(s) B_j({\bf x}_{s})d{\bf w}_{s}^{(j)},\ \ \
{\bf x}_0={\bf x}(0,\omega),
\end{equation}
                                               
\noindent
where $t\in [0, T],$ 
$\sigma_0(t)\in L_1([0, T],{\bf R}),$ $\sigma_j(t)\in L_2([0, T],{\bf R})$ $(j=1,\ldots,m),$
the functions ${\bf a}({\bf x}),$
$B_j({\bf x})$ $(j=1,\ldots,m):$ 
${\bf R}^n\to{\bf R}^n$
are measurable, satisfy the global Lipschitz condition, and the linear growth condition;
another notations 
appearing in (\ref{1.5.2sqqqx}) are the same as in (\ref{1.5.2}).

According to Theorem~1.1 from \cite{Jan}, under the above conditions, the It\^{o} SDE (\ref{1.5.2sqqqx})
has a unique strong solution.

Consider the following differential operators 
$$
L= 
\sigma_0(t) \sum^{n}_{l=1}a^{(l)}({\bf x}) 
{\partial \over  \partial  {\bf  x}^{(l)}} +
{1\over 2} \sum^{m}_{j=1} \sigma_j^2(t)\sum^{n}_{l,k=1}
B^{(lj)}({\bf x}) B^ {(kj)}({\bf x}) 
{\partial^{2} \over\partial{\bf x}^{(l)}\partial{\bf x}^{(k)}}\stackrel{\sf def}{=}
$$
$$
\stackrel{\sf def}{=}\sigma_0(t)L_0 + \sum^{m}_{j=1} \sigma_j^2(t)L_j,
$$
\begin{equation}
\label{oper10}
G_0^{(j)}=\sigma_j(t)G^{(j)},\ \ \ G^{(j)}=\sum^{n}_{k=1} B^{(kj)}({\bf x})
{\partial \over \partial {\bf x}^{(k)}}\ \ \
(j=1,\ldots,m),
\end{equation}

\noindent
where
\begin{equation}
\label{oper11}
L_0=\sum^{n}_{l=1}a^{(l)}({\bf x}) 
{\partial \over  \partial  {\bf  x}^{(l)}},\ \ \ 
L_j={1\over 2}\sum^{n}_{l,k=1}
B^{(lj)}({\bf x}) B^ {(kj)}({\bf x}) 
{\partial^{2} \over\partial{\bf x}^{(l)}\partial{\bf x}^{(k)}}\ \ \ (j=1,\ldots,m),
\end{equation}
${\bf x}^{(i)}$ is the $i$th component of
${\bf x}$, ${\bf a}^{(i)}({\bf x})$ is the $i$th component
of ${\bf a}({\bf x})$, and
$B^{(ij)}({\bf x})$ is the $ij$th element of the matrix $B({\bf x})$
with columns $B_j({\bf x})$.

Let us assume that the functions ${\bf a}({\bf x})$ and
$B_j({\bf x})$ $(j=1,\ldots,m)$ are sufficiently smooth with respect 
to ${\bf x}$.

From (\ref{1.5.2sqqqx}) we obtain for $t>\tau$ w.~p.~1
\begin{equation}
\label{1.2sqqqx}
{\bf x}_t={\bf x}_{\tau}+\int\limits_{\tau}^t \sigma_0(s){\bf a}({\bf x}_{s})ds+
\sum\limits_{j=1}^m
\int\limits_{\tau}^t \sigma_j(s) B_j({\bf x}_{s})d{\bf w}_{s}^{(j)}.
\end{equation}

Applying It\^{o}'s formula to ${\bf a}({\bf x}_{s})$ and $B_j({\bf x}_{s}),$
we get for $s>\tau$ w.~p.~1
$$
{\bf a}({\bf x}_{s})={\bf a}({\bf x}_{\tau})+
\int\limits_{\tau}^s \sigma_0(\theta)L_0{\bf a}({\bf x}_{\theta})d\theta+
\sum\limits_{l=1}^m\int\limits_{\tau}^s \sigma_l^2(\theta)L_l{\bf a}({\bf x}_{\theta})d\theta+
$$
\begin{equation}
\label{5.1000}
+
\sum\limits_{l=1}^m
\int\limits_{\tau}^s \sigma_l(\theta)G^{(l)}{\bf a}({\bf x}_{\theta})d{\bf w}_{\theta}^{(l)},
\end{equation}
$$
B_j({\bf x}_{s})=B_j({\bf x}_{\tau})+
\int\limits_{\tau}^s \sigma_0(\theta)L_0 B_j({\bf x}_{\theta})d\theta+
\sum\limits_{l=1}^m\int\limits_{\tau}^s \sigma_l^2(\theta)L_l B_j({\bf x}_{\theta})d\theta+
$$
\begin{equation}
\label{5.1001}
+
\sum\limits_{l=1}^m
\int\limits_{\tau}^s \sigma_l(\theta)G^{(l)} B_j({\bf x}_{\theta})d{\bf w}_{\theta}^{(l)}.
\end{equation}

Substituting (\ref{5.1000}) and (\ref{5.1001}) into (\ref{1.2sqqqx}), we obtain w.~p.~1
$$
{\bf x}_t={\bf x}_{\tau}+{\bf a}({\bf x}_{\tau})\int\limits_{\tau}^t \sigma_0(s)ds+
B_j({\bf x}_{\tau})\sum\limits_{j=1}^m\int\limits_{\tau}^t \sigma_j(s)d{\bf w}_{\theta}^{(j)}+
$$
$$
+
\int\limits_{\tau}^t \sigma_0(s)\int\limits_{\tau}^s \sigma_0(\theta)L_0{\bf a}({\bf x}_{\theta})d\theta ds+
\sum\limits_{l=1}^m
\int\limits_{\tau}^t \sigma_0(s)\int\limits_{\tau}^s \sigma_l^2(\theta)L_l{\bf a}({\bf x}_{\theta})
d\theta ds+
$$
$$
+\hspace{-0.7mm}
\sum\limits_{l=1}^m
\int\limits_{\tau}^t \sigma_0(s)
\int\limits_{\tau}^s \sigma_l(\theta)G^{(l)}{\bf a}({\bf x}_{\theta})d{\bf w}_{\theta}^{(l)}ds+
\sum\limits_{j=1}^m
\int\limits_{\tau}^t \sigma_j(s)\int\limits_{\tau}^s \sigma_0(\theta)L_0 B_j({\bf x}_{\theta})d\theta 
d{\bf w}_{\theta}^{(j)}\hspace{-0.7mm}+
$$
$$
+
\sum\limits_{l, j=1}^m
\int\limits_{\tau}^t \sigma_j(s)\int\limits_{\tau}^s \sigma_l^2(\theta)L_l B_j({\bf x}_{\theta})d\theta
d{\bf w}_{\theta}^{(j)}+
$$
\begin{equation}
\label{5.1002}
+
\sum\limits_{l,j=1}^m
\int\limits_{\tau}^t \sigma_j(s)
\int\limits_{\tau}^s \sigma_l(\theta)G^{(l)} B_j({\bf x}_{\theta})d{\bf w}_{\theta}^{(l)} 
d{\bf w}_{\theta}^{(j)}.
\end{equation}

\vspace{3mm}

If we now apply It\^{o}'s formula to the processes 
$$
L_0{\bf a}({\bf x}_{\theta}),\ \ L_l{\bf a}({\bf x}_{\theta}),\ \ G^{(l)}{\bf a}({\bf x}_{\theta}),\ \ 
L_0 B_j({\bf x}_{\theta}),\ \ L_l B_j({\bf x}_{\theta}),\ \ G^{(l)} B_j({\bf x}_{\theta}),
$$

\noindent
substitute 
the result into (\ref{5.1002}), etc., we can obtain the Taylor--It\^{o} expansion 
(see (\ref{9192}), (\ref{9192a}))
of the solution ${\bf x}_{t}$ to the It\^{o} SDE (\ref{1.5.2sqqqx}) in terms 
of iterated It\^{o} stochastic integrals of the form
$$
J_{(\lambda_k\ldots\lambda_1)t,\tau}^{\psi_{i_k,j_k}\ldots \psi_{i_1,j_1}
((i_k,j_k)\ldots (i_1,j_1))}=
$$
\begin{equation}
\label{int100}
=\int\limits_{\tau}^t \psi_{i_1,j_1}(t_1)\int\limits_{\tau}^{t_1} \psi_{i_2,j_2}(t_2)\ldots
\int\limits_{\tau}^{t_{k-1}}\psi_{i_k,j_k}(t_k)
d{\bf w}_{t_{k}}^{(i_k,j_k)}\ldots
d{\bf w}_{t_2}^{(i_2,j_2)}d{\bf w}_{t_1}^{(i_1,j_1)},
\end{equation}

\noindent
where $k\ge 1,$ $i_1,\ldots,i_k=0,1,\ldots,m,$ 
$j_1,\ldots,j_k=0,1,$
${\bf w}_{s}^{(i,1)}\stackrel{\sf def}{=}{\bf w}_{s}^{(i)}$ $(i=1,\ldots,m)$ are independent
standard Wiener processes, ${\bf w}_{s}^{(i,0)}\stackrel{\sf def}{=}s$ $(i=1,\ldots,m),$ 
${\bf w}_{s}^{(0,0)}\stackrel{\sf def}{=}s,$ $(\lambda_k,\ldots,\lambda_1)\in {M}_k,$
$$
{M}_k=\biggl\{(\lambda_k,\ldots,\lambda_1):\
\lambda_l=1\ \hbox{or}\ \lambda_l=0,\ l=1,\ldots,k\biggr\},
$$

\noindent
and
$$
\psi_{i_l,j_l}(s)=\left\{
\begin{matrix}
\sigma_{i_l}(s) \ \ & \hbox{for} \ \ i_l=1,\ldots,m \ \ \hbox{and} \ \  j_l=1\ \  (\hbox{if}\ \ \lambda_l=1)
\cr \cr
\sigma_{i_l}^2(s) \ \ & \hbox{for} \ \ i_l=1,\ldots,m \ \ \hbox{and} \ \  j_l=0\ \  (\hbox{if}\ \ \lambda_l=1)
\cr \cr
\sigma_0(s) \ \ & \hbox{for} \ \ i_l=0 \ \ \hbox{and} \ \  j_l=0\ \  (\hbox{if}\ \ \lambda_l=0) 
\end{matrix}\right.,
$$

\noindent
where $l=1,\ldots,k.$ 

\vspace{2mm}

Let us consider several examples of integrals of the form (\ref{int100})
$$
\int\limits_{\tau}^t \sigma_{i_1}(t_1)\int\limits_{\tau}^{t_1} \sigma_{i_2}(t_2)
d{\bf w}_{t_2}^{(i_2)}d{\bf w}_{t_1}^{(i_1)},\ \ \ 
\int\limits_{\tau}^t \sigma_{i_1}(t_1)\int\limits_{\tau}^{t_1} \sigma_{i_2}^2(t_2)
dt_2d{\bf w}_{t_1}^{(i_1)},\ \ \ 
$$
$$
\int\limits_{\tau}^t \sigma_{i_1}(t_1)\int\limits_{\tau}^{t_1} \sigma_{i_2}(t_2)
\int\limits_{\tau}^{t_{2}}\sigma_{i_3}^2(t_3)dt_3 
d{\bf w}_{t_2}^{(i_2)}d{\bf w}_{t_1}^{(i_1)},
$$
$$
\int\limits_{\tau}^t \sigma_0(t_1)\int\limits_{\tau}^{t_1} \sigma_{i_2}(t_2)
\int\limits_{\tau}^{t_{2}}\sigma_{i_3}(t_3)d{\bf w}_{t_3}^{(i_3)}
d{\bf w}_{t_2}^{(i_2)}dt_1,
$$
$$
\int\limits_{\tau}^t \sigma_{i_1}(t_1)\int\limits_{\tau}^{t_1} \sigma_{i_2}(t_2)
\int\limits_{\tau}^{t_{2}}\sigma_{i_3}(t_3)d{\bf w}_{t_3}^{(i_3)}
d{\bf w}_{t_2}^{(i_2)}d{\bf w}_{t_1}^{(i_1)},
$$
where $i_1,i_2,i_3=1,\ldots,m.$

Since 
$$
\sigma_0(t)\in L_1([0, T],{\bf R}),\ \ \ \sigma_j(t)\in L_2([0, T],{\bf R})\ \ \ (j=1,\ldots,m),
$$

\noindent
the iterated It\^{o} stochastic integrals (\ref{int100})
exist in the mean-square sense and can be approximated using Theorem~1
(provided that $t$ and $\tau$ are constant time instants).

Denote
$$
Q_{\lambda_l}^{(i_l,j_l)}=\left\{
\begin{matrix}
G^{(i_l)}\ \  & \hbox{for}\ \  i_l=1,\ldots,m\ \  \hbox{and}\ \   j_l=1\ \  (\hbox{if}\ \ \lambda_l=1)
\cr \cr 
L_{i_l}\ \ & \hbox{for}\ \ i_l=1,\ldots,m\ \ \hbox{and}\ \  j_l=0\ \ (\hbox{if}\ \ \lambda_l=1)
\cr \cr
L_0\ \ & \hbox{for}\ \ i_l=0\ \ \hbox{and}\ \  j_l=0\ \ (\hbox{if}\ \ \lambda_l=0) 
\end{matrix}\right.,
$$

\noindent
where $l=1,\ldots,k.$ 

Note that
\begin{equation}
\label{ww1003}
L_0{\bf x}_{\tau}={\bf a}({\bf x}_{\tau}),\ \ \ G^{(i_1)}{\bf x}_{\tau}=B_{i_1}({\bf x}_{\tau}),\ \ \ 
L_{i_1}{\bf x}_{\tau}\equiv {\bf 0},
\end{equation}
where ${\bf 0}$ is a column vector of length $n$ consisting of zeros, and
$i_1=1,\ldots,m$.

Suppose that the function $R({\bf x}): {\bf R}^n\to {\bf R}$ is 
sufficiently smooth with respect 
to ${\bf x}$.
Consider the stochastic process $R({\bf x}_t),$ where ${\bf x}_t$
is the strong solution to the It\^{o} SDE (\ref{1.5.2sqqqx}).

Iteratively 
applying It\^{o}'s formula (as above (see (\ref{5.1002})) to the process 
$R({\bf x}_t),$ we obtain the following 
Taylor--It\^{o} expansion
$$
R({\bf x}_t)=R({\bf x}_{\tau})
+
$$
$$
+\sum_{k=1}^r \sum_{(\lambda_{k},\ldots,\lambda_1)
\in{M}_k}
\sum_{i_1=\lambda_1}^{m\lambda_1}\sum_{j_1=0}^{\lambda_1}
\ldots 
\sum_{i_k=\lambda_k}^{m\lambda_k}\sum_{j_k=0}^{\lambda_k}
Q_{\lambda_k}^{(i_k,j_k)}\ldots Q_{\lambda_1}^{(i_1,j_1)}
R({\bf x}_{\tau})\times
$$

\begin{equation}
\label{9192}
\times 
J_{(\lambda_k\ldots\lambda_1)t,\tau}^{\psi_{i_k,j_k}\ldots \psi_{i_1,j_1}
((i_k,j_k)\ldots (i_1,j_1))}
+\left(D_{r+1}\right)_{t,\tau}
\end{equation}

\vspace{2mm}
\noindent
w.~p.~1, where $t>\tau,$ $r\in {\bf N},$ and
$$
\left(D_{r+1}\right)_{t,\tau}=
$$
$$
=
\sum_{(\lambda_{r+1},\ldots,\lambda_{1})\in{M}_{r+1}}
\sum_{i_1=\lambda_1}^{m\lambda_1}\sum_{j_1=0}^{\lambda_1}
\sum_{i_2=\lambda_2}^{m\lambda_2}\sum_{j_2=0}^{\lambda_2}
\ldots 
\sum_{i_{r+1}=\lambda_{r+1}}^{m\lambda_{r+1}}\sum_{j_{r+1}=0}^{\lambda_{r+1}}
\int\limits_{\tau}^{t}\psi_{i_1,j_1}(t_1)\int\limits_{\tau}^{t_1}\psi_{i_2,j_2}(t_2)\ldots
$$
$$
\ldots
\int\limits_{\tau}^{t_r}
\psi_{i_{r+1},j_{r+1}}(t_{r+1})Q_{\lambda_{r+1}}^{(i_{r+1},j_{r+1})}\ldots
Q_{\lambda_{2}}^{(i_2,j_2)}Q_{\lambda_{1}}^{(i_1,j_1)}
R({\bf x}_{t_{r+1}})\times
$$
\begin{equation}
\label{9192a}
\times
d{\bf w}_{t_{r+1}}^{(i_{r+1},j_{r+1})}
\ldots
d{\bf w}_{t_{2}}^{(i_2,j_2)}d{\bf w}_{t_{1}}^{(i_1,j_1)}.
\end{equation}

\vspace{2mm}

Note that for $R({\bf x}_{\tau})={\bf x}_{\tau},$ the expression
$Q_{\lambda_1}^{(i_1,j_1)}{\bf x}_{\tau}$ is defined by the equality (by virtue of (\ref{ww1003}))

\vspace{-4mm}
\begin{equation}
\label{ww1004}
Q_{\lambda_1}^{(i_1,j_1)}{\bf x}_{\tau}=\left\{
\begin{matrix}
B_{i_1}({\bf x}_{\tau})\ \  & \hbox{for}\ \  i_1=1,\ldots,m\ \  \hbox{and}\ \   j_1=1\ \  (\hbox{if}\ \ 
\lambda_l=1)
\cr \cr 
{\bf 0}\ \ & \hbox{for}\ \ i_1=1,\ldots,m\ \ \hbox{and}\ \  j_1=0\ \ (\hbox{if}\ \ \lambda_1=1)
\cr \cr
{\bf a}({\bf x}_{\tau})\ \ & \hbox{for}\ \ i_1=0\ \ \hbox{and}\ \  j_1=0\ \ (\hbox{if}\ \ \lambda_1=0) 
\end{matrix}\right..
\end{equation}

\vspace{2mm}

The existence in the mean-square sense of the remainder
$\left(D_{r+1}\right)_{t,\tau}$ in the integral form (\ref{9192a}) of the Taylor--It\^{o}
expansion (\ref{9192}) of the solution to the It\^{o} SDE (\ref{1.5.2sqqqx})
will depend on the growth conditions of the functions 
$R({\bf x}),$ ${\bf a}^{(i)}({\bf x}),$
$B^{(ij)}({\bf x})$ $(i=1,\ldots,n,\ j=1,\ldots,m)$
and their partial derivatives. Here 
${\bf a}^{(i)}({\bf x})$ is the $i$th component
of ${\bf a}({\bf x})$ and $B^{(ij)}({\bf x})$
is the $ij$th element of matrix $B({\bf x})$ with columns $B_j({\bf x})$ $(j=1,\ldots,m).$

For instance, if the functions $R({\bf x}),$ ${\bf a}^{(i)}({\bf x}),$
$B^{(ij)}({\bf x})$ $(i=1,\ldots,n,\ j=1,\ldots,m)$
together with all their partial 
derivatives are uniformly bounded on $[0, T],$ then analogously to the 
derivation of (\ref{gg111}) under condition (\ref{gg111ddd}), we obtain that the 
specified remainder $\left(D_{r+1}\right)_{t,\tau}$ exists in the mean-square sense provided that 
$\sigma_0(t)\in L_1([0, T],{\bf R}),$ $\sigma_j(t)\in L_2([0, T],{\bf R})$ $(j=1,\ldots,m).$

If we relax the growth conditions on the functions 
$R({\bf x}),$ ${\bf a}^{(i)}({\bf x}),$
$B^{(ij)}({\bf x})$ $(i=1,\ldots,n,\ j=1,\ldots,m)$ and their partial derivatives, 
it will be necessary to take into account the existence of certain 
moments ${\sf M}\left\{\left\Vert{\bf x}_t\right\Vert^r\right\}$  
of the solution ${\bf x}_t$ to the It\^{o} SDE (\ref{1.5.2sqqqx}) (see condition (\ref{12321})). 
This will significantly restrict the conditions on the functions 
$\sigma_0(t),$ $\sigma_j(t),$ $j=1,\ldots,m.$

For example, if the functions 
$Q_{\lambda_{r+1}}^{(i_{r+1},j_{r+1})}\ldots
Q_{\lambda_{2}}^{(i_2,j_2)}Q_{\lambda_{1}}^{(i_1,j_1)}
R({\bf x})$ included in $\left(D_{r+1}\right)_{t,\tau}$ satisfy the 
linear growth condition (\ref{qw11}), then ${\sf M}\bigl\{\left\Vert{\bf x}_t\right\Vert^2\bigr\}$   
must exist for the existence of $\left(D_{r+1}\right)_{t,\tau}$ in the form (\ref{9192a})
in the mean-square sense. This implies the conditions
$\sigma_0^{2+\varepsilon}(t)\in L_1([0, T],{\bf R}),$ 
$\sigma_j^{4+2\varepsilon}(t)\in L_1([0, T],{\bf R})$ ($j=1,\ldots,m$),
where $\varepsilon$ is an arbitrary small positive real number
(see  (\ref{12321dd})).

Suppose that $R({\bf x})\equiv {\bf x}^{(i)},$ where ${\bf x}^{(i)}$
is the $i$th component of ${\bf x}$ $(i=1,\ldots,m).$
Then, using (\ref{ww1003}), (\ref{ww1004}), 
and the truncated Taylor--It\^{o} expansion (\ref{9192}) 
(without the remainder term $\left(D_{r+1}\right)_{t,\tau}$), one can construct high-order 
stong numerical methods for It\^{o} SDE (\ref{1.5.2sqqqx}), similarly to usual It\^{o} SDEs
\cite{2},\cite{15}-\cite{15b}.

Consider a more general It\^{o} SDE than (\ref{1.5.2sqqqx}) in the form

\newpage
\noindent
\begin{equation}
\label{1600sqqqx}
{\bf x}_t={\bf x}_0+\int\limits_0^t \sigma_0(s){\bf a}({\bf x}_{s},s)ds+
\sum\limits_{j=1}^m
\int\limits_0^t \sigma_j(s) B_j({\bf x}_{s},s)d{\bf w}_{s}^{(j)},\ \ \
{\bf x}_0={\bf x}(0,\omega),
\end{equation}
                                               
\noindent
where $t\in [0, T],$ 
$\sigma_0(t)\in L_1([0, T],{\bf R}),$ $\sigma_j(t)\in L_2([0, T],{\bf R})$ $(j=1,\ldots,m),$
the functions ${\bf a}({\bf x},t),$
$B_j({\bf x},t)$ $(j=1,\ldots,m):$ 
${\bf R}^n\times [0, T]\to{\bf R}^n$
are measurable with respect to the collection of 
variables $({\bf x},t),$ satisfy the global Lipschitz condition, and the linear growth condition;
another notations 
appearing in (\ref{1600sqqqx}) are the same as in (\ref{1.5.2}).

According to Theorem~1.1 from \cite{Jan}, under the above conditions, the It\^{o} SDE (\ref{1.5.2sqqqx})
has a unique strong solution.

Suppose that the functions $R({\bf x},t): {\bf R}^n\times [0, T]\to {\bf R},$
${\bf a}({\bf x},t)$ and
$B_j({\bf x},t)$ $(j=1,\ldots,m)$ are sufficiently smooth with respect 
to the variables ${\bf x}$ and $t.$ 

Let us generalize the Taylor--It\^{o} expansion (\ref{9192}), (\ref{9192a})
to the case of the stochastic process
$R({\bf x}_t,t),$ where ${\bf x}_t$
is the strong solution to the It\^{o} SDE (\ref{1600sqqqx}).

Consider the following differential operators 
$$
L=\bar L_0+\sigma_0(t)L_0 + \sum^{m}_{j=1} \sigma_j^2(t)L_j,
$$
$$
G_0^{(j)}=\sigma_j(t)G^{(j)}\ \ \ 
(j=1,\ldots,m),
$$

\vspace{2mm}
\noindent
where $L_0,$ $L_j,$ $G^{(j)}$ $(j=1,\ldots,m)$ are defined by (\ref{oper10}), (\ref{oper11}), and
\begin{equation}
\label{Mil14}
\bar L_0=\frac{\partial}{\partial t}.
\end{equation}

Consider the following iterated It\^{o} stochastic integrals
$$
\bar J_{(\lambda_k\ldots\lambda_1)t,\tau}^{~\bar \psi_{i_k,j_k}\ldots \bar\psi_{i_1,j_1}
((i_k,j_k)\ldots (i_1,j_1))}=
$$
\begin{equation}
\label{int100aaa}
=\int\limits_{\tau}^t \bar\psi_{i_1,j_1}(t_1)\int\limits_{\tau}^{t_1} \bar\psi_{i_2,j_2}(t_2)\ldots
\int\limits_{\tau}^{t_{k-1}}\bar\psi_{i_k,j_k}(t_k)
d\bar {\bf w}_{t_{k}}^{(i_k,j_k)}\ldots
d\bar{\bf w}_{t_2}^{(i_2,j_2)}d\bar{\bf w}_{t_1}^{(i_1,j_1)},
\end{equation}
where $k\ge 1,$ $i_1,\ldots,i_k=0,1,\ldots,m,$ 
$j_1,\ldots,j_k=-1, 0 ,1,$
$\bar{\bf w}_{s}^{(i,1)}\stackrel{\sf def}{=}{\bf w}_{s}^{(i)}$ $(i=1,\ldots,m)$ are independent
standard Wiener processes, $\bar{\bf w}_{s}^{(i,0)}\stackrel{\sf def}{=}s$ $(i=1,\ldots,m),$ 
$\bar{\bf w}_{s}^{(0,0)}\stackrel{\sf def}{=}s,$ $\bar{\bf w}_{s}^{(0,-1)}\stackrel{\sf def}{=}s,$
$(\lambda_k,\ldots,\lambda_1)\in {M}_k,$
$$
{M}_k=\biggl\{(\lambda_k,\ldots,\lambda_1):\
\lambda_l=1\ \hbox{or}\ \lambda_l=0,\ l=1,\ldots,k\biggr\},
$$
and
$$
\bar \psi_{i_l,j_l}(s)=\left\{
\begin{matrix}
\sigma_{i_l}(s) \ \ & \hbox{for} \ \ i_l=1,\ldots,m \ \ \hbox{and} \ \  j_l=1\ \  (\hbox{if}\ \ \lambda_l=1)
\cr \cr
\sigma_{i_l}^2(s) \ \ & \hbox{for} \ \ i_l=1,\ldots,m \ \ \hbox{and} \ \  j_l=0\ \  (\hbox{if}\ \ \lambda_l=1)
\cr \cr
\sigma_0(s) \ \ & \hbox{for} \ \ i_l=0 \ \ \hbox{and} \ \  j_l=0\ \  (\hbox{if}\ \ \lambda_l=0)
\cr \cr
1 \ \ & \hbox{for} \ \ i_l=0 \ \ \hbox{and} \ \  j_l=-1\ \  (\hbox{if}\ \ \lambda_l=0) 
\end{matrix}\right.,
$$
where $l=1,\ldots,k.$

Let us consider several examples of integrals of the form (\ref{int100aaa})
$$
\int\limits_{\tau}^t \sigma_{i_1}(t_1)\int\limits_{\tau}^{t_1} \sigma_{i_2}(t_2)
d{\bf w}_{t_2}^{(i_2)}d{\bf w}_{t_1}^{(i_1)},\ \ 
\int\limits_{\tau}^t \sigma_{i_1}(t_1)\int\limits_{\tau}^{t_1}
dt_2d{\bf w}_{t_1}^{(i_1)},
$$
$$
\int\limits_{\tau}^t \sigma_0(t_1)\int\limits_{\tau}^{t_1}
dt_2dt_1,\ \ \
\int\limits_{\tau}^t \sigma_{0}(t_1)\int\limits_{\tau}^{t_1} \sigma_{i_2}(t_2)
d{\bf w}_{t_2}^{(i_2)}dt_1,
$$
$$
\int\limits_{\tau}^t \sigma_{i_1}(t_1)\int\limits_{\tau}^{t_1} \sigma_{i_2}^2(t_2)
dt_2d{\bf w}_{t_1}^{(i_1)},\ \ \ 
\int\limits_{\tau}^t \sigma_0(t_1)\int\limits_{\tau}^{t_1} 
\int\limits_{\tau}^{t_{2}}\sigma_{i_3}(t_3)d{\bf w}_{t_3}^{(i_3)}
dt_2dt_1,
$$
$$
\int\limits_{\tau}^t \sigma_{i_1}(t_1)\int\limits_{\tau}^{t_1} \sigma_{i_2}(t_2)
\int\limits_{\tau}^{t_{2}}\sigma_{i_3}^2(t_3)dt_3 
d{\bf w}_{t_2}^{(i_2)}d{\bf w}_{t_1}^{(i_1)},
$$
$$
\int\limits_{\tau}^t \sigma_{i_1}(t_1)\int\limits_{\tau}^{t_1} \sigma_{i_2}(t_2)
\int\limits_{\tau}^{t_{2}}\sigma_{i_3}(t_3)d{\bf w}_{t_3}^{(i_3)}
d{\bf w}_{t_2}^{(i_2)}d{\bf w}_{t_1}^{(i_1)},
$$

\noindent
where $i_1,i_2,i_3=1,\ldots,m.$

Note that 
since 
$$
\sigma_0(t)\in L_1([0, T],{\bf R}),\ \ \ \sigma_j(t)\in L_2([0, T],{\bf R})\ \ \ (j=1,\ldots,m),
$$

\noindent
the iterated It\^{o} stochastic integrals (\ref{int100aaa})
exist in the mean-square sense and can be approximated using Theorem~1
(provided that $t$ and $\tau$ are constant time instants).

Denote
$$
\bar Q_{\lambda_l}^{(i_l,j_l)}=\left\{
\begin{matrix}
G^{(i_l)}\ \  & \hbox{for}\ \  i_l=1,\ldots,m\ \  \hbox{and}\ \   j_l=1\ \  (\hbox{if}\ \ \lambda_l=1)
\cr \cr 
L_{i_l}\ \ & \hbox{for}\ \ i_l=1,\ldots,m\ \ \hbox{and}\ \  j_l=0\ \ (\hbox{if}\ \ \lambda_l=1)
\cr \cr
L_0\ \ & \hbox{for}\ \ i_l=0\ \ \hbox{and}\ \  j_l=0\ \ (\hbox{if}\ \ \lambda_l=0) \cr\cr
\bar L_0\ \ & \hbox{for}\ \ i_l=0\ \ \hbox{and}\ \  j_l=-1\ \ (\hbox{if}\ \ \lambda_l=0) 
\end{matrix}\right.,
$$

\noindent
where $l=1,\ldots,k.$ 
Note that
\begin{equation}
\label{ww1003aa}
\bar L_0{\bf x}_{\tau}\equiv {\bf 0},\ \ \ 
L_0{\bf x}_{\tau}={\bf a}({\bf x}_{\tau},\tau),\ \ \ 
G^{(i_1)}{\bf x}_{\tau}=B_{i_1}({\bf x}_{\tau},\tau),\ \ \ L_{i_1}{\bf x}_{\tau}\equiv {\bf 0},
\end{equation}
where ${\bf 0}$ is a column vector of length $n$ consisting of zeros, and
$i_1=1,\ldots,m$.

Suppose that the function $R({\bf x},t): {\bf R}^n\times[0, T]\to {\bf R}$ is 
sufficiently smooth with respect 
to the variables ${\bf x}$ and $t.$ 
Consider the stochastic process $R({\bf x}_t,t),$ where ${\bf x}_t$
is the strong solution to the It\^{o} SDE (\ref{1600sqqqx}).

Iteratively 
applying It\^{o}'s formula (as above (see (\ref{5.1002})) to the process 
$R({\bf x}_t,t),$ we obtain the following 
Taylor--It\^{o} expansion
$$
R({\bf x}_t,t)=R({\bf x}_{\tau},\tau)
+
$$
$$
+\sum_{k=1}^r \sum_{(\lambda_{k},\ldots,\lambda_1)
\in{M}_k}
\sum_{i_1=\lambda_1}^{m\lambda_1}\sum_{j_1=\lambda_1-1}^{\lambda_1}
\ldots 
\sum_{i_k=\lambda_k}^{m\lambda_k}\sum_{j_k=\lambda_k-1}^{\lambda_k}
\bar Q_{\lambda_k}^{(i_k,j_k)}\ldots \bar Q_{\lambda_1}^{(i_1,j_1)}
R({\bf x}_{\tau},\tau)\times
$$
\begin{equation}
\label{9192sss}
\times 
\bar J_{(\lambda_k\ldots\lambda_1)t,\tau}^{~\bar \psi_{i_k,j_k}\ldots \bar \psi_{i_1,j_1}
((i_k,j_k)\ldots (i_1,j_1))}
+\left(\bar D_{r+1}\right)_{t,\tau}
\end{equation}

\vspace{2mm}
\noindent
w.~p.~1, where $t>\tau,$ $r\in {\bf N},$ and
$$
\left(\bar D_{r+1}\right)_{t,\tau}=
$$
$$
=
\sum_{(\lambda_{r+1},\ldots,\lambda_{1})\in{M}_{r+1}}
\sum_{i_1=\lambda_1}^{m\lambda_1}\sum_{j_1=\lambda_1-1}^{\lambda_1}
\ldots 
\sum_{i_{r+1}=\lambda_{r+1}}^{m\lambda_{r+1}}\sum_{j_{r+1}=\lambda_{r+1}-1}^{\lambda_{r+1}}
\int\limits_{\tau}^{t}\bar \psi_{i_1,j_1}(t_1)\int\limits_{\tau}^{t_1}\bar \psi_{i_2,j_2}(t_2)\ldots
$$
$$
\ldots
\int\limits_{\tau}^{t_r}
\bar \psi_{i_{r+1},j_{r+1}}(t_{r+1})\bar Q_{\lambda_{r+1}}^{(i_{r+1},j_{r+1})}\ldots
\bar Q_{\lambda_{2}}^{(i_2,j_2)}\bar Q_{\lambda_{1}}^{(i_1,j_1)}
R({\bf x}_{t_{r+1}})\times
$$
\begin{equation}
\label{9192asss}
\times
d\bar{\bf w}_{t_{r+1}}^{(i_{r+1},j_{r+1})}
\ldots
d\bar{\bf w}_{t_{2}}^{(i_2,j_2)}d\bar{\bf w}_{t_{1}}^{(i_1,j_1)}.
\end{equation}

\vspace{2mm}

Note that for $R({\bf x}_{\tau},\tau)={\bf x}_{\tau}$ the expression
$\bar Q_{\lambda_1}^{(i_1,j_1)}{\bf x}_{\tau}$ is defined by the equality (by virtue of (\ref{ww1003aa}))
\begin{equation}
\label{ww1004aaaa}
Q_{\lambda_1}^{(i_1,j_1)}{\bf x}_{\tau}=\left\{
\begin{matrix}
B_{i_1}({\bf x}_{\tau},\tau)\ \  & \hbox{for}\ \  i_1=1,\ldots,m\ \  \hbox{and}\ \   j_1=1\ \  (\hbox{if}\ \ 
\lambda_l=1)
\cr \cr 
{\bf 0}\ \ & \hbox{for}\ \ i_1=1,\ldots,m\ \ \hbox{and}\ \  j_1=0\ \ (\hbox{if}\ \ \lambda_1=1)
\cr \cr
{\bf a}({\bf x}_{\tau},\tau)\ \ & \hbox{for}\ \ i_1=0\ \ \hbox{and}\ \  j_1=0\ \ (\hbox{if}\ \ \lambda_1=0)
\cr \cr
{\bf 0}\ \ & \hbox{for}\ \ i_1=0\ \ \hbox{and}\ \  j_1=-1\ \ (\hbox{if}\ \ \lambda_1=0)
\end{matrix}\right..
\end{equation}

\vspace{2mm}

The issue regarding the existence of the remainder term 
(\ref{9192asss})
in the mean-square sense is resolved analogously to the case
of the remainder term $\left(D_{r+1}\right)_{t,\tau}$ of the form (\ref{9192a}).

Suppose that $R({\bf x},t)\equiv {\bf x}^{(i)},$ where ${\bf x}^{(i)}$
is the $i$th component of ${\bf x}$ $(i=1,\ldots,m).$
Then, using (\ref{ww1003aa}), (\ref{ww1004aaaa}), 
and the truncated Taylor--It\^{o} expansion (\ref{9192sss}) 
(without the remainder term $\left(\bar D_{r+1}\right)_{t,\tau}$), one can construct high-order 
stong numerical methods for It\^{o} SDE (\ref{1600sqqqx}), 
by analogy to usual It\^{o} SDEs
\cite{2},\cite{15}-\cite{15b}.

\section{Some Useful Estimates}

Consider the following It\^{o} SDE 
\begin{equation}
\label{1600t}
{\bf x}_t={\bf x}_0+\int\limits_0^t \sigma_0(s){\bf a}({\bf x}_{s},s)ds+
\sum\limits_{j=1}^m
\int\limits_0^t \sigma_j(s) B_j({\bf x}_{s},s)d{\bf w}_{s}^{(j)},\ \ \
{\bf x}_0={\bf x}(0,\omega),
\end{equation}
                                               
\noindent
where $t\in [0, T],$ 
$\sigma_0(t)\in L_1([0, T],{\bf R}),$ $\sigma_j(t)\in L_2([0, T],{\bf R})$ $(j=1,\ldots,m),$
the functions ${\bf a}({\bf x},t),$
$B_j({\bf x},t)$ $(j=1,\ldots,m):$ 
${\bf R}^n\times [0, T]\to{\bf R}^n$
are measurable with respect to the collection of 
variables $({\bf x},t),$ satisfy the global Lipschitz condition, and the linear growth condition;
another notations 
appearing in (\ref{1600sqqqx}) are the same as in (\ref{1.5.2}).

According to Theorem~1.1 from \cite{Jan}, under the above conditions, the It\^{o} SDE (\ref{1600t})
has a unique strong solution. Moreover, suppose that
\begin{equation}
\label{123xyz}
\sigma_0^r(t)\in L_1([0, T],{\bf R}),\ \ \ \sigma_j^{2r}(t)\in L_1([0, T],{\bf R})\ \ \ (r>2),
\end{equation}
where $j=1,\ldots,m$.
Then, from Theorem~1.7 \cite{Jan}, it follows that 
\begin{equation}
\label{667}
{\sf M}\left\{\left\Vert{\bf x}_t\right\Vert^r\right\}\le K<\infty\ \ \ (r>2),
\end{equation}
where $t\in [0, T],$ 
${\bf x}_t$ is a strong solution to
the It\^{o} SDE (\ref{1600t}),
the constant $K$ is independent of $t.$

Using (\ref{667}), let us prove that under the above conditions on the functions 
${\bf a}({\bf x},t),$
$B_j({\bf x},t)$ $(j=1,\ldots,m),$ as well as under the condition (\ref{123xyz}),
the following estimate holds for $r\ge 2$ 
\begin{equation}
\label{667xx}
{\sf M}\left\{\left\Vert{\bf x}_t-{\bf x}_s\right\Vert^r\right\}\le C\left(
\left(\int\limits_s^t \left\vert\sigma_0(\tau)\right\vert d\tau\right)^{r}+
\sum\limits_{j=1}^m \left(\int\limits_s^t \sigma_j^2(\tau)d\tau\right)^{r/2}\right),
\end{equation}

\noindent
where $C$ is a constants that does not depend on $t, s,$
$0\le s<t\le T.$

First, consider the case $r>2.$ Recall that 
${\bf x}^{(i)}$ is the $i$th component of
${\bf x}$, ${\bf a}^{(i)}({\bf x},t)$ is the $i$th component
of ${\bf a}({\bf x},t)$, $B^{(ij)}({\bf x},t)$
is the $ij$th element of the matrix $B({\bf x}, t)$
with columns $B_j({\bf x},t)$
$(i=1,\ldots,n,\ j=1,\ldots,m).$

From (\ref{1600t}) we have for $t>s$ w.~p.~1
\begin{equation}
\label{66811}
{\bf x}_t^{(i)}-{\bf x}_s^{(i)}=\int\limits_s^t \sigma_0(\tau){\bf a}^{(i)}({\bf x}_{\tau},\tau)d\tau+
\sum\limits_{j=1}^m
\int\limits_s^t \sigma_j(\tau) B^{(ij)}({\bf x}_{\tau},\tau)d{\bf w}_{\tau}^{(j)}.
\end{equation}

Then,
$$
{\sf M}\left\{\left|{\bf x}_t^{(i)}-{\bf x}_s^{(i)}\right|^r\right\}\le (m+1)^{r-1}\times
$$
\begin{equation}
\label{668}
\times
\left({\sf M}\left\{\left|\int\limits_s^t \sigma_0(\tau){\bf a}^{(i)}({\bf x}_{\tau},\tau)d\tau
\right|^r\right\}+
\sum\limits_{j=1}^m{\sf M}\left\{\left|
\int\limits_s^t \sigma_j(\tau) B^{(ij)}({\bf x}_{\tau},\tau)d{\bf w}_{\tau}^{(j)}
\right|^r\right\}\right),
\end{equation}
where $i=1,\ldots,m.$

In obtaining (\ref{668}), we used the inequality
\begin{equation}
\label{901}
\left|a_1+\ldots+a_k\right|^p\le k^{p-1}\left(\left|a_1\right|^p+\ldots+\left|a_k\right|^p\right),
\end{equation}
where $a_1,\ldots,a_k\in {\bf R},$ $k\in{\bf N},$ $p\ge 1.$

Using the Burkholder--Davis--Gundy inequality \cite{Karatzas}, we get
\begin{equation}
\label{681}
{\sf M}\left\{\left|\int\limits_s^t \sigma_j(\tau) B^{(ij)}({\bf x}_{\tau},\tau)d{\bf w}_{\tau}^{(j)}
\right|^r\right\}\le C_r
{\sf M}\left\{\left|\int\limits_s^t \sigma_j^2(\tau) \left(B^{(ij)}({\bf x}_{\tau},\tau)
\right)^2 d\tau \right|^{r/2}\right\},
\end{equation}

\noindent
where $C_r$ is a constant.

Further, we have for $r>2$ w.~p.~1
\begin{equation}
\label{669}
\sigma_j^2(\tau) \left(B^{(ij)}({\bf x}_{\tau},\tau)
\right)^2=\left(\sigma_j^2(\tau) \left|B^{(ij)}({\bf x}_{\tau},\tau)
\right|^r\right)^{2/r}\left(\sigma_j^2(\tau)\right)^{(r-2)/r}.
\end{equation}

Applying  the integral version of H\"older's inequality to (\ref{669}) with
$u=r/2,$ $v=r/(r-2)$ ($1/u+1/v=1,$ $r>2$), we obtain w.~p.~1
$$
\int\limits_s^t \sigma_j^2(\tau) \left(B^{(ij)}({\bf x}_{\tau},\tau)
\right)^2 d\tau \le 
\left(\int\limits_s^t 
\left(\left(\sigma_j^2(\tau) \left|B^{(ij)}({\bf x}_{\tau},\tau)
\right|^r\right)^{2/r}\right)^{r/2}d\tau\right)^{2/r}\times
$$
$$
\times
\left(\int\limits_s^t \left(\left(\sigma_j^2(\tau)\right)^{(r-2)/r}\right)^{r/(r-2)}d\tau\right)^{(r-2)/r}=
$$
$$
=\left(\int\limits_s^t 
\sigma_j^2(\tau) \left|B^{(ij)}({\bf x}_{\tau},\tau)
\right|^r d\tau\right)^{2/r}
\left(\int\limits_s^t \sigma_j^2(\tau) d\tau\right)^{(r-2)/r}.
$$

Then,
$$
\left|\int\limits_s^t \sigma_j^2(\tau) \left(B^{(ij)}({\bf x}_{\tau},\tau)
\right)^2 d\tau\right|^{r/2} \le 
$$
\begin{equation}
\label{670}
\le \int\limits_s^t 
\sigma_j^2(\tau) \left|B^{(ij)}({\bf x}_{\tau},\tau)
\right|^r d\tau
\left(\int\limits_s^t \sigma_j^2(\tau) d\tau\right)^{(r-2)/2}.
\end{equation}

Using Fubini's Theorem, the linear growth condition on $B^{(ij)}({\bf x},\tau),$
and the estimate (\ref{667}), from (\ref{670}) we have
$$
{\sf M}\left\{\left|\int\limits_s^t \sigma_j^2(\tau) \left(B^{(ij)}({\bf x}_{\tau},\tau)
\right)^2 d\tau\right|^{r/2}\right\} 
\le K_1 \int\limits_s^t 
\sigma_j^2(\tau) d\tau
\left(\int\limits_s^t \sigma_j^2(\tau) d\tau\right)^{(r-2)/2}=
$$
\begin{equation}
\label{671}
=
K_1\left(\int\limits_s^t \sigma_j^2(\tau) d\tau\right)^{r/2},
\end{equation}

\noindent
where the constant $K_1$ is independent of $s, t.$

Further, we have w.~p.~1
$$
\left|\int\limits_s^t \sigma_0(\tau) {\bf a}^{(i)}({\bf x}_{\tau},\tau)d\tau
\right|\le 
\int\limits_s^t \left|\sigma_0(\tau){\bf a}^{(i)}({\bf x}_{\tau},\tau)\right|d\tau=
$$
\begin{equation}
\label{672}
=\int\limits_s^t \left|\sigma_0(\tau)\right|^{(r-1)/r}
\left(\left|\sigma_0(\tau)\right|^{1/r}\left|{\bf a}^{(i)}({\bf x}_{\tau},\tau)\right|\right)d\tau.
\end{equation}

Applying  the integral version of H\"older's inequality to (\ref{672}) with
$u=r/(r-1),$ $v=r$ ($1/u+1/v=1$), we obtain w.~p.~1
$$
\int\limits_s^t \left|\sigma_0(\tau){\bf a}^{(i)}({\bf x}_{\tau},\tau)\right|d\tau\le
\left(\int\limits_s^t
\left(\left|\sigma_0(\tau)\right|^{(r-1)/r}\right)^{r/(r-1)}d\tau\right)^{(r-1)/r}\times
$$
$$
\times \left(\int\limits_s^t
\left(\left|\sigma_0(\tau)\right|^{1/r}\left|{\bf a}^{(i)}({\bf x}_{\tau},\tau)\right|\right)^r d\tau
\right)^{1/r}=
$$
\begin{equation}
\label{673}
=\left(\int\limits_s^t
\left|\sigma_0(\tau)\right|d\tau\right)^{(r-1)/r}
\left(\int\limits_s^t
\left|\sigma_0(\tau)\right|\left|{\bf a}^{(i)}({\bf x}_{\tau},\tau)\right|^r d\tau
\right)^{1/r}.
\end{equation}

Using (\ref{672}) and (\ref{673}), we get w.~p.~1
\begin{equation}
\label{674}
\left|\int\limits_s^t \sigma_0(\tau) {\bf a}^{(i)}({\bf x}_{\tau},\tau)d\tau
\right|^r\le 
\left(\int\limits_s^t
\left|\sigma_0(\tau)\right|d\tau\right)^{r-1}
\int\limits_s^t
\left|\sigma_0(\tau)\right|\left|{\bf a}^{(i)}({\bf x}_{\tau},\tau)\right|^r d\tau.
\end{equation}

Applying Fubini's Theorem, the linear growth condition on ${\bf a}^{(i)}({\bf x},\tau),$
and the estimate (\ref{667}), from (\ref{674}) we have
$$
{\sf M}\left\{\left|\int\limits_s^t \sigma_0(\tau) {\bf a}^{(i)}({\bf x}_{\tau},\tau)d\tau
\right|^r\right\}\le K_2
\left(\int\limits_s^t
\left|\sigma_0(\tau)\right|d\tau\right)^{r-1}
\int\limits_s^t
\left|\sigma_0(\tau)\right|d\tau=
$$
\begin{equation}
\label{675}
=K_2\left(\int\limits_s^t
\left|\sigma_0(\tau)\right|d\tau\right)^{r},
\end{equation}

\noindent
where the constant $K_2$ is independent of $s, t.$

Combining (\ref{668}), (\ref{681}), (\ref{671}), and (\ref{675}), we obtain 
\begin{equation}
\label{676}
{\sf M}\left\{\left|{\bf x}_t^{(i)}-{\bf x}_s^{(i)}\right|^r\right\}\le C_1\left(
\left(\int\limits_s^t \left\vert\sigma_0(\tau)\right\vert d\tau\right)^{r}+
\sum\limits_{j=1}^m \left(\int\limits_s^t \sigma_j^2(\tau)d\tau\right)^{r/2}\right),
\end{equation}

\noindent
where $i=1,\ldots,m,$ and the constant $C_1$ is independent of $s, t.$

Applying the inequality (\ref{901}), we get
\begin{equation}
\label{902}
{\sf M}\left\{\left\Vert{\bf x}_t-{\bf x}_s\right\Vert^r\right\}\le
n^{r/2-1}\sum\limits_{i=1}^n {\sf M}\left\{\left|{\bf x}_t^{(i)}-{\bf x}_s^{(i)}\right|^r\right\}.
\end{equation}

From (\ref{676}), (\ref{902}) we have
(\ref{667xx}). The estimate (\ref{667xx}) is proved
for $r>2$.

It is easy to see that for $r=2,$ we can apply the Cauchy--Bunyakovsky
inequality. Indeed, from (\ref{66811}) we have
$$
{\sf M}\left\{\left|{\bf x}_t^{(i)}-{\bf x}_s^{(i)}\right|^2\right\}\le (m+1)\times
$$
$$
\le\hspace{-0.5mm}
\left({\sf M}\left\{\left|\int\limits_s^t \sigma_0(\tau){\bf a}^{(i)}({\bf x}_{\tau},\tau)d\tau
\right|^2\right\}+
\sum\limits_{j=1}^m{\sf M}\left\{\left|
\int\limits_s^t \sigma_j(\tau) B^{(ij)}({\bf x}_{\tau},\tau)d{\bf w}_{\tau}^{(j)}
\right|^2\right\}\right)\hspace{-0.5mm}\le
$$
$$
\le (m+1)
\left(\int\limits_s^t \left|\sigma_0(\tau)\right|d\tau\
{\sf M}\left\{\int\limits_s^t \left|\sigma_0(\tau)\right|\left(
{\bf a}^{(i)}({\bf x}_{\tau},\tau)\right)^2d\tau\right\}+\right.
$$
$$
\left.
+
\sum\limits_{j=1}^m{\sf M}\left\{
\int\limits_s^t \sigma_j^2(\tau) \left(B^{(ij)}({\bf x}_{\tau},\tau)\right)^2 d\tau
\right\}\right).
$$

The further proof for the case $r=2$ can be carried out analogously
to the case $r>2.$ Thus, the estimate 
(\ref{667xx}) also holds for $r=2.$

It should be noted that the estimates (\ref{667}), (\ref{667xx}) are necessary to 
justify the strong convergence rate of the new modification (\ref{11}) of 
the Milstein method, modified Euler method (numerical method (\ref{11}) 
without the last term), and more accurate numerical methods 
(constructed on the basis of the Taylor--It\^{o} expansion (\ref{9192sss})) 
for It\^{o} SDE (\ref{1600t}) with singular coefficients.

Suppose that 
$\sigma_0(s)=s^{\alpha},$
$\sigma_j(s)=t^{\beta_j},$ $\alpha\in (-1,0)$, $\beta_j\in (-1/2,0)$ $(j=1,\ldots,m)$,
i.e. $\sigma_0(s)=s^{-|\alpha|},$ $\sigma_j(s)=t^{-|\beta_j|},$
$|\alpha|\in (0,1)$, $|\beta_j|\in (0,1/2)$ $(j=1,\ldots,m)$.

Let us estimate the right-hand side of (\ref{667xx}) for such functions 
$\sigma_0(s),$
$\sigma_j(s)$ $(j=1,\ldots,m)$.

We have 
\begin{equation}
\label{1980}
\int\limits_s^t \left\vert\sigma_0(\tau)\right\vert d\tau=
\int\limits_s^t {\tau}^{-|\alpha|} d\tau=
\frac{t^{1-|\alpha|}-s^{1-|\alpha|}}{1-|\alpha|}\le \frac{(t-s)^{1-|\alpha|}}{1-|\alpha|},
\end{equation}
where $0\le s<t\le T.$ 

In obtaining (\ref{1980}), we used the inequality
$$
x^{\gamma}-1\le (x-1)^{\gamma}, \ \ \ x\ge 1,\ \ \ \gamma\in (0,1).
$$

Analogously, we obtain 
\begin{equation}
\label{1981}
\int\limits_s^t \sigma_j^2(\tau)d\tau=\int\limits_s^t {\tau}^{-2|\beta_j|}d\tau\le
\frac{t^{1-2|\beta_j|}-s^{1-2|\beta_j|}}{1-2|\beta_j|}\le \frac{(t-s)^{1-2|\beta_j|}}{1-2|\beta_j|}.
\end{equation}
            
Combining (\ref{667xx}), (\ref{1980}), (\ref{1981}), we get
\begin{equation}
\label{1982}
{\sf M}\bigl\{\left\Vert{\bf x}_t-{\bf x}_s\right\Vert^2\bigr\}\le C
\left(\frac{(t-s)^{2-2|\alpha|}}{(1-|\alpha|)^2}+
\sum\limits_{j=1}^m \frac{(t-s)^{1-2|\beta_j|}}{1-2|\beta_j|}\right),
\end{equation}
where $C$ is a constants that does not depend on $t, s.$

Consider a partition $\{\tau_i\}_{i=0}^N$ 
of the interval $[0, T]$ such that
\begin{equation}
\label{1984}
0=\tau_0<\tau_1<\ldots <\tau_N=T,\ \ \
\Delta_N=
\hbox{\vtop{\offinterlineskip\halign{
\hfil#\hfil\cr
{\rm max}\cr
$\stackrel{}{{}_{0\le i\le N-1}}$\cr
}} }\left|\tau_{i+1}-\tau_i\right|.
\end{equation}

Let $[s,t]\subseteq [\tau_i,\tau_{i+1}]$ for some $i$ $(0\le i\le N-1),$ i.e.
$|t-s|\le \tau_{i+1}-\tau_i \le \Delta_N.$ Then, from 
(\ref{1982}) we obtain
$$
{\sf M}\bigl\{\left\Vert{\bf x}_t-{\bf x}_s\right\Vert^2\bigr\}\le C
\left(\frac{(\Delta_N)^{2-2|\alpha|}}{(1-|\alpha|)^2}+
\sum\limits_{j=1}^m \frac{(\Delta_N)^{1-2|\beta_j|}}{1-2|\beta_j|}\right)\le
$$
$$
\le \bar C (\Delta_N)^{\min\limits_{1\le j\le m}\{2-2|\alpha|,\ 1-2|\beta_j|\}}=
\bar C (\Delta_N)^{\min\{2-2|\alpha|,\ 1-2\max\limits_{1\le j\le m}|\beta_j|\}}=
$$
\begin{equation}
\label{1983}
=\bar C (\Delta_N)^{2\min\{1-|\alpha|,\ 1/2-\max\limits_{1\le j\le m}|\beta_j|\}},
\end{equation}

\noindent
where $\bar C$ is a constants that does not depend on $\Delta_N.$

The estimate (\ref{1983}) will be used below.

\section{Convergence Rate of the Modified Euler Method for the It\^{o} SDE (\ref{1600t})}

Consider the It\^{o} SDE (\ref{1600t}).
We call the numerical scheme (\ref{11}) without the last term the modified Euler 
method, i.e, 
\begin{equation}
\label{euler1}
{\bf y}_{p+1}={\bf y}_{p}+
{\bf a}({\bf y}_{p},\tau_p) \int\limits_{\tau_p}^{\tau_{p+1}} \sigma_0(s)ds+
\sum\limits_{j=1}^m
B_j({\bf y}_{p},\tau_p)\int\limits_{\tau_p}^{\tau_{p+1}} \sigma_j(s)d{\bf w}_{s}^{(j)},
\end{equation}
where $p=0,1,\ldots,N-1,$ the partition $\{\tau_p\}_{p=0}^N$ 
of the interval $[0, T]$ satisfies the condition (\ref{1984}),
${\bf y}_{\tau_p}\stackrel{\sf def}{=}{\bf y}_p,$ 
another notations in (\ref{euler1}) are the same as in
(\ref{1600sqqqx}).

Let us estimate the convergence rate of the modified Euler method (\ref{euler1}).

{\bf Theorem 2}.\ {\it Let $\sigma_0(t)=t^{-|\alpha|},$ $\sigma_j(t)=t^{-|\beta_j|},$ 
$\alpha\in (-1/(2+\varepsilon), 0),$
$\beta_j\in (-1/(4+2\varepsilon), 0),$
where $\varepsilon$ is an arbitrary
small positive real number. Furthermore, 
let the following conditions be satisfied 
for all ${\bf x}, {\bf y}\in{\bf R}^n,$ $s, t\in[0,T]:$
\begin{equation}
\label{888.06}
\left({\sf M}\bigl\{\left\Vert{\bf x}_0-{\bf y}_0\right\Vert^2\bigr\}\right)^{1/2}\le
M\sqrt{\Delta_N},
\end{equation}
$$
\left\Vert{\bf a}({\bf x},t)-{\bf a}({\bf y},t)\right\Vert+\sum_{j=1}^m
\left\Vert B_j({\bf x},t)-B_j({\bf y},t)\right\Vert\le C\left\Vert{\bf x}-{\bf y}\right\Vert,
$$
$$
\left\Vert{\bf a}({\bf x},t)\right\Vert+\sum_{j=1}^m
\left\Vert B_j({\bf x},t)\right\Vert\le C(1+\left\Vert{\bf x}\right\Vert),
$$
$$
\left\Vert{\bf a}({\bf x},s)-{\bf a}({\bf x},t)\right\Vert+\sum_{j=1}^m
\left\Vert B_j({\bf x},s)-B_j({\bf x},t)\right\Vert\le 
C\sqrt{s-t}(1+\left\Vert{\bf x}\right\Vert),
$$
where $C, M$ are constants that are independent of 
$\Delta_N.$ 
Then
\begin{equation}
\label{res0}
{\sf M}\{\left\Vert{\bf x}_T-{\bf y}_N\right\Vert\}\le K
(\Delta_N)^{\min\{1-|\alpha|,\ 1/2-\max\limits_{1\le j\le m}|\beta_j|\}},
\end{equation}
where the constant $K$ is independent of $\Delta_N,$ 
and ${\bf x}_T$ is the solution of the It\^{o} SDE {\rm (\ref{1600t})} at $t=T.$}

{\bf Proof.}\ Under the assumptions of this theorem, 
according to Theorem 1.7 \cite{Jan} and the results 
of the previous section, we have 
\begin{equation}
\label{useful1}
{\sf M}\bigl\{\left\Vert{\bf x}_s\right\Vert^2\bigr\}\le 
K_1<\infty,
\end{equation}
\begin{equation}
\label{1986}
{\sf M}\bigl\{\left\Vert{\bf x}_t-{\bf x}_s\right\Vert^2\bigr\}
\le \bar C (\Delta_N)^{2\min\{1-|\alpha|,\ 1/2-\max\limits_{1\le j\le m}|\beta_j|\}},
\end{equation}
where $[s,t]\subseteq [\tau_p, \tau_{p+1}]$ for some $p,$ 
the constants $K_1, \bar C$ do not depend on $\Delta_N.$

Denote
$$
{\bf y}_t={\bf y}_p+\int\limits_{\tau_p}^t \sigma_0(s){\bf a}({\bf y}_p,\tau_p)ds+
\sum\limits_{j=1}^m \int\limits_{\tau_p}^t \sigma_j(s)B_j({\bf y}_p,\tau_p)d{\bf w}_s^{(j)},\ \ \
t\in[\tau_p,\tau_{p+1}].
$$

Obviously, the stochastic process ${\bf y}_t$ w.~p.~1 coincides 
with ${\bf y}_p$ at $t=\tau_p,$ $p=0, 1,\ldots,N.$

It is easy to see that
\begin{equation}
\label{888.07}
{\bf y}_t-{\bf x}_t={\bf y}_0-{\bf x}_0+U_t+V_t+R_t+G_t+Q_t+P_t,
\end{equation}
where ${\bf x}_t$ is a strong solution to the It\^{o} SDE (\ref{1600t}),
$$
U_t=\sum_{l=0}^{n_t-1}\int\limits_{\tau_l}^{\tau_{l+1}}
\sigma_0(s)\left({\bf a}({\bf x}_s,\tau_l)-{\bf a}({\bf x}_s,s)\right)
ds+
\int\limits_{\tau_{n_t}}^{t}
\sigma_0(s)\left({\bf a}({\bf x}_s,\tau_{n_t})-{\bf a}({\bf x}_s,s)\right)
ds,
$$
$$
V_t=
\sum\limits_{j=1}^m\Biggl(\sum_{l=0}^{n_t-1}\int\limits_{\tau_l}^{\tau_{l+1}}
\sigma_j(s)\left(B_j({\bf x}_s,\tau_l)-B_j({\bf x}_s,s)\right)d{\bf w}_s^{(j)}+\Biggr.
$$
$$
\Biggl.+
\int\limits_{\tau_{n_t}}^{t}
\sigma_j(s)\left(B_j({\bf x}_s,\tau_{n_t})-B_j({\bf x}_s,s)\right)
d{\bf w}_s^{(j)}\Biggr),
$$
$$
R_t=\sum_{l=0}^{n_t-1}\int\limits_{\tau_l}^{\tau_{l+1}}
\sigma_0(s)\left({\bf a}({\bf x}_{\tau_l},\tau_l)-{\bf a}({\bf x}_s,{\tau_l})\right)
ds+
$$
\begin{equation}
\label{start5}
+
\int\limits_{\tau_{n_t}}^{t}
\sigma_0(s)\left({\bf a}({\bf x}
_{\tau_{n_t}},\tau_{n_t})-{\bf a}({\bf x}_s,\tau_{n_t})\right)ds,
\end{equation}
$$
G_t=\sum\limits_{j=1}^m\Biggl(\sum_{l=0}^{n_t-1} \int\limits_{\tau_l}^{\tau_{l+1}}
\sigma_j(s)\left(B_j({\bf x}_{\tau_l},\tau_l)-B_j({\bf x}_s,
\tau_l)\right)d{\bf w}_s^{(j)}+\Biggr.
$$
\begin{equation}
\label{start6}
\Biggl.+
\int\limits_{\tau_{n_t}}^{t}
\sigma_j(s)\left(B_j({\bf x}_{\tau_{n_t}},\tau_{n_t})-B_j({\bf x}_s,
\tau_{n_t})\right)
d{\bf w}_s^{(j)}\Biggr),
\end{equation}
$$
Q_t=\sum_{l=0}^{n_t-1}\int\limits_{\tau_l}^{\tau_{l+1}}
\sigma_0(s)\left({\bf a}({\bf y}_l,\tau_l)-{\bf a}({\bf x}_{\tau_l},\tau_l)\right)
ds+
$$
$$
+
\int\limits_{\tau_{n_t}}^{t}
\sigma_0(s)\left({\bf a}({\bf y}_{n_t},\tau_{n_t})-{\bf a}({\bf x}_{\tau_{n_t}},\tau_{n_t})\right)
ds,
$$
$$
P_t=\sum\limits_{j=1}^m\Biggl(\sum_{l=0}^{n_t-1}\int\limits_{\tau_l}^{\tau_{l+1}}
\sigma_j(s)\left(B_j({\bf y}_l,\tau_l)-B_j({\bf x}_{\tau_l},\tau_l)\right)
d{\bf w}_s^{(j)}+\Biggr.
$$
$$
\Biggl.+
\int\limits_{\tau_{n_t}}^{t}
\sigma_j(s)\left(B_j({\bf y}_{n_t},\tau_{n_t})-B_j({\bf x}_{\tau_{n_t}},\tau_{n_t})\right)
d{\bf w}_s^{(j)}\Biggr),
$$
$$
n_t=
\max\limits_{j}\{j:\ \tau_j\le t\}.
$$

From (\ref{888.07}), we obtain for $u\in [0, T]$
$$
Z_u\stackrel{\sf def}{=}\sup\limits_{0\le t\le u}
{\sf M}\bigl\{\left\Vert{\bf y}_t-{\bf x}_t\right\Vert ^2\bigr\}\le 
$$
$$
\le
\sup\limits_{0\le t\le u}
K_2\left({\sf M}\bigl\{\left\Vert{\bf y}_0-{\bf x}_0\right\Vert^2\bigr\}+
{\sf M}\bigl\{\left\Vert U_t\right\Vert^2\bigr\}+
{\sf M}\bigl\{\left\Vert V_t\right\Vert^2\bigr\} +\right.
$$
\begin{equation}
\label{888.14}
\left.+
{\sf M}\bigl\{\left\Vert R_t\right\Vert^2\bigr\}+{\sf M}\bigl\{\left\Vert G_t\right\Vert^2\bigr\}+
{\sf M}\bigl\{\left\Vert Q_t\right\Vert^2\bigr\}+
{\sf M}\bigl\{\left\Vert P_t\right\Vert^2\bigr\}\right),
\end{equation}
where $K_2$ is a constant.

Let us estimate one by one
the terms on the right-hand side of (\ref{888.14}).
Using Minkowski's inequality for mathematical expectations, 
the Cauchy--Bunyakovsky inequality, the Fubini Theorem,
the conditions of Theorem~2, and (\ref{useful1}), we have

\vspace{-2mm}
$$
\left({\sf M}\left\{\bigl|U_t^{(k)}\bigr|^2\right\}\right)^{1/2}\le
$$

\vspace{-2mm}
$$
\le
\sum_{l=0}^{n_t-1}\Biggl({\sf M}\Biggl\{\Biggl(
\int\limits_{\tau_l}^{\tau_{l+1}}
\sigma_0(s)\left({\bf a}^{(k)}({\bf x}_s,\tau_l)-
{\bf a}^{(k)}({\bf x}_s,s)\right)ds\Biggr)^2\Biggr\}
\Biggr)^{1/2}+
$$
$$
+
\Biggl({\sf M}\Biggl\{\Biggl(~\int\limits_{\tau_{n_t}}^{t}
\sigma_0(s)\left({\bf a}^{(k)}({\bf x}_s,\tau_{n_t})-
{\bf a}^{(k)}({\bf x}_s,s)\right)ds\Biggr)^2\Biggr\}
\Biggr)^{1/2}\le
$$
$$
\le
\sum_{l=0}^{n_t-1}\Biggl(\int\limits_{\tau_l}^{\tau_{l+1}}
\left|\sigma_0(s)\right|ds
\int\limits_{\tau_l}^{\tau_{l+1}}
\left|\sigma_0(s)\right|{\sf M}\left\{\left({\bf a}^{(k)}({\bf x}_s,\tau_l)-
{\bf a}^{(k)}({\bf x}_s,s)\right)^2\right\}
ds
\Biggr)^{1/2}+
$$
$$
+
\Biggl(\int\limits_{\tau_{n_t}}^{t} \left|\sigma_0(s)\right|ds
\int\limits_{\tau_{n_t}}^{t} \left|\sigma_0(s)\right|{\sf M}\left\{
\left({\bf a}^{(k)}({\bf x}_s,\tau_{n_t})-{\bf a}^{(k)}({\bf x}_s,
s)\right)^2\right\}
ds\Biggr)^{1/2}\le
$$
$$
\le C_1\sqrt{\Delta_N}
\sum_{l=0}^{n_t-1}\Biggl(\int\limits_{\tau_l}^{\tau_{l+1}}
\left|\sigma_0(s)\right|ds
\int\limits_{\tau_l}^{\tau_{l+1}}
\left|\sigma_0(s)\right|{\sf M}\left\{(1+\left\Vert{\bf x}_s\right\Vert)^2\right\}
ds
\Biggr)^{1/2}+
$$
$$
+
\Biggl(\int\limits_{\tau_{n_t}}^{t} \left|\sigma_0(s)\right|ds
\int\limits_{\tau_{n_t}}^{t} \left|\sigma_0(s)\right|
{\sf M}\left\{(1+\left\Vert{\bf x}_s\right\Vert)^2\right\}
ds\Biggr)^{1/2}\le
$$
$$
\le \sqrt{\Delta_N} C_2 \Biggl(
\sum_{l=0}^{n_t-1}\int\limits_{\tau_l}^{\tau_{l+1}}
\left|\sigma_0(s)\right|ds+
\int\limits_{\tau_{n_t}}^{t} \left|\sigma_0(s)\right|ds\Biggr)=
\sqrt{\Delta_N} C_2 
\int\limits_{0}^{t}
\left|\sigma_0(s)\right|ds\le
$$
\begin{equation}
\label{888.15}
\le C_3\sqrt{\Delta_N},
\end{equation}

\noindent
where here and below, $C_i$ $(i=1,2,\ldots)$ will denote constants that are independent of $\Delta_N$
and $X^{(k)}$ will denote the $k$th component of the vector $X.$

Applying Minkowski's inequality for mathematical expectations, we obtain
$$
\left({\sf M}\left\{\bigl|V_t^{(k)}\bigr|^2\right\}\right)^{1/2}\le
$$
$$
\le
\sum_{j=1}^m
\Biggl({\sf M}\Biggl\{\Biggl(\sum_{l=0}^{n_t-1}
\int\limits_{\tau_l}^{\tau_{l+1}}
\sigma_j(s)\left(B^{(kj)}({\bf x}_s,\tau_l)-B^{(kj)}({\bf x}_s,s)\right)
d{\bf w}_s^{(j)}+\Biggr.\Biggr.\Biggr.
$$
\begin{equation}
\label{useful3}
\Biggl.\Biggl.\Biggl.+\int\limits_{\tau_{n_t}}^{t}
\sigma_j(s)\left(B^{(kj)}({\bf x}_s,\tau_{n_t})-B^{(kj)}({\bf x}_s,s)\right)
d{\bf w}_s^{(j)}\Biggr)^2\Biggr\}
\Biggr)^{1/2},
\end{equation}

\noindent
where $B^{(kj)}({\bf x},t)$ is the $kj$th element of the matrix $B({\bf x},t).$

Let us note that for It\^{o} stochastic integrals, one can write
$$
{\sf M}\left\{\int\limits_{t_3}^{t_4}\xi_{\tau}dw_{\tau}
\int\limits_{t_1}^{t_2}\eta_{\tau}dw_{\tau}\right\}=
{\sf M}\left\{{\sf M}\left\{\int\limits_{t_3}^{t_4}\xi_{\tau}dw_{\tau}
\int\limits_{t_1}^{t_2}\eta_{\tau}dw_{\tau}\bigl|\bigr.{\rm F}_{t_3}\right\}\right\}=
$$
\begin{equation}
\label{useful2}
={\sf M}\left\{\int\limits_{t_1}^{t_2}\eta_{\tau}dw_{\tau}
{\sf M}\left\{\int\limits_{t_3}^{t_4}\xi_{\tau}dw_{\tau}
\bigl|\bigr.{\rm F}_{t_3}\right\}\right\}=0\ \ \ \hbox{w.~p.~1,}
\end{equation}

\noindent
where $0\le t_1<t_2\le t_3<t_4\le T,$ $\xi_{\tau}$ and $\eta_{\tau}$ are some stochastic 
processes for which there exist It\^{o} stochastic integrals 
in the mean-square sense  that appear in the formula (\ref{useful2}), ${\rm F}_{t_3}\subset {\rm F},$
and $w_{\tau}$ is a standard scalar Wiener process.

Then, using the conditions of Theorem~2, It\^{o}'s isometry, (\ref{useful3}), 
(\ref{useful2}), and (\ref{useful1}), we obtain
$$
\left({\sf M}\left\{\bigl|V_t^{(k)}\bigr|^2\right\}\right)^{1/2}
\le
$$
$$
\le\sum_{j=1}^m
\Biggl(\sum_{l=0}^{n_t-1}
\int\limits_{\tau_l}^{\tau_{l+1}}
\sigma_j^2(s){\sf M}\left\{\left(B^{(kj)}({\bf x}_s,\tau_l)-
B^{(kj)}({\bf x}_s,s)\right)^2\right\}
ds+\Biggr.
$$
$$
\Biggl.+
\int\limits_{\tau_{n_t}}^{t}\sigma_j^2(s){\sf M}\left\{
\left(B^{(kj)}({\bf x}_s,\tau_{n_t})-
B^{(kj)}({\bf x}_s,s)\right)^2\right\}
ds\Biggr)^{1/2} \le
$$
$$
\le C_4\sum_{i=1}^m\Biggl(
\sum_{l=0}^{n_t-1}
\int\limits_{\tau_l}^{\tau_{l+1}}
\sigma_j^2(s){\sf M}\left\{(1+\left\Vert{\bf x}_s\right\Vert)^2\right\}(\tau_{l+1}-\tau_l)
ds
+\Biggr.
$$
$$
\Biggl.
+\int\limits_{\tau_{n_t}}^{t}\sigma_j^2(s)
{\sf M}\left\{(1+\left\Vert{\bf x}_s\right\Vert)^2\right\}
(t-\tau_{n_t})
ds\Biggr)^{1/2}\le
$$
$$
\le C_5\sqrt{\Delta_N}\sum_{j=1}^m\Biggl(
\sum_{l=0}^{n_t-1}
\int\limits_{\tau_l}^{\tau_{l+1}}
\sigma_j^2(s)ds
+\int\limits_{\tau_{n_t}}^{t}\sigma_j^2(s)ds\Biggr)^{1/2}=
$$
\begin{equation}
\label{useful11}
=
C_5\sqrt{\Delta_N}\sum_{j=1}^m
\Biggl(\int\limits_{0}^{t}
\sigma_j^2(s)ds\Biggr)^{1/2}
\le C_6\sqrt{\Delta_N}.
\end{equation}

Similarly to the derivation of (\ref{888.15}), and 
using the conditions of Theorem~2 and (\ref{1986}), we have
$$
\left({\sf M}\left\{\bigl|R_t^{(k)}\bigr|^2\right\}\right)^{1/2}\le
$$
$$
\le
\sum_{l=0}^{n_t-1}\Biggl({\sf M}\Biggl\{\Biggl(
\int\limits_{\tau_l}^{\tau_{l+1}}
\sigma_0(s)\left({\bf a}^{(k)}({\bf x}_{\tau_l},\tau_l)-
{\bf a}^{(k)}({\bf x}_s,\tau_l)\right)ds\Biggr)^2\Biggr\}
\Biggr)^{1/2}+
$$
$$
+
\Biggl({\sf M}\Biggl\{\Biggl(~\int\limits_{\tau_{n_t}}^{t}
\sigma_0(s)\left({\bf a}^{(k)}({\bf x}_{\tau_{n_t}},\tau_{n_t})-
{\bf a}^{(k)}({\bf x}_s,\tau_{n_t})\right)ds\Biggr)^2\Biggr\}
\Biggr)^{1/2}\le
$$
$$
\le
\sum_{l=0}^{n_t-1}\Biggl(\int\limits_{\tau_l}^{\tau_{l+1}}
\left|\sigma_0(s)\right|ds
\int\limits_{\tau_l}^{\tau_{l+1}}
\left|\sigma_0(s)\right|{\sf M}\left\{\left({\bf a}^{(k)}({\bf x}_{\tau_l},\tau_l)-
{\bf a}^{(k)}({\bf x}_s,\tau_l)\right)^2\right\}
ds
\Biggr)^{1/2}+
$$
$$
+
\Biggl(\int\limits_{\tau_{n_t}}^{t} \left|\sigma_0(s)\right|ds
\int\limits_{\tau_{n_t}}^{t} \left|\sigma_0(s)\right|{\sf M}\left\{
\left({\bf a}^{(k)}({\bf x}_{\tau_{n_t}},\tau_{n_t})-{\bf a}^{(k)}({\bf x}_s,
\tau_{n_t})\right)^2\right\}
ds\Biggr)^{1/2}\le
$$
$$
\le
C_7\Biggl(\sum_{l=0}^{n_t-1}\Biggl(\int\limits_{\tau_l}^{\tau_{l+1}}
\left|\sigma_0(s)\right|ds
\int\limits_{\tau_l}^{\tau_{l+1}}
\left|\sigma_0(s)\right|{\sf M}\left\{\left\Vert{\bf x}_{\tau_l}-{\bf x}_s\right\Vert^2\right\}
ds
\Biggr)^{1/2}+\Biggr.
$$
\begin{equation}
\label{setka15}
\Biggl.+
\Biggl(\int\limits_{\tau_{n_t}}^{t} \left|\sigma_0(s)\right|ds
\int\limits_{\tau_{n_t}}^{t} \left|\sigma_0(s)\right|
{\sf M}\left\{\left\Vert{\bf x}_{\tau_{n_t}}-{\bf x}_s\right\Vert^2\right\}
ds\Biggr)^{1/2}\Biggr)\le
\end{equation}
$$
\le
C_8 (\Delta_N)^{\min\{1-|\alpha|,\ 1/2-\max\limits_{1\le j\le m}|\beta_j|\}}
\int\limits_{0}^{t}\left|\sigma_0(s)\right|ds\le 
$$
\begin{equation}
\label{comb1}
\le
C_9 (\Delta_N)^{\min\{1-|\alpha|,\ 1/2-\max\limits_{1\le j\le m}|\beta_j|\}}.
\end{equation}

\vspace{3mm}

Furthermore, similarly to the derivation of (\ref{useful3}), (\ref{useful11}), 
and using (\ref{1986}) and the conditions of Theorem~2, 
we obtain

$$
\left({\sf M}\left\{\bigl|G_t^{(k)}\bigr|^2\right\}\right)^{1/2}\le
$$

\vspace{-4mm}
$$
\le
\sum_{j=1}^m
\Biggl({\sf M}\Biggl\{\Biggl(\sum_{l=0}^{n_t-1}
\int\limits_{\tau_l}^{\tau_{l+1}}
\sigma_j(s)\left(B^{(kj)}({\bf x}_{\tau_l},\tau_l)-B^{(kj)}({\bf x}_s,\tau_l)\right)
d{\bf w}_s^{(j)}+\Biggr.\Biggr.\Biggr.
$$
$$
\Biggl.\Biggl.\Biggl.+\int\limits_{\tau_{n_t}}^{t}
\sigma_j(s)\left(B^{(kj)}({\bf x}_{\tau_{n_t}},\tau_{n_t})-B^{(kj)}({\bf x}_s,\tau_{n_t})\right)
d{\bf w}_s^{(j)}\Biggr)^2\Biggr\}
\Biggr)^{1/2}=
$$
$$
=\sum_{j=1}^m
\Biggl(\sum_{l=0}^{n_t-1}
\int\limits_{\tau_l}^{\tau_{l+1}}
\sigma_j^2(s){\sf M}\left\{\left(B^{(kj)}({\bf x}_{\tau_l},\tau_l)-
B^{(kj)}({\bf x}_s,\tau_l)\right)^2\right\}
ds+\Biggr.
$$
$$
\Biggl.+
\int\limits_{\tau_{n_t}}^{t}\sigma_j^2(s){\sf M}\left\{
\left(B^{(kj)}({\bf x}_{\tau_{n_t}},\tau_{n_t})-
B^{(kj)}({\bf x}_s,\tau_{n_t})\right)^2\right\}
ds\Biggr)^{1/2} \le
$$
$$
\le C_{10}\sum_{j=1}^m
\Biggl(\sum_{l=0}^{n_t-1}
\int\limits_{\tau_l}^{\tau_{l+1}}
\sigma_j^2(s)
{\sf M}\left\{\left\Vert{\bf x}_{\tau_l}-{\bf x}_{s}\right\Vert^2\right\}
ds+\Biggr.
$$
\begin{equation}
\label{setka30}
\Biggl.+
\int\limits_{\tau_{n_t}}^{t}\sigma_j^2(s)
{\sf M}\left\{\left\Vert{\bf x}_{\tau_{n_t}}-{\bf x}_{s}\right\Vert^2\right\}
ds\Biggr)^{1/2} \le
\end{equation}

\vspace{1mm}
$$
\le
C_{11}(\Delta_N)^{\min\{1-|\alpha|,\ 1/2-\max\limits_{1\le j\le m}|\beta_j|\}}
\times
$$

\vspace{-1mm}
$$
\times
\sum_{j=1}^m\Biggl(\sum_{l=0}^{n_t-1}
\int\limits_{\tau_l}^{\tau_{l+1}}
\sigma_j^2(s)
ds+
\int\limits_{\tau_{n_t}}^{t}\sigma_j^2(s)
ds\Biggr)^{1/2}=
$$
$$
=
C_{11}(\Delta_N)^{\min\{1-|\alpha|,\ 1/2-\max\limits_{1\le j\le m}|\beta_j|\}}
\sum_{j=1}^m
\Biggl(\int\limits_{0}^{t}
\sigma_j^2(s)ds\Biggr)^{1/2}\le
$$

\vspace{2mm}
\begin{equation}
\label{comb2}
\le
C_{12}(\Delta_N)^{\min\{1-|\alpha|,\ 1/2-\max\limits_{1\le j\le m}|\beta_j|\}}.
\end{equation}

\vspace{5mm}

Finally, applying the above reasoning, 
we arrive at the following estimates

$$
\left({\sf M}\left\{\bigl|Q_t^{(k)}\bigr|^2\right\}\right)^{1/2}\le
$$

\vspace{-4mm}
$$
\le
\sum_{l=0}^{n_t-1}\Biggl({\sf M}\Biggl\{\Biggl(
\int\limits_{\tau_l}^{\tau_{l+1}}
\sigma_0(s)\left({\bf a}^{(k)}({\bf y}_l,\tau_l)-
{\bf a}^{(k)}({\bf x}_{\tau_l},\tau_l)\right)ds\Biggr)^2\Biggr\}
\Biggr)^{1/2}+
$$

\vspace{-4mm}
$$
+
\Biggl({\sf M}\Biggl\{\Biggl(~\int\limits_{\tau_{n_t}}^{t}
\sigma_0(s)\left({\bf a}^{(k)}({\bf y}_{\tau_{n_t}},\tau_{n_t})-
{\bf a}^{(k)}({\bf x}_{\tau_{n_t}},\tau_{n_t})\right)ds\Biggr)^2\Biggr\}
\Biggr)^{1/2}\le
$$

\vspace{-4mm}
$$
\le
\sum_{l=0}^{n_t-1}\Biggl(\int\limits_{\tau_l}^{\tau_{l+1}}
\left|\sigma_0(s)\right|ds
\int\limits_{\tau_l}^{\tau_{l+1}}
\left|\sigma_0(s)\right|{\sf M}\left\{\left({\bf a}^{(k)}({\bf y}_l,\tau_l)-
{\bf a}^{(k)}({\bf x}_{\tau_l},\tau_l)\right)^2\right\}
ds
\Biggr)^{1/2}+
$$
$$
+
\Biggl(\int\limits_{\tau_{n_t}}^{t} \left|\sigma_0(s)\right|ds
\int\limits_{\tau_{n_t}}^{t} \left|\sigma_0(s)\right|{\sf M}\left\{
\left({\bf a}^{(k)}({\bf y}_{\tau_{n_t}},\tau_{n_t})-
{\bf a}^{(k)}({\bf x}_{\tau_{n_t}},\tau_{n_t})\right)^2\right\}
ds\Biggr)^{1/2}\le
$$
$$
\le
C_{13}\Biggl(\sum_{l=0}^{n_t-1}\Biggl(\int\limits_{\tau_l}^{\tau_{l+1}}
\left|\sigma_0(s)\right|ds
\int\limits_{\tau_l}^{\tau_{l+1}}
\left|\sigma_0(s)\right|{\sf M}\bigl\{\left\Vert{\bf y}_l-{\bf x}_{\tau_l}\right\Vert^2\bigr\}
ds
\Biggr)^{1/2}+\Biggr.
$$
$$
\Biggl.+
\Biggl(\int\limits_{\tau_{n_t}}^{t} \left|\sigma_0(s)\right|ds
\int\limits_{\tau_{n_t}}^{t} \left|\sigma_0(s)\right|
{\sf M}\bigl\{\left\Vert{\bf y}_{n_t}-
{\bf x}_{\tau_{n_t}}\right\Vert^2\bigr\}
ds\Biggr)^{1/2}\Biggr)=
$$
$$
=
C_{13}\Biggl(\sum_{l=0}^{n_t-1} \int\limits_{\tau_l}^{\tau_{l+1}}
\left|\sigma_0(s)\right|ds
\sqrt{{\sf M}\bigl\{\left\Vert{\bf y}_l-{\bf x}_{\tau_l}\right\Vert^2\bigr\}}
+\Biggr.
$$
$$
\Biggl.
+
\int\limits_{\tau_{n_t}}^{t} \left|\sigma_0(s)\right|ds
\sqrt{{\sf M}\bigl\{\left\Vert{\bf y}_{n_t}-
{\bf x}_{\tau_{n_t}}\right\Vert^2\bigr\}}\Biggr)\le
$$
$$
\le
C_{13}\int\limits_{0}^{t}
\left|\sigma_0(s)\right|
\sqrt{Z_s}ds\le 
$$
$$
\le
C_{13}\Biggl(\int\limits_{0}^{t}
\left|\sigma_0(s)\right|ds\Biggr)^{1/2}
\Biggl(\int\limits_{0}^{t}
\left|\sigma_0(s)\right|Z_sds\Biggr)^{1/2}\le 
$$
$$
\le
C_{14}
\Biggl(\int\limits_{0}^{t}
\left|\sigma_0(s)\right|Z_sds\Biggr)^{1/2},
$$

\vspace{3mm}
\noindent
and

\vspace{-3mm}
$$
\left({\sf M}\left\{\bigl|P_t^{(k)}\bigr|^2\right\}\right)^{1/2}\le
$$

\vspace{-3mm}
$$
\le
\sum_{j=1}^m
\Biggl({\sf M}\Biggl\{\Biggl(\sum_{l=0}^{n_t-1}
\int\limits_{\tau_l}^{\tau_{l+1}}
\sigma_j(s)\left(B^{(kj)}({\bf y}_l,\tau_l)-B^{(kj)}({\bf x}_{\tau_l},\tau_l)\right)
d{\bf w}_s^{(j)}+\Biggr.\Biggr.\Biggr.
$$
$$
\Biggl.\Biggl.\Biggl.+\int\limits_{\tau_{n_t}}^{t}
\sigma_j(s)\left(B^{(kj)}({\bf y}_{n_t},\tau_{n_t})-B^{(kj)}({\bf x}_{\tau_{n_t}},\tau_{n_t})\right)
d{\bf w}_s^{(j)}\Biggr)^2\Biggr\}
\Biggr)^{1/2}=
$$
$$
=\sum_{j=1}^m
\Biggl(\sum_{l=0}^{n_t-1}
\int\limits_{\tau_l}^{\tau_{l+1}}
\sigma_j^2(s){\sf M}\left\{\left(B^{(kj)}({\bf y}_l,\tau_l)-
B^{(kj)}({\bf x}_{\tau_l},\tau_l)\right)^2\right\}
ds+\Biggr.
$$
$$
\Biggl.+
\int\limits_{\tau_{n_t}}^{t}\sigma_j^2(s){\sf M}\left\{
\left(B^{(kj)}({\bf y}_{n_t},\tau_{n_t})-B^{(kj)}({\bf x}_{\tau_{n_t}},\tau_{n_t})\right)^2\right\}
ds\Biggr)^{1/2} \le
$$
$$
\le C_{15}\sum_{j=1}^m
\Biggl(\sum_{l=0}^{n_t-1}
\int\limits_{\tau_l}^{\tau_{l+1}}
\sigma_j^2(s){\sf M}\bigl\{\left\Vert{\bf y}_l-{\bf x}_{\tau_l}\right\Vert^2\bigr\}
ds+\Biggr.
$$
$$
\Biggl.+
\int\limits_{\tau_{n_t}}^{t}\sigma_j^2(s){\sf M}\bigl\{\left\Vert{\bf y}_{n_t}-
{\bf x}_{\tau_{n_t}}\right\Vert^2\bigr\}
ds\Biggr)^{1/2} \le
$$
$$
\le C_{15}\sum_{j=1}^m
\Biggl(
\int\limits_{0}^{t}
\sigma_j^2(s) Z_s ds\Biggr)^{1/2}. 
$$

\vspace{3mm}

Then,
\begin{equation}
\label{comb3}
{\sf M}\left\{\bigl|Q_t^{(k)}\bigr|^2\right\}\le 
\left(C_{14}\right)^2
\int\limits_{0}^{t}
\left|\sigma_0(s)\right|Z_sds,
\end{equation}
$$
{\sf M}\left\{\bigl|P_t^{(k)}\bigr|^2\right\}\le
\left(C_{15}\right)^2 m\int\limits_{0}^{t} \left(\sum_{j=1}^m
\sigma_j^2(s)\right) Z_s ds = 
$$
\begin{equation}
\label{comb4}
=C_{16}
\int\limits_{0}^{t} \left(\sum_{j=1}^m
\sigma_j^2(s)\right) Z_s ds.
\end{equation}

\vspace{2mm}

Combining (\ref{888.14}), (\ref{888.15}), (\ref{useful11}), 
and (\ref{comb1}), (\ref{comb2}), (\ref{comb3}), (\ref{comb4}), 
and taking into account condition (\ref{888.06}), 
we arrive at the inequality
$$
Z_u\le K_3 (\Delta_N)^{2\min\{1-|\alpha|,\ 1/2-\max\limits_{1\le j\le m}|\beta_j|\}}
+K_4\int\limits_0^u \left(\left|\sigma_0(s)\right|+
\sum_{j=1}^m
\sigma_j^2(s)\right) Z_{\theta}d\theta,
$$

\noindent
from which, using Gronwall's inequality \cite{webb}, we find
\begin{equation}
\label{888.240}
Z_u\le K_5(\Delta_N)^{2\min\{1-|\alpha|,\ 1/2-\max\limits_{1\le j\le m}|\beta_j|\}},\ \ \ u\in[0,T],
\end{equation}

\noindent
where the constants $K_3, K_4,$ and $K_5$ do not depend on $\Delta_N.$ 
Further, taking into account (\ref{888.240}), 
we have
$$
{\sf M}\{\left\Vert{\bf x}_T-{\bf y}_N\right\Vert\}\le K
(\Delta_N)^{\min\{1-|\alpha|,\ 1/2-\max\limits_{1\le j\le m}|\beta_j|\}},
$$

\noindent
where the constant $K$ is independent of $\Delta_N.$
Theorem~2 is proved.\ $\Box$

If we assume more restrictive conditions on the functions 
${\bf a}({\bf x},t),$ $B_j({\bf x}_{s},s)$ $(j=1,\ldots,m),$ the conditions 
on the singularities of the It\^{o} SDE (\ref{1600t}) can be 
relaxed (see Sect.~4). 
More precisely, let us
consider the following non-commutative linear It\^{o} SDE
\begin{equation}
\label{1.5.2svvxx1}
x_t=x_0+\int\limits_0^t \sigma_0(s)\left(a(s)x_{s}+c(s)\right)ds+
\sum\limits_{j=1}^m
\int\limits_0^t \sigma_j(s) \left(b_j(s)x_{s}+f_j(s)\right)d{\bf w}_{s}^{(j)},
\end{equation}

\noindent
where $x_0=x(0,\omega),$ $t\in [0, T],$ 
$\sigma_0(t)\in L_1([0, T],{\bf R}),$ $\sigma_j(t)\in L_2([0, T],{\bf R}),$
$a(s),$ $c(s),$ $b_j(s),$ 
$f_j(s):$ $[0, T]\to{\bf R}$ $(j=1,\ldots,m)$ are non-random functions, another 
notations are the same as in (\ref{1.5.2}).

By analogy with the proof of Theorem~2, we obtain the following theorem.

{\bf Theorem~3.}\ {\it Suppose that
$\sigma_0(t)=t^{-|\alpha|},$
$\sigma_j(t)=t^{-|\beta_j|},$ $\alpha\in (-1,0)$, $\beta_j\in (-1/2,0)$ $(j=1,\ldots,m)$,
the nonrandom functions $a(s),$ $c(s),$ $b_j(s),$ $f_j(s):$ $[0, T]\to{\bf R}$ $(j=1,\ldots,m)$
are continuous, and their first-order derivatives are uniformly bounded.
Furthermore, assume that
$$
\left({\sf M}\bigl\{\left|x_0-y_0\right|^2\bigr\}\right)^{1/2}\le
M\sqrt{\Delta_N}.
$$

\noindent
Then, the following estimate
\begin{equation}
\label{res2}
{\sf M}\{\left|x_T-y_N\right|\}\le K(\Delta_N)^{\min\{1-|\alpha|,\ 1/2-\max\limits_{1\le j\le m}|\beta_j|\}}
\end{equation}

\noindent 
holds$,$ where the constant $K$ is independent of $\Delta_N,$ 
and $x_T$ is the solution of the It\^{o} SDE {\rm (\ref{1.5.2svvxx1})} at $t=T.$
}

{\bf Remark~2.}\ {\it 
Note that the estimates {\rm (\ref{res0}), (\ref{res2})} are quite natural.
First$,$ when $\beta_1^2+\ldots+\beta_m^2>0$ or $1/2<|\alpha|<1,$ the strong convergence order of
method {\rm (\ref{euler1})} is less than $1/2,$ whereas 
for non-singular It\^{o} SDEs under standard conditions, the 
strong convergence order of the Euler method is $1/2$.
If $\alpha^2+\beta_1^2+\ldots+\beta_m^2=0,$ 
the strong convergence order of method {\rm (\ref{euler1})} is $1/2,$
which agrees with the theory of numerical methods
for It\^{o} SDEs without singularities.
}

\section{Convergence Rate of the New Modification (\ref{11}) of the Milstein Method 
for the It\^{o} SDE (\ref{1600t})}

Consider the It\^{o} SDE (\ref{1600t}) and recall 
the new modification (\ref{11}) of the Milstein method.
We have 
$$
{\bf y}_{p+1}={\bf y}_{p}+
{\bf a}({\bf y}_{p},\tau_p)\int\limits_{\tau_p}^{\tau_{p+1}} \sigma_0(s)ds+
\sum\limits_{j=1}^m
B_j({\bf y}_{p},\tau_p)\int\limits_{\tau_p}^{\tau_{p+1}} \sigma_j(s)d{\bf w}_{s}^{(j)}
+
$$
\begin{equation}
\label{Mil11}
+\sum\limits_{l,j=1}^m 
G^{(l)}B_j({\bf y}_p,\tau_p)
\int\limits_{\tau_p}^{\tau_{p+1}} \sigma_j(s)\int\limits_{\tau_p}^s
\sigma_l(\theta)d{\bf w}_{\theta}^{(l)}d{\bf w}_{s}^{(j)},
\end{equation}

\noindent
where the partition $\{\tau_p\}_{p=0}^N$ 
of the interval $[0, T]$ satisfies the condition (\ref{1984}),
$G^{(l)}$ $(j=1,\ldots,m)$ is defined by (\ref{oper10}),
${\bf y}_{\tau_p}\stackrel{\sf def}{=}{\bf y}_p.$

Let us estimate the convergence rate of the method (\ref{Mil11}).

{\bf Theorem 4}.\ {\it Let $\sigma_0(t)=t^{-|\alpha|},$ $\sigma_j(t)=t^{-|\beta_j|},$ 
$\alpha\in (-1/(2+\varepsilon), 0),$
$\beta_j\in (-1/(4+2\varepsilon), 0),$
where $\varepsilon$ is an arbitrary
small positive real number. Furthermore, 
let the following conditions be satisfied 
for all ${\bf x}, {\bf y}\in{\bf R}^n,$ $s, t\in[0,T]:$

\vspace{-6mm}
\begin{equation}
\label{Mil12}
\left({\sf M}\bigl\{\left\Vert{\bf x}_0-{\bf y}_0\right\Vert^2\bigr\}\right)^{1/2}\le
M \Delta_N,
\end{equation}

\vspace{-8mm}
$$
\left\Vert{\bf a}({\bf x},t)-{\bf a}({\bf y},t)\right\Vert+
\left\Vert B_j({\bf x},t)-B_j({\bf y},t)\right\Vert
+
\left\Vert G^{(l)}B_j({\bf x},t)-G^{(l)}B_j({\bf y},t)\right\Vert\le
$$
$$
\le C\left\Vert{\bf x}-{\bf y}\right\Vert,
$$

\vspace{-6mm}
\begin{equation}
\label{nomer1}
\left\Vert{\bf a}({\bf x},t)\right\Vert+
\left\Vert B_j({\bf x},t)\right\Vert+
\left\Vert G^{(l)}B_j({\bf x},t)\right\Vert
\le C(1+\left\Vert{\bf x}\right\Vert),
\end{equation}

\vspace{-6mm}
\begin{equation}
\label{nomer2}
\left\Vert \bar L_0{\bf a}({\bf x},t)\right\Vert +
\left\Vert L_0{\bf a}({\bf x},t)\right\Vert+
\left\Vert L_l{\bf a}({\bf x},t)\right\Vert
\le C(1+\left\Vert{\bf x}\right\Vert),
\end{equation}

\vspace{-6mm}
\begin{equation}
\label{nomer3}
\left\Vert \bar L_0 B_j({\bf x},t)\right\Vert +
\left\Vert L_0 B_j({\bf x},t)\right\Vert+
\left\Vert L_l B_j({\bf x},t)\right\Vert
\le C(1+\left\Vert{\bf x}\right\Vert),
\end{equation}

\vspace{-6mm}
\begin{equation}
\label{nomer4}
\left\Vert G^{(l)}{\bf a}({\bf x},t)\right\Vert
\le C(1+\left\Vert{\bf x}\right\Vert),
\end{equation}

\vspace{-6mm}
$$
\left\Vert G^{(l)}B_j({\bf x},s)-G^{(l)}B_j({\bf x},t)\right\Vert\le 
C\sqrt{s-t}(1+\left\Vert{\bf x}\right\Vert),
$$

\newpage
\noindent
where $L_0,$ $L_l,$ $G^{(l)},$ $\bar L_0$ are
defined by {\rm (\ref{oper10}), (\ref{oper11}), (\ref{Mil14}),}
$l, j=1,\ldots,m,$ $C, M$ are constants that are independent of 
$\Delta_N$ 
{\rm (}we assume that all partial derivatives involved in the 
inequalities listed above exist{\rm )}.
Then
\begin{equation}
\label{res0qqqq}
{\sf M}\{\left\Vert{\bf x}_T-{\bf y}_N\right\Vert\}\le
K(\Delta_N)^{\min\{1-|\alpha|,\
1-2\max\limits_{1\le j\le m}|\beta_j|\}}.
\end{equation}

\noindent
where the constant $K$ is independent of $\Delta_N,$ 
and ${\bf x}_T$ is the solution of the It\^{o} SDE {\rm (\ref{1600t})} at $t=T.$}

{\bf Proof.}\ Under the assumptions of this theorem, 
according to Theorem 1.7 \cite{Jan} and the results 
of Sect.~10, we have the estimates (\ref{useful1}), (\ref{1986}).

Consider (\ref{6s}) for $\tau=\tau_p,$ $t=\tau_{p+1}$
$$
{\bf x}_{\tau_{p+1}}={\bf x}_{\tau_p}+{\bf a}({\bf x}_{\tau_p},\tau_p)
\int\limits_{\tau_p}^{\tau_{p+1}} \sigma_0(s)ds+
\sum\limits_{j=1}^m
B_j({\bf x}_{\tau_p},\tau_p)\int\limits_{\tau_p}^{\tau_{p+1}} \sigma_j(s)d{\bf w}_{s}^{(j)}+
$$
\begin{equation}
\label{8121}
+\sum\limits_{l,j=1}^m 
\int\limits_{\tau_p}^{\tau_{p+1}} \sigma_j(s)\int\limits_{\tau_p}^s
\sigma_l(\theta)G^{(l)}B_j({\bf x}_{\theta},\theta)d{\bf w}_{\theta}^{(l)}d{\bf w}_{s}^{(j)}+
R_{{\tau_{p+1}},\tau_p}
\end{equation}

\noindent
w.~p.~1, where
\begin{equation}
\label{arc10}
R_{t,\tau}=R_{(1)t,\tau}+R_{(2)t,\tau}+R_{(3)t,\tau}+R_{(4)t,\tau}+R_{(5)t,\tau}+R_{(6)t,\tau}+
R_{(7)t,\tau},
\end{equation}

\noindent
and
$$
R_{(1)t,\tau}=
\int\limits_{\tau}^t \sigma_0(s)\int\limits_{\tau}^s
\bar L_0{\bf a}({\bf x}_{\theta},\theta)d\theta ds,
$$
\begin{equation}
\label{setka900}
R_{(2)t,\tau}=
\int\limits_{\tau}^t \sigma_0(s)\int\limits_{\tau}^s
\sigma_0(\theta)L_0{\bf a}({\bf x}_{\theta},\theta)d\theta ds,
\end{equation}
$$
R_{(3)t,\tau}=
\sum^{m}_{j=1}\int\limits_{\tau}^t \sigma_0(s)\int\limits_{\tau}^s
\sigma_j^2(\theta)L_j{\bf a}({\bf x}_{\theta},\theta)d\theta ds,
$$
$$
R_{(4)t,\tau}=\sum\limits_{l=1}^m
\int\limits_{\tau}^t \sigma_0(s)\int\limits_{\tau}^s
\sigma_l(\theta)G^{(l)}{\bf a}({\bf x}_{\theta},\theta)d{\bf w}_{\theta}^{(l)}ds,
$$
$$
R_{(5)t,\tau}=\sum\limits_{j=1}^m
\int\limits_{\tau}^t \sigma_j(s)\int\limits_{\tau}^s
\bar L_0 B_j({\bf x}_{\theta},\theta)d\theta d{\bf w}_{s}^{(j)},
$$
$$
R_{(6)t,\tau}=
\sum\limits_{j=1}^m
\int\limits_{\tau}^t \sigma_j(s)\int\limits_{\tau}^s
\sigma_0(\theta)L_0 B_j({\bf x}_{\theta},\theta)d\theta d{\bf w}_{s}^{(j)},
$$
$$
R_{(7)t,\tau}=
\sum\limits_{l,j=1}^m
\int\limits_{\tau}^t \sigma_j(s)\int\limits_{\tau}^s
\sigma_l^2(\theta)L_l B_j({\bf x}_{\theta},\theta)d\theta d{\bf w}_{s}^{(j)}.
$$

\vspace{3mm}

Define the stochastic process ${\bf y}_s,$ $s\in[\tau_p,\tau_{p+1}]$
$(p=0, 1,\ldots,N-1)$
by equality (\ref{Mil11}), in which it is necessary
to replace $\tau_{p+1}$ with $s.$ Consider also the representation
of the stochastic process ${\bf x}_s,$ $s\in[\tau_p,\tau_{p+1}]$
$(p=0, 1,\ldots,N-1)$
by equality (\ref{8121}), in which we
replace $\tau_{p+1}$ with $s.$

For all $u\in[0,T]$, let us set
$$
Z_u\stackrel{\sf def}{=}
\sup\limits_{0\le s\le u}{\sf M}\bigl\{
\left\Vert{\bf y}_s-{\bf x}_s\right\Vert^2\bigr\}.
$$

\vspace{2mm}

It is easy to see that
\begin{equation}
\label{333x}
{\bf y}_t-{\bf x}_t={\bf y}_0-{\bf x}_0+Q_t+P_t+R_t+S_t+V_t+D_t,
\end{equation}

\noindent
where 
$$
Q_t=\sum_{l=0}^{n_t-1}\int\limits_{\tau_l}^{\tau_{l+1}}
\sigma_0(s)\left({\bf a}({\bf y}_l,\tau_l)-{\bf a}({\bf x}_{\tau_l},\tau_l)\right)
ds+
$$
$$
+
\int\limits_{\tau_{n_t}}^{t}
\sigma_0(s)\left({\bf a}({\bf y}_{n_t},\tau_{n_t})-{\bf a}({\bf x}_{\tau_{n_t}},\tau_{n_t})\right)
ds,
$$

\vspace{2mm}
$$
P_t=\sum\limits_{j=1}^m\Biggl(\sum_{l=0}^{n_t-1}\int\limits_{\tau_l}^{\tau_{l+1}}
\sigma_j(s)\left(B_j({\bf y}_l,\tau_l)-B_j({\bf x}_{\tau_l},\tau_l)\right)
d{\bf w}_s^{(j)}+\Biggr.
$$
$$
\Biggl.+
\int\limits_{\tau_{n_t}}^{t}
\sigma_j(s)\left(B_j({\bf y}_{n_t},\tau_{n_t})-B_j({\bf x}_{\tau_{n_t}},\tau_{n_t})\right)
d{\bf w}_s^{(j)}\Biggr),
$$
$$
R_t=\sum\limits_{i,j=1}^m\Biggl(\sum_{l=0}^{n_t-1}\int\limits_{\tau_l}^{\tau_{l+1}}
\sigma_j(s)\int\limits_{\tau_l}^{s}\sigma_i(\theta)
\left(G^{(i)}B_j({\bf y}_l,\tau_l)-G^{(i)}B_j({\bf x}_{\tau_l},\tau_l)\right)
d{\bf w}_{\theta}^{(i)}d{\bf w}_s^{(j)}+\Biggr.
$$
$$
\Biggl.+
\int\limits_{\tau_{n_t}}^{t}
\sigma_j(s)
\int\limits_{\tau_{n_t}}^{s}
\sigma_i(\theta)
\left(G^{(i)}B_j({\bf y}_{n_t},\tau_{n_t})-G^{(i)}B_j({\bf x}_{\tau_{n_t}},\tau_{n_t})\right)
d{\bf w}_{\theta}^{(i)}d{\bf w}_s^{(j)}\Biggr),
$$

\vspace{-2mm}
$$
S_t=\sum\limits_{i,j=1}^m\Biggl(\sum_{l=0}^{n_t-1}\int\limits_{\tau_l}^{\tau_{l+1}}
\sigma_j(s)\int\limits_{\tau_l}^{s}\sigma_i(\theta)
\left(G^{(i)}B_j({\bf x}_{\tau_l},\tau_l)-G^{(i)}B_j({\bf x}_{\theta},\tau_l)\right)
d{\bf w}_{\theta}^{(i)}d{\bf w}_s^{(j)}+\Biggr.
$$
\begin{equation}
\label{setka400}
\Biggl.+
\int\limits_{\tau_{n_t}}^{t}
\sigma_j(s)
\int\limits_{\tau_{n_t}}^{s}
\sigma_i(\theta)
\left(G^{(i)}B_j({\bf x}_{\tau_{n_t}},\tau_{n_t})-G^{(i)}B_j({\bf x}_{\theta},\tau_{n_t})\right)
d{\bf w}_{\theta}^{(i)}d{\bf w}_s^{(j)}\Biggr),
\end{equation}

\vspace{-2mm}
$$
V_t=\sum\limits_{i,j=1}^m\Biggl(\sum_{l=0}^{n_t-1}\int\limits_{\tau_l}^{\tau_{l+1}}
\sigma_j(s)\int\limits_{\tau_l}^{s}\sigma_i(\theta)
\left(G^{(i)}B_j({\bf x}_{\theta},\tau_l)-G^{(i)}B_j({\bf x}_{\theta},\theta)\right)
d{\bf w}_{\theta}^{(i)}d{\bf w}_s^{(j)}+\Biggr.
$$
\begin{equation}
\label{setka401}
\Biggl.+
\int\limits_{\tau_{n_t}}^{t}
\sigma_j(s)
\int\limits_{\tau_{n_t}}^{s}
\sigma_i(\theta)
\left(G^{(i)}B_j({\bf x}_{\theta},\tau_{n_t})-G^{(i)}B_j({\bf x}_{\theta},\theta)\right)
d{\bf w}_{\theta}^{(i)}d{\bf w}_s^{(j)}\Biggr),
\end{equation}

\begin{equation}
\label{one5}
D_t=-\left(\sum_{l=0}^{n_t-1}R_{{\tau_{l+1}},\tau_l}+
R_{t,\tau_{n_t}}\right),
\end{equation}

\vspace{1mm}
\noindent
where
$$
n_t=
\max\limits_{j}\{j:\ \tau_j\le t\}.
$$

\vspace{1mm}

From (\ref{333x}), we obtain for $u\in [0, T]$
$$
Z_u\stackrel{\sf def}{=}\sup\limits_{0\le t\le u}
{\sf M}\bigl\{\left\Vert{\bf y}_t-{\bf x}_t\right\Vert ^2\bigr\}\le 
$$
$$
\le
\sup\limits_{0\le t\le u}
K_1\left({\sf M}\bigl\{\left\Vert{\bf y}_0-{\bf x}_0\right\Vert^2\bigr\}+
{\sf M}\bigl\{\left\Vert Q_t\right\Vert^2\bigr\}+
{\sf M}\bigl\{\left\Vert P_t\right\Vert^2\bigr\} +\right.
$$
\begin{equation}
\label{999a}
\left.+
{\sf M}\bigl\{\left\Vert R_t\right\Vert^2\bigr\}+{\sf M}\bigl\{\left\Vert S_t\right\Vert^2\bigr\}+
{\sf M}\bigl\{\left\Vert V_t\right\Vert^2\bigr\}+
{\sf M}\bigl\{\left\Vert D_t\right\Vert^2\bigr\}\right),
\end{equation}

\vspace{2mm}
\noindent
where $K_1$ is a constant.

Let us estimate one by one
the terms on the right-hand side of (\ref{999a}).
Using (\ref{comb3}), (\ref{comb4}), we obtain 
\begin{equation}
\label{comb30}
{\sf M}\bigl\{\left\Vert Q_t\right\Vert^2\bigr\}\le 
K_{2}
\int\limits_{0}^{t}
\left|\sigma_0(s)\right|Z_sds,
\end{equation}
\begin{equation}
\label{comb40}
{\sf M}\bigl\{\left\Vert P_t\right\Vert^2\bigr\}\le
K_{3}
\int\limits_{0}^{t} \left(\sum_{j=1}^m
\sigma_j^2(s)\right) Z_s ds,
\end{equation}

\noindent
where here and below, $K_i$ $(i=2,3,\ldots)$ 
will denote constants that are independent of $\Delta_N$.

Applying Minkowski's inequality for mathematical expectations,
the conditions of Theorem~4, It\^{o}'s isometry, and
(\ref{useful2}), we obtain
$$
\left({\sf M}\left\{\bigl|R_t^{(k)}\bigr|^2\right\}\right)^{1/2}\le
$$
$$
\le \sum\limits_{i,j=1}^m \Biggl({\sf M}\Biggl\{\Biggl(\sum_{l=0}^{n_t-1}\int\limits_{\tau_l}^{\tau_{l+1}}
\sigma_j(s)\int\limits_{\tau_l}^{s}\sigma_i(\theta)
\left(G^{(i)}B^{(kj)}({\bf y}_l,\tau_l)-G^{(i)}B^{(kj)}({\bf x}_{\tau_l},\tau_l)\right)
\times\Biggr.\Biggr.
$$
$$
\times
d{\bf w}_{\theta}^{(i)}d{\bf w}_s^{(j)}+
$$
$$
\Biggl.\Biggl.\Biggl.+
\int\limits_{\tau_{n_t}}^{t}
\sigma_j(s)
\int\limits_{\tau_{n_t}}^{s}
\sigma_i(\theta)
\left(G^{(i)}B^{(kj)}({\bf y}_{n_t},\tau_{n_t})-G^{(i)}B^{(kj)}({\bf x}_{\tau_{n_t}},\tau_{n_t})\right)
\times
\Biggl.\Biggl.\Biggl.
$$
$$
\times
d{\bf w}_{\theta}^{(i)}d{\bf w}_s^{(j)}\Biggr)^2\Biggr\}\Biggr)^{1/2}=
$$
$$
=\sum\limits_{i,j=1}^m \Biggl(\sum_{l=0}^{n_t-1}\int\limits_{\tau_l}^{\tau_{l+1}}
\sigma_j^2(s)\int\limits_{\tau_l}^{s}\sigma_i^2(\theta)
{\sf M}\left\{\left(G^{(i)}B^{(kj)}({\bf y}_l,\tau_l)-G^{(i)}B^{(kj)}({\bf x}_{\tau_l},\tau_l)\right)^2
\right\}
d\theta ds +
$$
$$
+
\int\limits_{\tau_{n_t}}^{t}
\sigma_j^2(s)
\int\limits_{\tau_{n_t}}^{s}
\sigma_i^2(\theta)
{\sf M}\left\{\left(G^{(i)}B^{(kj)}({\bf y}_{n_t},\tau_{n_t})
-G^{(i)}B^{(kj)}({\bf x}_{\tau_{n_t}},\tau_{n_t})\right)^2\right\}
d\theta ds \Biggr)^{1/2}\le
$$
$$
\le K_4\sum\limits_{i,j=1}^m \Biggl(\sum_{l=0}^{n_t-1}\int\limits_{\tau_l}^{\tau_{l+1}}
\sigma_j^2(s)\int\limits_{\tau_l}^{s}\sigma_i^2(\theta)
{\sf M}\left\{\left\Vert {\bf y}_l-{\bf x}_{\tau_l}\right\Vert^2\right\} d\theta ds +
\Biggr.
$$
$$
\Biggl.+
\int\limits_{\tau_{n_t}}^{t}
\sigma_j^2(s)
\int\limits_{\tau_{n_t}}^{s}
\sigma_i^2(\theta)
{\sf M}\left\{\left\Vert{\bf y}_{n_t}-
{\bf x}_{\tau_{n_t}}\right\Vert^2\right\}
d\theta ds \Biggr)^{1/2}\le
$$
$$
\le K_4\sum\limits_{i,j=1}^m \Biggl(\sum_{l=0}^{n_t-1}\int\limits_{\tau_l}^{\tau_{l+1}}
\sigma_j^2(s)ds \int\limits_{\tau_l}^{\tau_{l+1}}\sigma_i^2(\theta)
{\sf M}\left\{\left\Vert {\bf y}_l-{\bf x}_{\tau_l}\right\Vert^2\right\} d\theta  +
\Biggr.
$$
$$
\Biggl.+
\int\limits_{\tau_{n_t}}^{t}
\sigma_j^2(s)ds
\int\limits_{\tau_{n_t}}^{t}
\sigma_i^2(\theta)
{\sf M}\left\{\left\Vert{\bf y}_{n_t}-
{\bf x}_{\tau_{n_t}}\right\Vert^2\right\}
d\theta \Biggr)^{1/2}\le
$$
$$
\le K_5\sum\limits_{i=1}^m \Biggl(\sum_{l=0}^{n_t-1}
\int\limits_{\tau_l}^{\tau_{l+1}}\sigma_i^2(\theta)
{\sf M}\left\{\left\Vert {\bf y}_l-{\bf x}_{\tau_l}\right\Vert^2\right\} d\theta  +
\Biggr.
$$
$$
\Biggl.
+\int\limits_{\tau_{n_t}}^{t}
\sigma_i^2(\theta)
{\sf M}\left\{\left\Vert{\bf y}_{n_t}-
{\bf x}_{\tau_{n_t}}\right\Vert^2\right\}
d\theta \Biggr)^{1/2}\le
$$
$$
\le K_5\sum\limits_{i=1}^m \Biggl(
\int\limits_{0}^{t}\sigma_i^2(\theta) Z_{\theta}
d\theta \Biggr)^{1/2}.
$$

Then,
$$
{\sf M}\left\{\bigl|R_t^{(k)}\bigr|^2\right\}\le
\left(K_{5}\right)^2 m\int\limits_{0}^{t} \left(\sum_{i=1}^m
\sigma_i^2(s)\right) Z_s ds = 
$$
\begin{equation}
\label{comb45}
=\bar K_6
\int\limits_{0}^{t} \left(\sum_{i=1}^m
\sigma_i^2(s)\right) Z_s ds.
\end{equation}

Applying Minkowski's inequality for mathematical expectations,
the conditions of Theorem~4, It\^{o}'s isometry,
(\ref{useful2}), (\ref{1981}), and (\ref{1986}), we obtain
$$
\left({\sf M}\left\{\bigl|S_t^{(k)}\bigr|^2\right\}\right)^{1/2}\le
$$
$$
\le \sum\limits_{i,j=1}^m \Biggl({\sf M}\Biggl\{\Biggl(\sum_{l=0}^{n_t-1}\int\limits_{\tau_l}^{\tau_{l+1}}
\sigma_j(s)\int\limits_{\tau_l}^{s}\sigma_i(\theta)
\left(G^{(i)}B^{(kj)}({\bf x}_{\tau_l},\tau_l)-G^{(i)}B^{(kj)}({\bf x}_{\theta},\tau_l)\right)
\times\Biggr.\Biggr.
$$
$$
\times
d{\bf w}_{\theta}^{(i)}d{\bf w}_s^{(j)}+
$$
$$
\Biggl.\Biggl.\Biggl.+
\int\limits_{\tau_{n_t}}^{t}
\sigma_j(s)
\int\limits_{\tau_{n_t}}^{s}
\sigma_i(\theta)
\left(G^{(i)}B^{(kj)}({\bf x}_{\tau_{n_t}},\tau_{n_t})-G^{(i)}B^{(kj)}({\bf x}_{\theta},\tau_{n_t})\right)
\times
\Biggl.\Biggl.\Biggl.
$$
$$
\times
d{\bf w}_{\theta}^{(i)}d{\bf w}_s^{(j)}\Biggr)^2\Biggr\}\Biggr)^{1/2}=
$$
$$
=\hspace{-0.5mm}\sum\limits_{i,j=1}^m \Biggl(\sum_{l=0}^{n_t-1}\int\limits_{\tau_l}^{\tau_{l+1}}
\hspace{-0.5mm}\sigma_j^2(s)\int\limits_{\tau_l}^{s}\hspace{-0.5mm}\sigma_i^2(\theta)
{\sf M}\left\{\left(G^{(i)}B^{(kj)}({\bf x}_{\tau_l},\tau_l)-G^{(i)}B^{(kj)}({\bf x}_{\theta},\tau_l)\right)^2
\right\}
d\theta ds +
$$
$$
+
\int\limits_{\tau_{n_t}}^{t}
\sigma_j^2(s)
\int\limits_{\tau_{n_t}}^{s}
\sigma_i^2(\theta)
{\sf M}\left\{\left(G^{(i)}B^{(kj)}({\bf x}_{\tau_{n_t}},\tau_{n_t})
-G^{(i)}B^{(kj)}({\bf x}_{\theta},\tau_{n_t})\right)^2\right\}
d\theta ds \Biggr)^{1/2}\le
$$
$$
\le K_6\sum\limits_{i,j=1}^m \Biggl(\sum_{l=0}^{n_t-1}\int\limits_{\tau_l}^{\tau_{l+1}}
\sigma_j^2(s)\int\limits_{\tau_l}^{s}\sigma_i^2(\theta)
{\sf M}\left\{\left\Vert {\bf x}_{\tau_l}-{\bf x}_{\theta}\right\Vert^2\right\} d\theta ds +
\Biggr.
$$
\begin{equation}
\label{setka405}
\Biggl.+
\int\limits_{\tau_{n_t}}^{t}
\sigma_j^2(s)
\int\limits_{\tau_{n_t}}^{s}
\sigma_i^2(\theta)
{\sf M}\left\{\left\Vert {\bf x}_{\tau_{n_t}}-
{\bf x}_{\theta}\right\Vert^2\right\}
d\theta ds \Biggr)^{1/2}\le
\end{equation}
$$
\le K_7 (\Delta_N)^{\min\{1-|\alpha|,\ 1/2-\max\limits_{1\le j\le m}|\beta_j|\}}\times
$$
$$
\times
\sum\limits_{i,j=1}^m \Biggl(\sum_{l=0}^{n_t-1}\int\limits_{\tau_l}^{\tau_{l+1}}
\sigma_j^2(s)ds \int\limits_{\tau_l}^{\tau_{l+1}}\sigma_i^2(\theta)d\theta
+\int\limits_{\tau_{n_t}}^{t}
\sigma_j^2(s)ds
\int\limits_{\tau_{n_t}}^{t}
\sigma_i^2(\theta)d\theta\Biggr)^{1/2}\le
$$
$$
\le K_7 (\Delta_N)^{\min\{1-|\alpha|,\ 1/2-\max\limits_{1\le j\le m}|\beta_j|\}}\times
$$
$$
\times
\sum\limits_{i,j=1}^m \Biggl(\Biggl( \max\limits_{0\le l\le n_t-1}\Biggl\{
\int\limits_{\tau_l}^{\tau_{l+1}}
\sigma_j^2(s)ds\Biggr\} + \int\limits_{\tau_{n_t}}^{t}\sigma_j^2(s)ds\Biggr)
\int\limits_{0}^{t}\sigma_i^2(\theta)d\theta
\Biggr)^{1/2}\le
$$
$$
\le K_8 (\Delta_N)^{\min\{1-|\alpha|,\ 1/2-\max\limits_{1\le j\le m}|\beta_j|\}}
\sum\limits_{j=1}^m \biggl( 
(\Delta_N)^{1-2|\beta_j|}+(\Delta_N)^{1-2|\beta_j|}\biggr)^{1/2}\le
$$
$$
\le K_9 (\Delta_N)^{\min\{1-|\alpha|,\ 1/2-\max\limits_{1\le j\le m}|\beta_j|\}}
\sum\limits_{j=1}^m 
(\Delta_N)^{1/2-|\beta_j|}\le
$$                                                                              
$$
\le K_{10} (\Delta_N)^{\min\{1-|\alpha|,\ 1/2-\max\limits_{1\le j\le m}|\beta_j|\}}
(\Delta_N)^{\min\limits_{1\le j\le m}\{1/2-|\beta_j|\}}=
$$
$$
= K_{10} (\Delta_N)^{\min\{1-|\alpha|,\ 1/2-\max\limits_{1\le j\le m}|\beta_j|\}+
\min\limits_{1\le j\le m}\{1/2-|\beta_j|\}}=
$$
\begin{equation}
\label{one2}
= K_{10} (\Delta_N)^{\min\{1-|\alpha|,\ 1/2-\max\limits_{1\le j\le m}|\beta_j|\}+
1/2-\max\limits_{1\le j\le m}|\beta_j|}=
\end{equation}
\begin{equation}
\label{one3}
= K_{10} (\Delta_N)^{\min\{3/2-|\alpha|-\max\limits_{1\le j\le m}|\beta_j|,\
1-2\max\limits_{1\le j\le m}|\beta_j|\}},
\end{equation}

\vspace{1mm}
\noindent
where the transition from (\ref{one2}) to (\ref{one3})
is carried out on the basis of the following relation:
$\min\left\{x,y\right\}+C=\min\left\{C+x,C+y\right\},$ where $C$ is a constant,
$x,y,C\in {\bf R}.$

By analogy with the above reasoning and using the conditions of Theorem~4
and (\ref{useful1}), we have
$$
\left({\sf M}\left\{\bigl|V_t^{(k)}\bigr|^2\right\}\right)^{1/2}\le
$$
$$
\le\sum\limits_{i,j=1}^m \Biggl(\sum_{l=0}^{n_t-1}\int\limits_{\tau_l}^{\tau_{l+1}}
\sigma_j^2(s)\int\limits_{\tau_l}^{s}\sigma_i^2(\theta)
{\sf M}\left\{\left(G^{(i)}B^{(kj)}({\bf x}_{\theta},\tau_l)-G^{(i)}B^{(kj)}({\bf x}_{\theta},\theta)\right)^2
\right\}
d\theta ds +
$$
$$
+
\int\limits_{\tau_{n_t}}^{t}
\sigma_j^2(s)
\int\limits_{\tau_{n_t}}^{s}
\sigma_i^2(\theta)
{\sf M}\left\{\left(G^{(i)}B^{(kj)}({\bf x}_{\theta},\tau_{n_t})
-G^{(i)}B^{(kj)}({\bf x}_{\theta},\theta)\right)^2\right\}
d\theta ds \Biggr)^{1/2}\le
$$
$$
\le K_{11}(\Delta_N)^{1/2}
\sum\limits_{i,j=1}^m \Biggl(\sum_{l=0}^{n_t-1}\int\limits_{\tau_l}^{\tau_{l+1}}
\sigma_j^2(s)\int\limits_{\tau_l}^{s}\sigma_i^2(\theta)
{\sf M}\left\{\left(1+\left\Vert{\bf x}_{\theta}\right\Vert\right)^2
\right\}
d\theta ds +
$$
$$
+
\int\limits_{\tau_{n_t}}^{t}
\sigma_j^2(s)
\int\limits_{\tau_{n_t}}^{s}
\sigma_i^2(\theta)
{\sf M}\left\{\left(1+\left\Vert{\bf x}_{\theta}\right\Vert\right)^2\right\}
d\theta ds \Biggr)^{1/2}\le
$$
\begin{equation}
\label{setka500}
\le K_{12}(\Delta_N)^{1/2}
\sum\limits_{i,j=1}^m \Biggl(\sum_{l=0}^{n_t-1}\int\limits_{\tau_l}^{\tau_{l+1}}
\sigma_j^2(s)ds\int\limits_{\tau_l}^{\tau_{l+1}}\sigma_i^2(\theta)
d\theta  +
\int\limits_{\tau_{n_t}}^{t}
\sigma_j^2(s)ds
\int\limits_{\tau_{n_t}}^{t}
\sigma_i^2(\theta)
d\theta  \Biggr)^{1/2}\le
\end{equation}
\begin{equation}
\label{one9000}
\le K_{12}(\Delta_N)^{1/2}
(\Delta_N)^{1/2-\max\limits_{1\le j\le m}|\beta_j|}=
K_{12}(\Delta_N)^{1-\max\limits_{1\le j\le m}|\beta_j|}.
\end{equation}

\vspace{3mm}

Further, let us estimate $\left({\sf M}\left\{\bigl|D_t^{(k)}\bigr|^2\right\}\right)^{1/2}.$
We have (see (\ref{one5}))
\begin{equation}
\label{setka402}
D_t^{(k)}=-\sum_{i=1}^7 \left(\sum_{l=0}^{n_t-1}R^{(k)}_{(i)\tau_{l+1},\tau_l}+
R^{(k)}_{(i)t,\tau_{n_t}}\right),
\end{equation}

\noindent
where $R^{(k)}_{t,\tau}$ is the $k$th component of $R_{t,\tau}$, and
$R_{t,\tau}$ is defined by (\ref{arc10}).

Then,
\begin{equation}
\label{one1}
\left({\sf M}\left\{\bigl|D_t^{(k)}\bigr|^2\right\}\right)^{1/2}\le
\sum_{i=1}^7 \left({\sf M}\left\{\left(\sum_{l=0}^{n_t-1}R^{(k)}_{(i)\tau_{l+1},\tau_l}+
R^{(k)}_{(i)t,\tau_{n_t}}\right)^2\right\}\right)^{1/2}.
\end{equation}

Let us estimate one by one
the terms on the right-hand side of (\ref{one1}). In doing so, 
let us begin with the seventh term in (\ref{one1}).
Using Minkowski's inequality for mathematical expectations, 
It\^{o}'s isometry,
the Cauchy--Bunyakovsky inequality, the Fubini Theorem,
the conditions of Theorem~4, (\ref{useful1}), and (\ref{1981}), we have
$$
\left({\sf M}\left\{\left(\sum_{l=0}^{n_t-1}R^{(k)}_{(7)\tau_{l+1},\tau_l}+
R^{(k)}_{(7)t,\tau_{n_t}}\right)^2\right\}\right)^{1/2}=
$$
$$
=\Biggl({\sf M}\Biggl\{\Biggl(\sum\limits_{i,j=1}^m\Biggl(\sum_{l=0}^{n_t-1}
\int\limits_{\tau_l}^{\tau_{l+1}} \sigma_j(s)\int\limits_{\tau_l}^s
\sigma_i^2(\theta)L_i B^{(kj)}({\bf x}_{\theta},\theta)d\theta d{\bf w}_{s}^{(j)}+
\Biggr.\Biggr.\Biggr.\Biggr.
$$
$$
\Biggl.\Biggl.\Biggl.\Biggl.+
\int\limits_{\tau_{n_t}}^{t} \sigma_j(s)\int\limits_{\tau_{n_t}}^s
\sigma_i^2(\theta)L_i B^{(kj)}({\bf x}_{\theta},\theta)d\theta d{\bf w}_{s}^{(j)}
\Biggr)\Biggr)^2\Biggr\}\Biggr)^{1/2}\le
$$
$$
\le\sum\limits_{i,j=1}^m\Biggl({\sf M}\Biggl\{\Biggl(\sum_{l=0}^{n_t-1}
\int\limits_{\tau_l}^{\tau_{l+1}} \sigma_j(s)\int\limits_{\tau_l}^s
\sigma_i^2(\theta)L_i B^{(kj)}({\bf x}_{\theta},\theta)d\theta d{\bf w}_{s}^{(j)}+
\Biggr.\Biggr.\Biggr.
$$
$$
\Biggl.\Biggl.\Biggl.+
\int\limits_{\tau_{n_t}}^{t} \sigma_j(s)\int\limits_{\tau_{n_t}}^s
\sigma_i^2(\theta)L_i B^{(kj)}({\bf x}_{\theta},\theta)d\theta d{\bf w}_{s}^{(j)}
\Biggr)^2\Biggr\}\Biggr)^{1/2}=
$$
$$
=\sum\limits_{i,j=1}^m\Biggl(\sum_{l=0}^{n_t-1}
\int\limits_{\tau_l}^{\tau_{l+1}} \sigma_j^2(s){\sf M}\Biggl\{\Biggl(\int\limits_{\tau_l}^s
\sigma_i^2(\theta)L_i B^{(kj)}({\bf x}_{\theta},\theta)d\theta\Biggr)^2\Biggr\} ds+
\Biggr.
$$
$$
\Biggl.+
\int\limits_{\tau_{n_t}}^{t} \sigma_j^2(s){\sf M}\Biggl\{\Biggl(\int\limits_{\tau_{n_t}}^s
\sigma_i^2(\theta)L_i B^{(kj)}({\bf x}_{\theta},\theta)d\theta\Biggr)^2\Biggr\} ds
\Biggr)^{1/2}\le
$$
$$
\le\sum\limits_{i,j=1}^m\Biggl(\sum_{l=0}^{n_t-1}
\int\limits_{\tau_l}^{\tau_{l+1}} \sigma_j^2(s)\int\limits_{\tau_l}^s
\sigma_i^2(\theta)d\theta \int\limits_{\tau_l}^s \sigma_i^2(\theta)
{\sf M}\left\{\left(L_i B^{(kj)}({\bf x}_{\theta},\theta)\right)^2\right\}d\theta ds+
\Biggr.
$$                                            
$$
\Biggl.+
\int\limits_{\tau_{n_t}}^{t} \sigma_j^2(s)\int\limits_{\tau_{n_t}}^s
\sigma_i^2(\theta) d\theta 
\int\limits_{\tau_{n_t}}^s
\sigma_i^2(\theta) {\sf M}\left\{\left(L_i B^{(kj)}({\bf x}_{\theta},\theta)\right)^2\right\}
d\theta ds
\Biggr)^{1/2}\le
$$                                        
\begin{equation}
\label{setka600}
\le K_{13}\sum\limits_{i,j=1}^m\Biggl(\sum_{l=0}^{n_t-1}
\int\limits_{\tau_l}^{\tau_{l+1}} \sigma_j^2(s)ds 
\Biggl(\int\limits_{\tau_l}^{\tau_{l+1}}
\sigma_i^2(\theta)d\theta\Biggr)^2 
+
\int\limits_{\tau_{n_t}}^{t} \sigma_j^2(s)ds \Biggl(\int\limits_{\tau_{n_t}}^t
\sigma_i^2(\theta) d\theta\Biggr)^2
\Biggr)^{1/2}\le
\end{equation}
$$
\le K_{13}\sum\limits_{i,j=1}^m\Biggl(\Biggl(
\max\limits_{0\le l\le n_t-1}\Biggl\{
\Biggl(\int\limits_{\tau_l}^{\tau_{l+1}}
\sigma_i^2(\theta)d\theta\Biggr)^2\Biggr\} 
+
\Biggl(\int\limits_{\tau_{n_t}}^t
\sigma_i^2(\theta) d\theta\Biggr)^2\Biggr)\int\limits_{0}^{t} \sigma_j^2(s)ds
\Biggr)^{1/2}\le
$$                                        
$$
\le K_{14}\sum\limits_{i=1}^m\Biggl(
\max\limits_{0\le l\le n_t-1}\Biggl\{
\Biggl(\int\limits_{\tau_l}^{\tau_{l+1}}
\sigma_i^2(\theta)d\theta\Biggr)^2\Biggr\} 
+
\Biggl(\int\limits_{\tau_{n_t}}^t
\sigma_i^2(\theta) d\theta\Biggr)^2
\Biggr)^{1/2}\le
$$                                        
$$
\le K_{15}\sum\limits_{i=1}^m\Biggl(
(\Delta_N)^{2-4|\beta_i|}+(\Delta_N)^{2-4|\beta_i|}
\Biggr)^{1/2}\le
$$                                        
\begin{equation}
\label{one3007}
\le K_{16}\sum\limits_{i=1}^m
(\Delta_N)^{1-2|\beta_i|}\le K_{17}
(\Delta_N)^{\min\limits_{1\le i\le m}\{1-2|\beta_i|\}}=
K_{17}
(\Delta_N)^{1-2\max\limits_{1\le i\le m}|\beta_i|}.
\end{equation}

\vspace{2mm}

Similarly, using (\ref{1980}), we get
$$
\left({\sf M}\left\{\left(\sum_{l=0}^{n_t-1}R^{(k)}_{(6)\tau_{l+1},\tau_l}+
R^{(k)}_{(6)t,\tau_{n_t}}\right)^2\right\}\right)^{1/2}=
$$
$$
=\Biggl({\sf M}\Biggl\{\Biggl(\sum\limits_{j=1}^m\Biggl(\sum_{l=0}^{n_t-1}
\int\limits_{\tau_l}^{\tau_{l+1}} \sigma_j(s)\int\limits_{\tau_l}^s
\sigma_0(\theta)L_0 B^{(kj)}({\bf x}_{\theta},\theta)d\theta d{\bf w}_{s}^{(j)}+
\Biggr.\Biggr.\Biggr.\Biggr.
$$
\begin{equation}
\label{setka800}
\Biggl.\Biggl.\Biggl.\Biggl.+
\int\limits_{\tau_{n_t}}^{t} \sigma_j(s)\int\limits_{\tau_{n_t}}^s
\sigma_0(\theta)L_0 B^{(kj)}({\bf x}_{\theta},\theta)d\theta d{\bf w}_{s}^{(j)}
\Biggr)\Biggr)^2\Biggr\}\Biggr)^{1/2}\le
\end{equation}
$$
\le K_{18}\sum\limits_{j=1}^m\Biggl(\Biggl(
\max\limits_{0\le l\le n_t-1}\Biggl\{
\Biggl(\int\limits_{\tau_l}^{\tau_{l+1}}
\left|\sigma_0(\theta)\right|d\theta\Biggr)^{\hspace{-1mm}2}\Biggr\} 
+
\Biggl(\int\limits_{\tau_{n_t}}^t
\left|\sigma_0(\theta)\right| d\theta\Biggr)^{\hspace{-1mm}2}\Biggr)\int\limits_{0}^{t} \sigma_j^2(s)ds
\Biggr)^{1/2}\hspace{-2.5mm}\le
$$                                        
\begin{equation}
\label{one3006}
\le K_{19}
(\Delta_N)^{1-|\alpha|},
\end{equation}
$$
\left({\sf M}\left\{\left(\sum_{l=0}^{n_t-1}R^{(k)}_{(5)\tau_{l+1},\tau_l}+
R^{(k)}_{(5)t,\tau_{n_t}}\right)^2\right\}\right)^{1/2}=
$$
$$
=\Biggl({\sf M}\Biggl\{\Biggl(\sum\limits_{j=1}^m\Biggl(\sum_{l=0}^{n_t-1}
\int\limits_{\tau_l}^{\tau_{l+1}} \sigma_j(s)\int\limits_{\tau_l}^s
\bar L_0 B^{(kj)}({\bf x}_{\theta},\theta)d\theta d{\bf w}_{s}^{(j)}+
\Biggr.\Biggr.\Biggr.\Biggr.
$$
$$
\Biggl.\Biggl.\Biggl.\Biggl.+
\int\limits_{\tau_{n_t}}^{t} \sigma_j(s)\int\limits_{\tau_{n_t}}^s
\bar L_0 B^{(kj)}({\bf x}_{\theta},\theta)d\theta d{\bf w}_{s}^{(j)}
\Biggr)\Biggr)^2\Biggr\}\Biggr)^{1/2}\le
$$
$$
\le K_{20}\sum\limits_{j=1}^m\Biggl(\Biggl(
\max\limits_{0\le l\le n_t-1}\Biggl\{
\Biggl(\int\limits_{\tau_l}^{\tau_{l+1}}
d\theta\Biggr)^2\Biggr\} 
+
\Biggl(\int\limits_{\tau_{n_t}}^t
d\theta\Biggr)^2\Biggr)\int\limits_{0}^{t} \sigma_j^2(s)ds
\Biggr)^{1/2}\le
$$                                        

\begin{equation}
\label{one3005}
\le K_{21}
\Delta_N,
\end{equation}                                        

\vspace{2mm}

\begin{equation}
\label{one100}
\left({\sf M}\left\{\left(\sum_{l=0}^{n_t-1}R^{(k)}_{(4)\tau_{l+1},\tau_l}+
R^{(k)}_{(4)t,\tau_{n_t}}\right)^2\right\}\right)^{1/2}=
\end{equation}
$$
=\Biggl({\sf M}\Biggl\{\Biggl(\sum\limits_{j=1}^m\Biggl(\sum_{l=0}^{n_t-1}
\int\limits_{\tau_l}^{\tau_{l+1}} \sigma_j(\theta)G^{(j)}{\bf a}^{(k)}({\bf x}_{\theta},\theta)
\int\limits_{\theta}^{\tau_{l+1}}
\sigma_0(s)
ds d{\bf w}_{\theta}^{(j)}+
\Biggr.\Biggr.\Biggr.\Biggr.
$$
\begin{equation}
\label{one101}
\Biggl.\Biggl.\Biggl.\Biggl.+
\int\limits_{\tau_{n_t}}^{t} \sigma_j(\theta)
G^{(j)}{\bf a}^{(k)}({\bf x}_{\theta},\theta)
\int\limits_{\theta}^t
\sigma_0(s)ds d{\bf w}_{\theta}^{(j)}
\Biggr)\Biggr)^2\Biggr\}\Biggr)^{1/2}\le
\end{equation}
$$
\le K_{22}\sum\limits_{j=1}^m\Biggl(\Biggl(
\max\limits_{0\le l\le n_t-1}\Biggl\{
\Biggl(\int\limits_{\tau_l}^{\tau_{l+1}}
\left|\sigma_0(s)\right|ds\Biggr)^{\hspace{-1mm}2}\Biggr\} 
+
\Biggl(\int\limits_{\tau_{n_t}}^t
\left|\sigma_0(s)\right|ds\Biggr)^{\hspace{-1mm}2}\Biggr)\int\limits_{0}^{t} \sigma_j^2(\theta)d\theta
\Biggr)^{1/2}\hspace{-2.5mm}\le
$$                                        

\begin{equation}
\label{one3004}
\le K_{23}
(\Delta_N)^{1-|\alpha|},
\end{equation}

\vspace{5mm}
\noindent
where 
the transition from (\ref{one100}) to (\ref{one101})
is carried out on the basis of the It\^{o} formula
(also see the integration order replacement technique for 
iterated It\^{o} stochastic integrals \cite{15} (Chapter~3)).

Further, 
using Minkowski's inequality for mathematical expectations, 
the Cauchy--Bunyakovsky inequality, the Fubini Theorem,
the conditions of Theorem~4, (\ref{useful1}), and (\ref{1980}), we have
$$
\left({\sf M}\left\{\left(\sum_{l=0}^{n_t-1}R^{(k)}_{(3)\tau_{l+1},\tau_l}+
R^{(k)}_{(3)t,\tau_{n_t}}\right)^2\right\}\right)^{1/2}=
$$
$$
=\Biggl({\sf M}\Biggl\{\Biggl(\sum\limits_{j=1}^m\Biggl(\sum_{l=0}^{n_t-1}
\int\limits_{\tau_l}^{\tau_{l+1}} \sigma_0(s)\int\limits_{\tau_l}^s
\sigma_j^2(\theta)L_j {\bf a}^{(k)}({\bf x}_{\theta},\theta)d\theta ds+
\Biggr.\Biggr.\Biggr.\Biggr.
$$
$$
\Biggl.\Biggl.\Biggl.\Biggl.+
\int\limits_{\tau_{n_t}}^{t} \sigma_0(s)\int\limits_{\tau_{n_t}}^s
\sigma_j^2(\theta)L_j {\bf a}^{(k)}({\bf x}_{\theta},\theta)d\theta ds
\Biggr)\Biggr)^2\Biggr\}\Biggr)^{1/2}\le
$$
$$
\le \sum\limits_{j=1}^m\Biggl({\sf M}\Biggl\{\Biggl(\sum_{l=0}^{n_t-1}
\int\limits_{\tau_l}^{\tau_{l+1}} \sigma_0(s)\int\limits_{\tau_l}^s
\sigma_j^2(\theta)L_j {\bf a}^{(k)}({\bf x}_{\theta},\theta)d\theta ds+
\Biggr.\Biggr.\Biggr.
$$
$$
\Biggl.\Biggl.\Biggl.+
\int\limits_{\tau_{n_t}}^{t} \sigma_0(s)\int\limits_{\tau_{n_t}}^s
\sigma_j^2(\theta)L_j {\bf a}^{(k)}({\bf x}_{\theta},\theta)d\theta ds
\Biggr)^2\Biggr\}\Biggr)^{1/2}\le
$$
$$
\le \sum\limits_{j=1}^m\Biggl(\sum_{l=0}^{n_t-1}\Biggl({\sf M}\Biggl\{\Biggl(
\int\limits_{\tau_l}^{\tau_{l+1}} \sigma_0(s)\int\limits_{\tau_l}^s
\sigma_j^2(\theta)L_j {\bf a}^{(k)}({\bf x}_{\theta},\theta)d\theta ds\Biggr)^2\Biggr\}\Biggr)^{1/2}+
\Biggr.
$$
$$
\Biggl.+
\Biggl({\sf M}\Biggl\{\Biggl(\int\limits_{\tau_{n_t}}^{t} \sigma_0(s)\int\limits_{\tau_{n_t}}^s
\sigma_j^2(\theta)L_j {\bf a}^{(k)}({\bf x}_{\theta},\theta)d\theta ds
\Biggr)^2\Biggr\}\Biggr)^{1/2}\Biggr)\le
$$
$$
\le \sum\limits_{j=1}^m\Biggl(\sum_{l=0}^{n_t-1}\Biggl(
\int\limits_{\tau_l}^{\tau_{l+1}} \left|\sigma_0(s)\right|ds 
\int\limits_{\tau_l}^{\tau_{l+1}} \left|\sigma_0(s)\right|
{\sf M}\Biggl\{\Biggl(\int\limits_{\tau_l}^s
\sigma_j^2(\theta)L_j {\bf a}^{(k)}({\bf x}_{\theta},\theta)d\theta\Biggr)^2\Biggr\} ds\Biggr)^{1/2}+
\Biggr.
$$
\begin{equation}
\label{one2000}
\Biggl.+
\Biggl(\int\limits_{\tau_{n_t}}^{t} \left|\sigma_0(s)\right|ds
\int\limits_{\tau_{n_t}}^{t} \left|\sigma_0(s)\right|
{\sf M}\Biggl\{\Biggl(\int\limits_{\tau_{n_t}}^s
\sigma_j^2(\theta)L_j {\bf a}^{(k)}({\bf x}_{\theta},\theta)d\theta \Biggr)^2\Biggr\}ds
\Biggr)^{1/2}\Biggr),
\end{equation}

\noindent
and
$$
{\sf M}\Biggl\{\Biggl(\int\limits_{\tau_l}^s
\sigma_j^2(\theta)L_j {\bf a}^{(k)}({\bf x}_{\theta},\theta)d\theta\Biggr)^2\Biggr\}\le
$$
$$
\le
\int\limits_{\tau_l}^s
\sigma_j^2(\theta) d\theta 
\int\limits_{\tau_l}^s
\sigma_j^2(\theta)
{\sf M}\left\{\left(
L_j {\bf a}^{(k)}({\bf x}_{\theta},\theta)\right)^2\right\}d\theta \le
$$
$$
\le K_{24}
\int\limits_{\tau_l}^s
\sigma_j^2(\theta) d\theta 
\int\limits_{\tau_l}^s
\sigma_j^2(\theta)
{\sf M}\left\{\left(
1+\left\Vert ({\bf x}_{\theta}\right\Vert \right)^2\right\}d\theta \le
$$
\begin{equation}
\label{one2001}
\le K_{25}
\Biggl(\int\limits_{\tau_l}^{\tau_{l+1}}
\sigma_j^2(\theta) d\theta \Biggr)^2,
\end{equation}

\begin{equation}
\label{one2002}
{\sf M}\Biggl\{\Biggl(\int\limits_{\tau_{n_t}}^s
\sigma_j^2(\theta)L_j {\bf a}^{(k)}({\bf x}_{\theta},\theta)d\theta\Biggr)^2\Biggr\}\le
K_{25}
\Biggl(\int\limits_{\tau_{n_t}}^t
\sigma_j^2(\theta) d\theta \Biggr)^2.
\end{equation}

Combining (\ref{one2000}), (\ref{one2001}), and (\ref{one2002}), we get
$$
\left({\sf M}\left\{\left(\sum_{l=0}^{n_t-1}R^{(k)}_{(3)\tau_{l+1},\tau_l}+
R^{(k)}_{(3)t,\tau_{n_t}}\right)^2\right\}\right)^{1/2}\le
$$
$$
\le K_{26}\sum\limits_{j=1}^m\Biggl(\sum_{l=0}^{n_t-1}\Biggl(
\Biggl(\int\limits_{\tau_l}^{\tau_{l+1}} \left|\sigma_0(s)\right|ds\Biggr)^2
\Biggl(\int\limits_{\tau_l}^{\tau_{l+1}}
\sigma_j^2(\theta) d\theta \Biggr)^2
\Biggr)^{1/2}+
\Biggr.
$$
$$
\Biggl.+
\Biggl(\Biggl(\int\limits_{\tau_{n_t}}^{t} \left|\sigma_0(s)\right|ds\Biggr)^2
\Biggl(\int\limits_{\tau_{n_t}}^t
\sigma_j^2(\theta) d\theta \Biggr)^2
\Biggr)^{1/2}\Biggr)=
$$
\begin{equation}
\label{setka700}
= K_{26}\sum\limits_{j=1}^m\Biggl(\sum_{l=0}^{n_t-1}
\int\limits_{\tau_l}^{\tau_{l+1}} \left|\sigma_0(s)\right|ds
\int\limits_{\tau_l}^{\tau_{l+1}}
\sigma_j^2(\theta) d\theta 
+
\int\limits_{\tau_{n_t}}^{t} \left|\sigma_0(s)\right|ds
\int\limits_{\tau_{n_t}}^t
\sigma_j^2(\theta) d\theta 
\Biggr)=
\end{equation}
$$
\le K_{26}\sum\limits_{j=1}^m\Biggl(
\max\limits_{0\le l\le n_t-1}\Biggl\{
\int\limits_{\tau_l}^{\tau_{l+1}}
\left|\sigma_0(s)\right|ds\Biggr\} 
+
\int\limits_{\tau_{n_t}}^t
\left|\sigma_0(s)\right|ds\Biggr)\int\limits_{0}^{t} \sigma_j^2(\theta)d\theta
\le
$$
\begin{equation}
\label{one3002}
\le K_{27}
(\Delta_N)^{1-|\alpha|}.
\end{equation}

Similarly, we have
$$
\left({\sf M}\left\{\left(\sum_{l=0}^{n_t-1}R^{(k)}_{(2)\tau_{l+1},\tau_l}+
R^{(k)}_{(2)t,\tau_{n_t}}\right)^2\right\}\right)^{1/2}\le
$$
\begin{equation}
\label{one3001}
\le K_{28}
(\Delta_N)^{1-|\alpha|}.
\end{equation}

$$
\left({\sf M}\left\{\left(\sum_{l=0}^{n_t-1}R^{(k)}_{(1)\tau_{l+1},\tau_l}+
R^{(k)}_{(1)t,\tau_{n_t}}\right)^2\right\}\right)^{1/2}\le
$$
$$
\le K_{29}\Biggl(
\max\limits_{0\le l\le n_t-1}\Biggl\{
\int\limits_{\tau_l}^{\tau_{l+1}}
ds\Biggr\} 
+
\int\limits_{\tau_{n_t}}^t
ds\Biggr)\int\limits_{0}^{t}\left|\sigma_0(s)\right| d\theta
\le
$$                                        
\begin{equation}
\label{one3000}
\le K_{30}
\Delta_N.
\end{equation}

\vspace{2mm}

Combining (\ref{one1}), (\ref{one3007}), (\ref{one3006}), (\ref{one3005}), (\ref{one3004}),
(\ref{one3002}), (\ref{one3001}), (\ref{one3000}), we obtain 
$$
\left({\sf M}\left\{\bigl|D_t^{(k)}\bigr|^2\right\}\right)^{1/2}\le
\sum_{i=1}^7 \left({\sf M}\left\{\left(\sum_{l=0}^{n_t-1}R^{(k)}_{(i)\tau_{p+1},\tau_p}+
R^{(k)}_{(i)t,\tau_{n_t}}\right)^2\right\}\right)^{1/2}\le
$$
$$
\le K_{31}\left(
(\Delta_N)^{1-2\max\limits_{1\le i\le m}|\beta_i|}
+(\Delta_N)^{1-|\alpha|}+\Delta_N
\right)\le K_{32} (\Delta_N)^{\min\{1-|\alpha|,\ 1-2\max\limits_{1\le i\le m}|\beta_i|\}}.
$$

Then,
\begin{equation}
\label{one4000}
{\sf M}\bigl\{\left\Vert|D_t\right\Vert^2\bigr\}\le K_{33}
(\Delta_N)^{2\min\{1-|\alpha|,\ 1-2\max\limits_{1\le i\le m}|\beta_i|\}}.
\end{equation}

\vspace{3mm}

Since $1-\left|\alpha\right|>0$ and $1/2-\left|\beta_j\right|>0$ $(j=1,\ldots,m),$ then
$$
2\min\biggl\{1-|\alpha|,\ 1-2\max\limits_{1\le i\le m}|\beta_i|\biggr\}\le
$$
$$
\le
2\min\biggl\{(1-|\alpha|) + \left(1/2- \max\limits_{1\le i\le m}|\beta_i|\right),\ 
1-2\max\limits_{1\le i\le m}|\beta_i|\biggr\}=
$$
$$
=
2\min\biggl\{3/2-|\alpha|- \max\limits_{1\le i\le m}|\beta_i|,\ 
1-2\max\limits_{1\le i\le m}|\beta_i|\biggr\}\le
$$
\begin{equation}
\label{one5000}
\le
2\left(1-2\max\limits_{1\le i\le m}|\beta_i|\right)\le 
2\left(1-\max\limits_{1\le i\le m}|\beta_i|\right)\le 2.
\end{equation}

\vspace{2mm}

Thus, using (\ref{999a}) and (\ref{Mil12}), (\ref{comb30}), (\ref{comb40}),
(\ref{comb45}), (\ref{one3}), (\ref{one9000}), (\ref{one4000}), (\ref{one5000}), we obtain
the inequality
$$
Z_u\le K_{34} (\Delta_N)^{2\min\{1-|\alpha|,\ 1-2\max\limits_{1\le i\le m}|\beta_i|\}}
+K_{35}\int\limits_0^u \left(\left|\sigma_0(s)\right|+
\sum_{j=1}^m
\sigma_j^2(s)\right) Z_{\theta}d\theta,
$$

\noindent
from which, using Gronwall's inequality \cite{webb}, we find
\begin{equation}
\label{one6000}
Z_u\le K_{36}(\Delta_N)^{2\min\{1-|\alpha|,\ 1-2\max\limits_{1\le i\le m}|\beta_i|\}},\ \ \ u\in[0,T],
\end{equation}
where the constants $K_{34}, K_{35},$ and $K_{36}$ do not depend on $\Delta_N.$ 
Further, taking into account (\ref{one6000}), 
we have
\begin{equation}
\label{one6001}
{\sf M}\{\left\Vert{\bf x}_T-{\bf y}_N\right\Vert\}\le K
(\Delta_N)^{\min\{1-|\alpha|,\ 1-2\max\limits_{1\le i\le m}|\beta_i|\}},
\end{equation}

\noindent
where the constant $K$ is independent of $\Delta_N.$
Theorem~4 is proved.\ $\Box$

\vspace{1mm}

Note that
$$
\min\{x,y\}=\frac{x+y-\left|x-y\right|}{2},\ \ \ \max\{x,y\}=\frac{x+y+\left|x-y\right|}{2}.
$$

Then 
$$
\min\{C-x,C-y\}=\frac{C-x+C-y-\left|C-x-(C-y)\right|}{2}=
$$
$$
=
\frac{2C-(x+y)-\left|x-y\right|}{2}=C-\frac{x+y+\left|x-y\right|}{2}=
$$
\begin{equation}
\label{one10000}
=C-\max\{x,y\},
\end{equation}

\noindent
where $C, x,y\in{\bf R}.$

From (\ref{one6001}) and (\ref{one10000}), we get
$$
{\sf M}\{\left\Vert{\bf x}_T-{\bf y}_N\right\Vert\}\le K
(\Delta_N)^{1-\max\{|\alpha|,\ 2\max\limits_{1\le i\le m}|\beta_i|\}}.
$$

Recall that for the modified Euler method (\ref{euler1}), we obtained the estimate
(see (\ref{res0}))
\begin{equation}
\label{one7777}
{\sf M}\{\left\Vert{\bf x}_T-{\bf y}_N\right\Vert\}\le K
(\Delta_N)^{\min\{1-|\alpha|,\ 1/2-\max\limits_{1\le j\le m}|\beta_j|\}},
\end{equation}

\noindent
where the constant $K$ is independent of $\Delta_N.$

Using (\ref{one10000}), we can write the estimate 
(\ref{one7777}) in the following form
$$
{\sf M}\{\left\Vert{\bf x}_T-{\bf y}_N\right\Vert\}\le K
(\Delta_N)^{1/2-\max\{|\alpha|-1/2,\ \max\limits_{1\le j\le m}|\beta_j|\}},
$$

\noindent
where the constant $K$ is independent of $\Delta_N.$

If we assume more restrictive conditions on the functions 
${\bf a}({\bf x},t),$ $B_j({\bf x}_{s},s)$ $(j=1,\ldots,m),$ the conditions 
on the singularities of the It\^{o} SDE (\ref{1600t}) can be 
relaxed (see Sect.~4). 
More precisely, let us
consider the following non-commutative linear It\^{o} SDE
\begin{equation}
\label{1.5.2svvxx1s}
x_t=x_0+\int\limits_0^t \sigma_0(s)\left(a(s)x_{s}+c(s)\right)ds+
\sum\limits_{j=1}^m
\int\limits_0^t \sigma_j(s) \left(b_j(s)x_{s}+f_j(s)\right)d{\bf w}_{s}^{(j)},
\end{equation}

\noindent
where $x_0=x(0,\omega),$
$t\in [0, T],$ 
$\sigma_0(t)\in L_1([0, T],{\bf R}),$ $\sigma_j(t)\in L_2([0, T],{\bf R}),$
$a(s),$ $c(s),$ $b_j(s),$ $f_j(s):$ $[0, T]\to{\bf R}$ $(j=1,\ldots,m)$
are non-random functions, another 
notations are the same as in (\ref{1.5.2}).

Then, by analogy with the proof of Theorem~4, we obtain the following theorem.

{\bf Theorem 5}. {\it Let $\sigma_0(t)=t^{-|\alpha|},$
$\sigma_j(t)=t^{-|\beta_j|},$ $\alpha\in (-1,0)$, $\beta_j\in (-1/2,0)$ $(j=1,\ldots,m)$.
Furthermore, 
the non-random functions $a(s),$ $c(s),$ $b_j(s),$ $f_j(s):$ $[0, T]\to{\bf R}$ $(j=1,\ldots,m)$
are continuous, and their first-order derivatives are uniformly bounded.
Moreover, let
$$
\left({\sf M}\bigl\{\left|x_0-y_0\right|^2\bigr\}\right)^{1/2}\le
M\Delta_N.
$$
Then, the following estimate 
\begin{equation}
\label{res2qqqq}
{\sf M}\{\left|x_T-y_N\right|\}\le 
K(\Delta_N)^{\min\{1-|\alpha|,\
1-2\max\limits_{1\le j\le m}|\beta_j|\}},
\end{equation}

\vspace{2mm}
\noindent 
holds$,$ where the constant $K$ is independent of $\Delta_N,$ 
and $x_T$ is the solution of the It\^{o} SDE {\rm (\ref{1.5.2svvxx1s})} at $t=T.$
}

{\bf Remark~3.}\ {\it 
Note that the estimates {\rm (\ref{res0qqqq}),  (\ref{res2qqqq})} are quite natural.
First$,$ when $\beta_1^2+\ldots+\beta_m^2>0$ or $0<|\alpha|<1,$ the strong convergence order of
method {\rm (\ref{Mil11})} is less than $1,$ whereas 
for non-singular It\^{o} SDEs under standard conditions, the 
strong convergence order of the Milstein method is $1$.
If $\alpha^2+\beta_1^2+\ldots+\beta_m^2=0,$ 
the strong convergence order of method {\rm (\ref{Mil11})} is $1,$
which agrees with the theory of numerical methods
for It\^{o} SDEs without singularities.
}

{\bf Remark~4.}\ {\it If we set $\sigma_j(s)B^{(kj)}(s)\equiv 0$ $(k=1,\ldots,n,\ j=1,\ldots,m)$ in 
{\rm (\ref{1600t}),} i.e.$,$ carry out the transition from the It\^{o} SDE 
{\rm (\ref{1600t})} to the system of ordinary differential equations$,$ then
the order of strong convergence of the modified Euler method
{\rm (\ref{euler1})} and the new modification {\rm (\ref{11})} of the Milstein method will be equal
to $1-\left|\alpha\right|.$ 
This corresponds to the fact that 
$$
1/2-\max\limits_{1\le j\le m}|\beta_j|\ \ \ \hbox{and}\ \ \ 
1-2\max\limits_{1\le i\le m}|\beta_i|
$$
must be removed from the right-hand sides of {\rm (\ref{res0})} and 
{\rm (\ref{res0qqqq}),} respectively.}

{\bf Remark~5.}\ {\it It is easy to see that Theorems {\rm 2--5} remain valid if we set
$$
\sigma_0(\tau)=\tau^{-\left|\alpha\right|}\left(
\prod\limits_{i=1}^M \left|\tau-\theta_i\right|^{-\left|\alpha\right|}\right)
(T-\tau)^{-\left|\alpha\right|},
$$
$$
\label{nov1001}
\sigma_j(\tau)=\tau^{-\left|\beta_j\right|}\left(
\prod\limits_{i=1}^M \left|\tau-\theta_i\right|^{-\left|\beta_j\right|}\right)
(T-\tau)^{-\left|\beta_j\right|},
$$
where $\alpha, \beta_j$ are the same as in Theorems {\rm 2--5,} $\tau\in [0, T],$ 
$\theta_1<\ldots<\theta_M,$ 
$\theta_1,\ldots,\theta_M\in (0, T),$ $M\in {\bf N},$ $j=1,\ldots,m.$
For the case under consideration$,$
we can prove {\rm(}see below{\rm )} the following estimates {\rm (see (\ref{1980}), (\ref{1981}))}
\begin{equation}
\label{start1}
\int\limits_s^t \left\vert\sigma_0(\tau)\right\vert d\tau
\le \frac{C(t-s)^{1-|\alpha|}}{1-|\alpha|},\ \ \ 
\int\limits_s^t \sigma_j^2(\tau) d\tau
\le \frac{K(t-s)^{1-2|\beta_j|}}{1-2|\beta_j|},
\end{equation}
where $j=1,\ldots,m,$ 
the constants $C$ and $K$ do not depend on $t, s, i, N$
and $[s, t] \subseteq [\tau_i, \tau_{i+1}]$ for some $i$ $(0\le i\le N-1),$ i.e.
$|t-s|\le \tau_{i+1}-\tau_i \le \Delta_N,$ $i=0,1,\ldots,N-1$
{\rm (see (\ref{1984}))}.
Therefore, the estimate {\rm (\ref{1983})} remains true.
}

Let us prove the first estimate in (\ref{start1}) (the second estimate 
in (\ref{start1}) is proved similarly).

{\bf Lemma~1.}\ {\it Let $0 < \left|\alpha\right| < 1$, let 
$0 = \theta_0 < \theta_1 < \dots < \theta_M < \theta_{M+1} = T$ 
be fixed points, and consider the function
$$
\sigma_0(\tau) = \tau^{-\left|\alpha\right|} 
\left( \prod_{i=1}^M |\tau - \theta_i|^{-\left|\alpha\right|} \right) 
(T - \tau)^{-\left|\alpha\right|}.
$$
If $\Delta_N = \max\limits_{0\le j\le N-1} \left|\tau_{j+1} - \tau_j\right|$ 
{\rm (see (\ref{1984}))} satisfies
$\Delta_N < \delta/2,$ where $\delta=\min\limits_{0 \le k \le M}\left|\theta_{k+1} - \theta_k\right|,$
then there exists a constant $C > 0$, independent of $s, t, j,$
and $N,$ such that for any subinterval $[\tau_j, \tau_{j+1}]$ of
$[0, T]$ and any $[s, t] \subseteq [\tau_j, \tau_{j+1}],$ the following 
estimates hold$:$
\begin{equation}
\label{start3}
\int\limits_s^t \left|\sigma_0(\tau)\right|d\tau 
\le C (t - s)^{1 - \left|\alpha\right|} 
\le C (\tau_{j+1} - \tau_j)^{1 - \left|\alpha\right|} 
\le C \left(\Delta_N\right)^{1 - \left|\alpha\right|},
\end{equation}
and$,$ in particular {\rm (}for $s = \tau_j$ and $t = \tau_{j+1}${\rm ),}
$$
\int\limits_{\tau_j}^{\tau_{j+1}} \left|\sigma_0(\tau)\right|d\tau 
\le C (\tau_{j+1} - \tau_j)^{1 - \left|\alpha\right|} 
\le C \left(\Delta_N\right)^{1 - \left|\alpha\right|}.
$$
}

{\bf Proof.} By the condition $\Delta_N < \delta/2$, any interval 
$[s, t] \subseteq [\tau_j, \tau_{j+1}]$ (having length 
$t - s \le \tau_{j+1} - \tau_j \le \Delta_N < \delta/2$) contains at most 
one singular point $\theta_k$ (either in its interior or at an 
endpoint). Indeed, suppose, on the contrary, that $\theta_k,\theta_{k+1}\in [s,t].$ 
Then, 
$$
\delta\le \left|\theta_{k+1}-\theta_k\right|\le t-s\le \Delta_N\le \delta/2.
$$

\noindent
Hence, the above statement holds.

{\it Case {\rm 1.} The interval $[s, t]$ contains no singular points.} Then,
let $\theta_k$ be the singular point closest to $[s, t]$. All other 
singularities $\theta_i$ ($i \neq k$) are separated from $[s, t]$ 
by a distance of at least $\delta / 4$. 
Indeed, we have
$$
\delta\le \left|\theta_{k+1}-\theta_k\right|\le \left|\theta_{k+1}-x_{k+1}\right|+
\left|x_{k+1}-x_k\right|+\left|x_k-\theta_k\right|=
$$
$$
={\rm dist}([s,t],\theta_{k+1})+\left|x_{k+1}-x_k\right|+{\rm dist}([s,t],\theta_k)\le
2{\rm dist}([s,t],\theta_{k+1})+\left|t-s\right|\le
$$
\begin{equation}
\label{start4}
\le 2{\rm dist}([s,t],\theta_{k+1})+\Delta_N\le 2{\rm dist}([s,t],\theta_{k+1})+\delta/2,
\end{equation}

\noindent
where 
$$
\left|x_{k+1}-\theta_{k+1}\right|={\rm dist}([s,t],\theta_{k+1}),\ \ \ 
\left|x_k-\theta_k\right|={\rm dist}([s,t],\theta_k),
$$ 

\noindent
and ${\rm dist}([s,t],\theta)=\min\limits_{x\in [s,t]}\left|x-\theta\right|,$ $x_k, x_{k+1}\in [s,t].$

Then, from (\ref{start4}), we obtain ${\rm dist}([s,t],\theta_{k+1})\ge \delta/4.$
Replacing $\theta_{k+1}$ with $\theta_{k-1},$ we get ${\rm dist}([s,t],\theta_{k-1})\ge \delta/4.$
Therefore, ${\rm dist}([s,t],\theta_{i})\ge \delta/4$ $(i\ne k).$
Furthermore, 
$\left|\sigma_0(\tau)\right| \le C_1 |\tau - \theta_k|^{-\left|\alpha\right|}$ on $[s, t]$, where the 
constant $C_1$ is independent of $s, t, j$, and $N$. Applying 
the inequality 
\begin{equation}
\label{start2}
\left|B^{1-\left|\alpha\right|} - A^{1-\left|\alpha\right|}\right|
\le |B - A|^{1-\left|\alpha\right|}\ \ \ 
(A, B \ge 0,\  0 < \left|\alpha\right| < 1),
\end{equation}
we obtain
$$
\int\limits_s^t \left|\sigma_0(\tau)\right|d\tau 
\le C_1 \int\limits_s^t |\tau - \theta_k|^{-\left|\alpha\right|}d\tau=
$$
$$
= C_1\ \frac{(t - \theta_k)^{1-\left|\alpha\right|} 
- (s - \theta_k)^{1-\left|\alpha\right|}}{1 - \left|\alpha\right|} 
\le C(t - s)^{1-\left|\alpha\right|}\ \ \ (\theta_k<s<t),
$$
$$
\int\limits_s^t \left|\sigma_0(\tau)\right|d\tau 
\le C_1 \int\limits_s^t |\tau - \theta_k|^{-\left|\alpha\right|}d\tau=
$$
$$
= C_1\ \frac{(\theta_k-s)^{1-\left|\alpha\right|} 
- (\theta_k-t)^{1-\left|\alpha\right|}}{1 - \left|\alpha\right|} 
\le C(t - s)^{1-\left|\alpha\right|}\ \ \ (s<t<\theta_k),
$$

\noindent
where $C=C_1/(1 - \left|\alpha\right|)$.

{\it Case {\rm 2}. The interval $[s, t]$ contains a single singular point 
$\theta_k$ {\rm (}in its interior or coinciding with one of the endpoints 
$s$ or $t${\rm )}.}
We have $\left|\sigma_0(\tau)\right| \le C_1 |\tau - \theta_k|^{-\left|\alpha\right|}$ on $[s, t]$, 
where the constant $C_1$ is independent of $s, t, j$, and $N$.
Further, we have
$$
\int\limits_s^t \left|\sigma_0(\tau)\right|d\tau=
\int\limits_s^{\theta_k} \left|\sigma_0(\tau)\right|d\tau+
\int\limits_{\theta_k}^t \left|\sigma_0(\tau)\right|d\tau\le
$$
$$
\le C_1\int\limits_s^{\theta_k} \left|\tau-\theta_k\right|^{-\left|\alpha\right|}d\tau+
C_1\int\limits_{\theta_k}^t \left|\tau - \theta_k\right|^{-\left|\alpha\right|}d\tau=
$$
$$
= C_1\int\limits_s^{\theta_k} (\theta_k-\tau)^{-\left|\alpha\right|}d\tau+
C_1\int\limits_{\theta_k}^t (\tau - \theta_k)^{-\left|\alpha\right|}d\tau=
$$
$$
=C_1\ \frac{(\theta_k-s)^{1-\left|\alpha\right|} 
+ (t - \theta_k)^{1-\left|\alpha\right|}}{1 - \left|\alpha\right|} \le C
(t-s)^{1-\left|\alpha\right|},
$$

\vspace{1mm}
\noindent
where $C = 2C_1/(1 - \left|\alpha\right|).$ 

For the cases $\theta_k = s$ and $\theta_k = t$, we obviously have
$$
\int\limits_s^t \left|\sigma_0(\tau)\right|d\tau \le C_1 \int\limits_s^t 
\left|\tau-\theta_k\right|^{-\left|\alpha\right|}d\tau =
$$
$$
=
C(t-s)^{1-\left|\alpha\right|},
$$

\noindent
where $C=C_1/(1 - \left|\alpha\right|).$

Applying the inequalities
$$
t - s \le \tau_{j+1} - \tau_j \le \Delta_N,
$$

\noindent
we arrive at all the 
desired estimates (\ref{start3}). $\Box$

\section{Refinement of Convergence Orders in Theorems 2--5}

We note that the convergence rate estimates 
in Theorems 2--5 are rather rough and can be refined in the 
case of a single singularity point $\tau=0$ and
\begin{equation}
\label{setka1}
\tau_p=p\Delta_N,\ \ \ p=0,1,\ldots,N,\ \ \ N\Delta_N=T,\ \ \ \Delta_N\in (0, 1),\ \ \ N>1.
\end{equation}

Namely, 
the following theorems hold.

{\bf Theorem 6}.\ {\it Suppose that the conditions of Theorem~{\rm 2} as well as 
the condition {\rm (\ref{setka1})}
are satisfied. Then$,$ for the modified Euler method {\rm (\ref{euler1})}
we have the following estimate
\begin{equation}
\label{setka50}
{\sf M}\{\left\Vert{\bf x}_T-{\bf y}_N\right\Vert\}\le K
(\Delta_N)^{\min\{1/2-\varepsilon,\ 2-2|\alpha|,\ 1-2\max\limits_{1\le j\le m}|\beta_j|\}},
\end{equation}

\noindent
where the constant $K$ is independent of $\Delta_N,$ 
${\bf x}_T$ is the solution of the It\^{o} SDE {\rm (\ref{1600t})} at $t=T,$
and $\varepsilon>0$ is an arbitrary
small positive real number.}

{\bf Theorem~7.}\ {\it 
Assume that the conditions of Theorem~{\rm 3}
as well as 
the condition {\rm (\ref{setka1})}
are satisfied. Then$,$ for the modified Euler method {\rm (\ref{euler1})}
the following estimate
$$
{\sf M}\{\left|x_T-y_N\right|\}\le 
K(\Delta_N)^{\min\{1/2-\varepsilon,\ 2-2|\alpha|,\ 1-2\max\limits_{1\le j\le m}|\beta_j|\}}
$$
holds$,$ where the constant $K$ is independent of $\Delta_N,$ 
$x_T$ is the solution of the It\^{o} SDE {\rm (\ref{1.5.2svvxx1})} at $t=T,$
and $\varepsilon>0$ is an arbitrary
small positive real number.}

{\bf Theorem 8}.\ {\it Suppose that the conditions of Theorem~{\rm 4}
as well as 
the condition {\rm (\ref{setka1})}
are satisfied. Then$,$ for the new modification {\rm (\ref{11})} of the Milstein  method 
we have the following estimate
$$
{\sf M}\{\left\Vert{\bf x}_T-{\bf y}_N\right\Vert\}\le 
$$

\vspace{-5mm}
$$
\le K
(\Delta_N)^{\min\{1-\varepsilon,\ 2-2|\alpha|,\ 3/2-3\max\limits_{1\le j\le m}|\beta_j|,\ 3/2-|\alpha|-
\max\limits_{1\le j\le m}|\beta_j|\}},
$$

\vspace{2mm}
\noindent
where the constant $K$ is independent of $\Delta_N,$ 
${\bf x}_T$ is the solution of the It\^{o} SDE {\rm (\ref{1600t})} at $t=T,$
and $\varepsilon>0$ is an arbitrary
small positive real number.}

{\bf Theorem~9.}\ {\it 
Assume that the conditions of Theorem~{\rm 5}
as well as 
the condition {\rm (\ref{setka1})}
are satisfied. Then$,$ for the new modification {\rm (\ref{11})} of the Milstein  method 
the following estimate
$$
{\sf M}\{\left|x_T-y_N\right|\}\le 
$$

\vspace{-7mm}
$$
\le
K(\Delta_N)^{\min\{1-\varepsilon,\ 2-2|\alpha|,\ 3/2-3\max\limits_{1\le j\le m}|\beta_j|,\ 3/2-|\alpha|-
\max\limits_{1\le j\le m}|\beta_j|\}}
$$

\vspace{2mm}
\noindent
holds$,$ where the constant $K$ is independent of $\Delta_N,$ 
$x_T$ is the solution of the It\^{o} SDE {\rm (\ref{1.5.2svvxx1s})} at $t=T,$
and $\varepsilon>0$ is an arbitrary
small positive real number.}

Let us prove Theorems 6 and 8 (the proofs of Theorems 7 and 9
are carried out similarly).

\subsection{Proof of Theorem 6}

The proof of Theorem~2 still holds for the most part
in proving Theorem~6; we only obtain some new estimates.

Throughout this proof, we assume 
that $\varepsilon > 0$ is an arbitrarily small positive constant, 
and $C$ denotes various positive constants independent of $\Delta_N$, $l$, $s\in [\tau_l, \tau_{l+1}]$.

{\bf Lemma~2.} {\it 
Let the condition {\rm (\ref{setka1})}
is satisfied and let
$\left|\alpha\right| \in (0, 1),$ $\left|\beta_j\right| \in (0, 1/2)$ $(j = 1, \dots, m)$. 
Then for any $l \ge 1$ and any $s \in [\tau_l, \tau_{l+1}],$ the following estimate 
\begin{equation}
\label{setka4}
{\sf M}\bigl\{\left\Vert {\bf x}_{\tau_l} - {\bf x}_s\right\Vert^2\bigr\} 
\le C \left(l^{-2\left|\alpha\right|} (\Delta_N)^{2 - 2\left|\alpha\right|}
+ l^{-2\max\limits_{1 \le j \le m}\left|\beta_j\right|} 
(\Delta_N)^{1 - 2\max\limits_{1 \le j \le m} \left|\beta_j\right|} \right)
\end{equation}
holds.
}

{\bf Proof.}\ Using (\ref{667xx}) for $r = 2$, we have
$$
{\sf M}\bigl\{\left\Vert{\bf x}_{\tau_l} - {\bf x}_s\right\Vert^2\bigr\} 
\le C \left( \left(\int\limits_{\tau_l}^s \tau^{-\left|\alpha\right|} d\tau\right)^2 + 
\sum_{j=1}^m \int\limits_{\tau_l}^s \tau^{-2\left|\beta_j\right|} d\tau \right).
$$

Applying (\ref{setka1}) and $0 < \tau_l \le \tau \le \tau_{l+1}$ for $l \ge 1,$ we obtain
\begin{equation}
\label{setka5}
\left(\int\limits_{\tau_l}^s \tau^{-\left|\alpha\right|} d\tau\right)^2 \le 
(\tau_l)^{-2\left|\alpha\right|} (\Delta_N)^2 = 
(l\Delta_N)^{-2\left|\alpha\right|} (\Delta_N)^2 = l^{-2\left|\alpha\right|} 
(\Delta_N)^{2 - 2\left|\alpha\right|},
\end{equation}
\begin{equation}
\label{setka2}
\int\limits_{\tau_l}^s \tau^{-2\left|\beta_j\right|} d\tau \le (\tau_l)^{-2\left|\beta_j\right|} \Delta_N.
\end{equation}

Since $0 < \tau_l/T \le 1$ and $-2\max\limits_{1 \le j \le m} \left|\beta_j\right| \le 
-2\left|\beta_j\right|,$ by the monotonicity of the exponential function 
with base at most $1,$ it follows that
$$
\left(\frac{\tau_l}{T}\right)^{-2|\beta_j|}  
\le \left(\frac{\tau_l}{T}\right)^{-2\max\limits_{1 \le j \le m} \left|\beta_j\right|}, 
$$
\begin{equation}
\label{setka3}
(\tau_l)^{-2\left|\beta_j\right|}  
\le T^{2\max\limits_{1 \le j \le m} \left|\beta_j\right| - 2\left|\beta_j\right|} 
(\tau_l)^{-2\max\limits_{1 \le j \le m} 
\left|\beta_j\right|} \le C (\tau_l)^{-2\max\limits_{1 \le j \le m} \left|\beta_j\right|}.
\end{equation}

Combining (\ref{setka2}), (\ref{setka3}), and the relation $\tau_l=l\Delta_N,$ we get
\begin{equation}
\label{setka6}
\sum_{j=1}^m \int\limits_{\tau_l}^s \tau^{-2\left|\beta_j\right|} d\tau 
\le C (\tau_l)^{-2\max\limits_{1 \le j \le m} \left|\beta_j\right|} \Delta_N
= C l^{-2\max\limits_{1 \le j \le m} \left|\beta_j\right|} 
(\Delta_N)^{1 - 2\max\limits_{1 \le j \le m} \left|\beta_j\right|}.
\end{equation}

\noindent
Using (\ref{setka5}) and (\ref{setka6}), we obtain (\ref{setka4}).
Lemma 2 is proved. $\Box$

{\bf Lemma~3.} {\it 
Let the condition {\rm (\ref{setka1})}
is satisfied. Then$,$
for any $p > 0$ and arbitrarily small $\varepsilon > 0,$ the following inequality holds$:$
\begin{equation}
\sum_{l=1}^{N-1} l^{-p} \le 
\left\{
\begin{matrix}
C (\Delta_N)^{p - 1} & \hbox{for}\ \ \ 0<p<1\cr
C (\Delta_N)^{-\varepsilon} & \hbox{for}\ \ \ p=1\cr
C & \hbox{for}\ \ \ p>1
\end{matrix}\right..
\end{equation}
}

{\bf Proof.}\ Since $f(x) = x^{-p}$ is decreasing for $x > 0,$ we have
$$ 
\sum\limits_{l=1}^{N-1} l^{-p} = 1 + \sum\limits_{l=2}^{N-1} l^{-p} = 
1 + \sum\limits_{l=2}^{N-1} \int\limits_{l-1}^{l} l^{-p}\,dx \le 
$$
\begin{equation}
\label{setka10}
\le 1 + \sum\limits_{l=2}^{N-1} \int\limits_{l-1}^{l} x^{-p} dx = 
1 + \int\limits_{1}^{N-1} x^{-p} dx \le 1 + \int\limits_{1}^{N} x^{-p} dx.
\end{equation}

For $0 < p < 1,$ we have
\begin{equation}
\label{setka11}
1 + \int\limits_{1}^{N} x^{-p} dx = \frac{N^{1-p} - p}{1 - p} < \frac{N^{1-p}}{1 - p} = 
\frac{T^{1-p}}{1 - p} (\Delta_N)^{p-1} = C(\Delta_N)^{p-1}. 
\end{equation}

Then, combining (\ref{setka10}) and (\ref{setka11}), we get for $0<p<1$
$$
\sum_{l=1}^{N-1} l^{-p} \le 
C (\Delta_N)^{p - 1}.
$$

For $p = 1,$ we obtain
\begin{equation}
\label{setka12}
1 + \int\limits_1^N x^{-1} dx = 1 + \ln N = 1 + \ln\left(\frac{T}{\Delta_N}\right).
\end{equation}

Using the elementary inequality $\ln u < u$ for $u = (T/\Delta_N)^{\varepsilon} > 0,$ we get
$$
\varepsilon \ln\left(\frac{T}{\Delta_N}\right) < 
\frac{T^{\varepsilon}}{(\Delta_N)^{\varepsilon}} 
\implies \ln\left(\frac{T}{\Delta_N}\right) < 
\frac{T^{\varepsilon}}{\varepsilon} (\Delta_N)^{-\varepsilon},
$$

\noindent
which together with $1 \le T^{\varepsilon} (\Delta_N)^{-\varepsilon}$ and (\ref{setka10}),
(\ref{setka12}) gives
$$
\sum_{l=1}^{N-1} l^{-1} \le 1 + \ln\left(\frac{T}{\Delta_N}\right) \le C (\Delta_N)^{-\varepsilon}.
$$

For $p > 1,$ applying (\ref{setka10}), we have
$$
\sum_{l=1}^{N-1} l^{-p} \le 1 + \int\limits_1^N x^{-p} dx = 1 + \frac{1 - N^{1 - p}}{p - 1} 
\le 1 + \frac{1}{p - 1} = \frac{p}{p - 1} = C.
$$

\noindent
Lemma 3 is proved. $\Box$

Let us return to the proof of Theorem~6. 
To this end, we need to derive new estimates for
$\left({\sf M}\left\{\bigl|R_t^{(k)}\bigr|^2\right\}\right)^{1/2}$
and $\left({\sf M}\left\{\bigl|G_t^{(k)}\bigr|^2\right\}\right)^{1/2}$
(see (\ref{start5}) and (\ref{start6})).

Further, the incomplete last subinterval $[\tau_{n_t}, t] \subseteq [\tau_{n_t}, \tau_{n_t+1}]$ 
($n_t \le N - 1$) is absorbed into the sums as the term with $l = n_t$.

Using (\ref{setka15}), we obtain
$$
\left({\sf M}\left\{\bigl|R_t^{(k)}\bigr|^2\right\}\right)^{1/2} \le 
C \sum_{l=0}^{n_t} \left( \int\limits_{\tau_l}^{\tau_{l+1}} \tau^{-|\alpha|} d\tau 
\int\limits_{\tau_l}^{\tau_{l+1}} s^{-|\alpha|} {\sf M}\bigl\{\left\Vert{\bf x}_{\tau_l} - 
{\bf x}_s\right\Vert^2\bigr\} ds \right)^{1/2}.
$$

For $l = 0$, using (\ref{1983}) on the first step 
and the estimate
$$
\int\limits_0^{\Delta_N} \tau^{-\left|\alpha\right|} d\tau \le C (\Delta_N)^{1 - \left|\alpha\right|},
$$ 
we obtain
$$
\Biggl( \int\limits_0^{\Delta_N} \tau^{-\left|\alpha\right|} d\tau
\int\limits_0^{\Delta_N} s^{-\left|\alpha\right|} 
{\sf M}\bigl\{\left\Vert{\bf x}_0 - {\bf x}_s\right\Vert^2\bigr\} ds \Biggr)^{1/2} \le
$$
$$
\le C (\Delta_N)^{1 - \left|\alpha\right| + 
\min\{1 - \left|\alpha\right|,\  1/2 - \max\limits_{1 \le j \le m} \left|\beta_j\right|\}}=
$$
\begin{equation}
\label{setka20}
= C (\Delta_N)^{\min\{2 - 2\left|\alpha\right|,\ 
3/2 - \left|\alpha\right| - \max\limits_{1 \le j \le m} \left|\beta_j\right|\}}.
\end{equation}

\vspace{2mm}

For $l \ge 1$, we have (see (\ref{setka5}))
$$
\int_{\tau_l}^{\tau_{l+1}} \tau^{-|\alpha|} d\tau \le l^{-|\alpha|} (\Delta_N)^{1 - |\alpha|}.
$$ 

Applying the inequality $\sqrt{u+v}\le \sqrt{u}+\sqrt{v}$ $(u,v\ge 0)$ as well as
Lemma~2 and Lemma~3, we deduce
$$
\sum_{l=1}^{n_t} \left( \int\limits_{\tau_l}^{\tau_{l+1}} \tau^{-|\alpha|} d\tau 
\int\limits_{\tau_l}^{\tau_{l+1}} s^{-|\alpha|} {\sf M}\bigl\{\left\Vert{\bf x}_{\tau_l} - 
{\bf x}_s\right\Vert^2\bigr\} ds \right)^{1/2}\le
$$
$$
\le C (\Delta_N)^{2 - 2\left|\alpha\right|} \sum_{l=1}^{N-1} l^{-2\left|\alpha\right|} +
C (\Delta_N)^{3/2 - \left|\alpha\right| - \max\limits_{1 \le j \le m} \left|\beta_j\right|} 
\sum_{l=1}^{N-1} l^{-\left|\alpha\right| - \max\limits_{1 \le j \le m} \left|\beta_j\right|} \le
$$
\begin{equation}
\label{setka21}
\le C (\Delta_N)^{\min\{1 - \varepsilon,\ 2 - 2\left|\alpha\right|\}} + 
C (\Delta_N)^{\min\{1/2 - \varepsilon,\ 3/2 - \left|\alpha\right| - 
\max\limits_{1 \le j \le m} \left|\beta_j\right|\}}.
\end{equation}

Combining (\ref{setka20}) and (\ref{setka21}), we obtain ($\Delta_N < 1$)
$$
\left({\sf M}\left\{\bigl|R_t^{(k)}\bigr|^2\right\}\right)^{1/2} \le 
C (\Delta_N)^{\min\{2 - 2\left|\alpha\right|,\ 
3/2 - \left|\alpha\right| - \max\limits_{1 \le j \le m} \left|\beta_j\right|\}}+
$$
$$
+
C (\Delta_N)^{\min\{1 - \varepsilon,\ 2 - 2\left|\alpha\right|\}} + 
C (\Delta_N)^{\min\{1/2 - \varepsilon,\ 3/2 - \left|\alpha\right| - 
\max\limits_{1 \le j \le m} \left|\beta_j\right|\}}\le
$$
\begin{equation}
\label{setka27}
\le
C (\Delta_N)^{\min\{1/2 - \varepsilon,\ 2 - 2\left|\alpha\right|,\ 
3/2 - \left|\alpha\right| - \max\limits_{1 \le j \le m} \left|\beta_j\right|\}}.
\end{equation}

\vspace{2mm}

Let us estimate $\left({\sf M}\left\{\bigl|G_t^{(k)}\bigr|^2\right\}\right)^{1/2}.$
Replacing $\tau_l$ with $s$ in (\ref{setka3}), we obtain
\begin{equation}
\label{setka303x}
s^{-2\left|\beta_j\right|} \le C s^{-2\max\limits_{1 \le j \le m} \left|\beta_j\right|}\ \ \ (0 < s/T \le 1).
\end{equation}

Then (see (\ref{setka30})), 
\begin{equation}
\label{setka303}
\left({\sf M}\left\{\bigl|G_t^{(k)}\bigr|^2\right\}\right)^{1/2} \le C 
\sum_{j=1}^m \left( \sum_{l=0}^{n_t} \int\limits_{\tau_l}^{\tau_{l+1}} 
s^{-2\max\limits_{1 \le j \le m} \left|\beta_j\right|} 
{\sf M}\bigl\{\left\Vert{\bf x}_{\tau_l} - {\bf x}_s\right\Vert^2\bigr\} ds \right)^{1/2}.
\end{equation}

Under the square root in (\ref{setka303}), for $l = 0$ (see (\ref{1983})):
$$
\int\limits_0^{\Delta_N} s^{-2\max\limits_{1 \le j \le m} \left|\beta_j\right|} 
{\sf M}\bigl\{\left\Vert {\bf x}_0 - {\bf x}_s\right\Vert^2\bigr\} ds \le
$$
$$
\le C (\Delta_N)^{1 - 2\max\limits_{1 \le j \le m} \left|\beta_j\right| + 
2\min\{1 - \left|\alpha\right|,\ 1/2 - \max\limits_{1 \le j \le m} \left|\beta_j\right|\}} 
=
$$
\begin{equation}
\label{setka28}
= C (\Delta_N)^{\min\{3 - 2\left|\alpha\right| - 
2\max\limits_{1 \le j \le m} \left|\beta_j\right|,\ 2 - 4\max\limits_{1 \le j \le m} 
\left|\beta_j\right|\}}.
\end{equation}

\vspace{2mm}

Under the square root in (\ref{setka303}), 
for $l \ge 1$, using (see (\ref{setka2}))
$$
\int\limits_{\tau_l}^{\tau_{l+1}} s^{-2\max\limits_{1\le j \le m}\left|\beta_j\right|} ds \le 
l^{-2\max_j \left|\beta_j\right|} (\Delta_N)^{1 - 2\max\limits_{1\le j \le m} \left|\beta_j\right|},
$$ 
Lemma~2, and Lemma~3, we get
$$
\sum_{l=1}^{n_t} 
\int\limits_{\tau_l}^{\tau_{l+1}} 
s^{-2\max\limits_{1 \le j \le m} \left|\beta_j\right|} 
{\sf M}\bigl\{\left\Vert{\bf x}_{\tau_l} - {\bf x}_s\right\Vert^2\bigr\}\le 
$$
$$
\le C (\Delta_N)^{3 - 2\left|\alpha\right| - 2\max\limits_{1 \le j \le m} \left|\beta_j\right|} 
\sum_{l=1}^{N-1} l^{-2\left|\alpha\right| - 2\max\limits_{1 \le j \le m} \left|\beta_j\right|} +
$$
$$
+
C (\Delta_N)^{2 - 4\max\limits_{1 \le j \le m} \left|\beta_j\right|} 
\sum_{l=1}^{N-1} l^{-4\max\limits_{1 \le j \le m} \left|\beta_j\right|} \le
$$
\begin{equation}
\label{setka35}
\le C (\Delta_N)^{\min\{2 - \varepsilon,\
3 - 2\left|\alpha\right| - 2\max\limits_{1 \le j \le m} \left|\beta_j\right|\}} + 
C (\Delta_N)^{\min\{1 - \varepsilon,\ 2 - 4\max\limits_{1 \le j \le m} \left|\beta_j\right|\}}.
\end{equation}

\vspace{2mm}

Combining (\ref{setka303}), (\ref{setka28}), and (\ref{setka35}), we obtain ($\Delta_N < 1$)
\begin{equation}
\label{setka36}
\left({\sf M}\left\{\bigl|G_t^{(k)}\bigr|^2\right\}\right)^{1/2} \le  
C (\Delta_N)^{\min\{1/2 - \varepsilon/2,\ 1 - 2\max\limits_{1 \le j \le m} 
\left|\beta_j\right|,\ 3/2 - \left|\alpha\right| - \max\limits_{1 \le j \le m} 
\left|\beta_j\right|\}}.
\end{equation}

\vspace{2mm}

Combining (\ref{888.14}), (\ref{888.15}), (\ref{useful11}), 
and (\ref{setka27}), (\ref{setka36}), (\ref{comb3}), (\ref{comb4}), 
and taking into account condition (\ref{888.06}), 
we arrive at the inequality
$$
Z_u \le C (\Delta_N)^{2\min\{1/2 - \varepsilon,\ 2 - 2\left|\alpha\right|,\ 
1 - 2\max\limits_{1 \le j \le m} \left|\beta_j\right|,\
3/2 - \left|\alpha\right| - \max\limits_{1 \le j \le m} \left|\beta_j\right|\}} +
$$
$$
+ C \int\limits_0^u \left(|\sigma_0(s)| + \sum_{j=1}^m \sigma_j^2(s)\right) Z_s ds, \ \ \ u \in [0, T],
$$

\noindent
from which, using Gronwall's inequality \cite{webb}
and Lyapunov's 
inequality, we find
\begin{equation}
\label{setka40}
{\sf M}\bigl\{\left\Vert{\bf x}_T - {\bf y}_N\right\Vert\bigr\} \le C 
(\Delta_N)^{\min\{1/2 - \varepsilon,\ 2 - 2\left|\alpha\right|,\  
1 - 2\max\limits_{1 \le j \le m} \left|\beta_j\right|,\ 
3/2 - \left|\alpha\right| - \max\limits_{1 \le j \le m} \left|\beta_j\right|\}}.
\end{equation}

Observe that 
$$
3 - 2\left|\alpha\right| - 2\max\limits_{1 \le j \le m} \left|\beta_j\right| = 
(2 - 2|\alpha|) + \left(1 - 2\max\limits_{1 \le j \le m} |\beta_j|\right) \ge
$$
$$
\ge 2\min\left\{2 - 2|\alpha|,\ 1 - 2\max_{1 \le j \le m} |\beta_j|\right\},
$$
i.e.
$$
\frac{3}{2} - \left|\alpha\right| - \max\limits_{1 \le j \le m} \left|\beta_j\right| 
\ge
\min\left\{2 - 2|\alpha|,\ 1 - 2\max_{1 \le j \le m} |\beta_j|\right\}.
$$

Then, we have 
$$
\min\left\{\frac{1}{2} - \varepsilon,\ 2 - 2\left|\alpha\right|,\ 
1 - 2\max\limits_{1 \le j \le m} \left|\beta_j\right|,\ 
\frac{3}{2} - \left|\alpha\right| - \max\limits_{1 \le j \le m} \left|\beta_j\right|\right\} 
=
$$
$$
= \min\left\{\frac{1}{2} - \varepsilon,\ 
\min\left\{2 - 2\left|\alpha\right|,\ 1 - 2\max\limits_{1 \le j \le m} \left|\beta_j\right|\right\},\ 
\frac{3}{2} - \left|\alpha\right| - \max\limits_{1 \le j \le m} \left|\beta_j\right|\right\}=
$$
\begin{equation}
\label{setka41}
= \min\left\{\frac{1}{2} - \varepsilon,\ 2 - 2\left|\alpha\right|,\ 
1 - 2\max_{1 \le j \le m} \left|\beta_j\right|\right\},
\end{equation}

From (\ref{setka40}) and (\ref{setka41}) we obtain (\ref{setka50}). Theorem~6 is proved. $\Box$

\subsection{Proof of Theorem 8}

The proof of Theorem~4 still holds for the most part
in proving Theorem~8; we only obtain some new estimates.

Throughout this proof, we assume 
that $\varepsilon > 0$ is an arbitrarily small positive real number, 
and $C$ denotes various positive constants independent of $\Delta_N$, $l$, $s\in [\tau_l, \tau_{l+1}]$.

As before, the incomplete last subinterval 
$[\tau_{n_t}, t] \subseteq [\tau_{n_t}, \tau_{n_t+1}]$ ($n_t \le N - 1$) 
is included in the sums as the term with $l = n_t$, 
and $\sigma_j^2(s)$, $\sigma_i^2(\theta)$ are bounded 
above by $C s^{-2\max\limits_{1 \le j \le m} \left|\beta_j\right|}$ and 
$C \theta^{-2\max\limits_{1 \le j \le m} \left|\beta_j\right|}$
(see (\ref{setka303x})).

To prove Theorem~8, we need to derive new estimates for
$\left({\sf M}\left\{\bigl|S_t^{(k)}\bigr|^2\right\}\right)^{1/2},$
$\left({\sf M}\left\{\bigl|V_t^{(k)}\bigr|^2\right\}\right)^{1/2},$
$\left({\sf M}\left\{\bigl|D_t^{(k)}\bigr|^2\right\}\right)^{1/2}$
(see (\ref{setka400}), (\ref{setka401}), (\ref{one5}), (\ref{setka402}), 
(\ref{one1})).

Let us estimate $\left({\sf M}\left\{\bigl|S_t^{(k)}\bigr|^2\right\}\right)^{1/2}.$
We have (see (\ref{setka405}))
$$
\left({\sf M}\left\{\bigl|S_t^{(k)}\bigr|^2\right\}\right)^{1/2} \le 
$$
\begin{equation}
\label{setka406}
\le C \sum_{i,j=1}^m \left( \sum_{l=0}^{n_t} \int\limits_{\tau_l}^{\tau_{l+1}} 
s^{-2\max\limits_{1 \le j \le m} \left|\beta_j\right|} \int\limits_{\tau_l}^s 
\theta^{-2\max_{1 \le j \le m} \left|\beta_j\right|} 
{\sf M}\bigl\{\left\Vert {\bf x}_{\tau_l}-{\bf x}_{\theta}\right\Vert^2\bigr\} 
d\theta ds \right)^{1/2}.
\end{equation}

Under the square root in (\ref{setka406}), for $l = 0$ (see (\ref{1983})):
$$
\int\limits_0^{\Delta_N} s^{-2\max\limits_{1 \le j \le m} \left|\beta_j\right|} 
\int\limits_0^s \theta^{-2\max\limits_{1 \le j \le m} 
\left|\beta_j\right|} {\sf M}\bigl\{\left\Vert {\bf x}_{0}-{\bf x}_{\theta}\right\Vert^2\bigr\}  
d\theta ds \le 
$$
$$
\le C (\Delta_N)^{2 - 4\max\limits_{1 \le j \le m} \left|\beta_j\right| + 
2\min\{1 - \left|\alpha\right|,\ 
1/2 - \max\limits_{1 \le j \le m} \left|\beta_j\right|\}} =
$$
\begin{equation}
\label{setka407}
= C (\Delta_N)^{\min\{4 - 2\left|\alpha\right| - 4\max\limits_{1 \le j \le m} \left|\beta_j\right|,\
3 - 6\max\limits_{1 \le j \le m} \left|\beta_j\right|\}}.
\end{equation}

\vspace{2mm}

Under the square root in (\ref{setka406}), for $l \ge 1$, using the relations
\begin{equation}
\label{setka510}
\int\limits_{\tau_l}^{\tau_{l+1}} \int\limits_{\tau_l}^s d\theta ds = \frac{(\Delta_N)^2}{2},\ \ \ 
s^{-2\max\limits_{1 \le j \le m} \left|\beta_j\right|} \le 
(\tau_l)^{-2\max\limits_{1 \le j \le m} \left|\beta_j\right|}\ \
(s\in [\tau_l,\tau_{l+1}]),
\end{equation}
Lemma~2 and Lemma~3, we obtain
$$
\sum_{l=1}^{n_t} \int\limits_{\tau_l}^{\tau_{l+1}} s^{-2\max\limits_{1 \le j \le m} 
\left|\beta_q\right|} 
\int\limits_{\tau_l}^s \theta^{-2\max\limits_{1 \le j \le m} \left|\beta_j\right|} 
{\sf M}\bigl\{\left\Vert {\bf x}_{\tau_l}-{\bf x}_{\theta}\right\Vert^2\bigr\} 
d\theta ds \le 
$$
$$
\le C (\Delta_N)^{4 - 2\left|\alpha\right| - 4\max\limits_{1 \le j \le m} \left|\beta_j\right|} 
\sum_{l=1}^{N-1} l^{-2\left|\alpha\right| - 4\max\limits_{1 \le j \le m} \left|\beta_j\right|} +
$$
$$
\quad + C (\Delta_N)^{3 - 6\max\limits_{1 \le j \le m} \left|\beta_j\right|} 
\sum_{l=1}^{N-1} l^{-6\max\limits_{1 \le j \le m} \left|\beta_j\right|} \le
$$
\begin{equation}
\label{setka408}
\le C (\Delta_N)^{\min\{3 - \varepsilon,\ 
4 - 2\left|\alpha\right| - 4\max\limits_{1 \le j \le m} \left|\beta_j\right|\}} + 
C (\Delta_N)^{\min\{2 - \varepsilon,\ 3 - 6\max\limits_{1 \le j \le m} \left|\beta_j\right|\}}.
\end{equation}

\vspace{2mm}

Combining (\ref{setka406}), (\ref{setka407}), and (\ref{setka408}), we get
$$
\left({\sf M}\left\{\bigl|S_t^{(k)}\bigr|^2\right\}\right)^{1/2} \le 
$$
\begin{equation}
\label{setka2000}
\le C (\Delta_N)^{\min\{1 - \varepsilon/2,\ 3/2 - 3\max\limits_{1 \le j \le m} |\beta_j|,\ 
2 - |\alpha| - 2\max\limits_{1 \le j \le m} |\beta_j|\}}.
\end{equation}

\vspace{2mm}

Let us estimate $\left({\sf M}\left\{\bigl|V_t^{(k)}\bigr|^2\right\}\right)^{1/2}.$
We have (see (\ref{setka500}))
$$
\left({\sf M}\left\{\bigl|V_t^{(k)}\bigr|^2\right\}\right)^{1/2}\le 
$$
\begin{equation}
\label{setka501}
\le C(\Delta_N)^{1/2}
\sum\limits_{i,j=1}^m \Biggl(\sum_{l=0}^{n_t}\int\limits_{\tau_l}^{\tau_{l+1}}
s^{-2\max\limits_{1 \le j \le m} \left|\beta_j\right|} ds\int\limits_{\tau_l}^{\tau_{l+1}}
\theta^{-2\max\limits_{1 \le j \le m}\left|\beta_j\right|}  
d\theta 
\Biggr)^{1/2}.
\end{equation}

Under the square root in (\ref{setka501}), for $l = 0$:
\begin{equation}
\label{setka520}
\int\limits_{0}^{\Delta_N}
s^{-2\max\limits_{1 \le j \le m} \left|\beta_j\right|} ds\int\limits_{0}^{\Delta_N}
\theta^{-2\max\limits_{1 \le j \le m}\left|\beta_j\right|}  
d\theta \le
C (\Delta_N)^{2 - 4\max\limits_{1 \le j \le m} \left|\beta_j\right|},
\end{equation}

\noindent
and for $l \ge 1$ (see (\ref{setka510}) and Lemma~3): 
$$
\sum_{l=1}^{n_t}\int\limits_{\tau_l}^{\tau_{l+1}}
s^{-2\max\limits_{1 \le j \le m} \left|\beta_j\right|} ds\int\limits_{\tau_l}^{\tau_{l+1}}
\theta^{-2\max\limits_{1 \le j \le m}\left|\beta_j\right|} 
d\theta \le
$$
\begin{equation}
\label{setka521}
\le C (\Delta_N)^{2 - 4\max\limits_{1 \le j \le m} \left|\beta_j\right|} 
\sum_{l=1}^{N-1} l^{-4\max\limits_{1 \le j \le m} \left|\beta_j\right|} 
\le C (\Delta_N)^{\min\{1 - \varepsilon,\ 2 - 4\max\limits_{1 \le j \le m} 
\left|\beta_j\right|\}}.
\end{equation}

\vspace{2mm}

Combining (\ref{setka501}), (\ref{setka520}), and (\ref{setka521}), we obtain 
\begin{equation}
\label{setka2001}
\left({\sf M}\left\{\bigl|V_t^{(k)}\bigr|^2\right\}\right)^{1/2} \le 
C (\Delta_N)^{\min\{1 - \varepsilon/2,\ 3/2 - 2\max\limits_{1 \le j \le m} \left|\beta_j\right|\}}.
\end{equation}

\vspace{1mm}

Let us estimate the terms on the right-hand side of the following inequality (see (\ref{one1}))
\begin{equation}
\label{setka530}
\left({\sf M}\left\{\bigl|D_t^{(k)}\bigr|^2\right\}\right)^{1/2}\le
\sum_{i=1}^7 \left({\sf M}\left\{\left(\sum_{l=0}^{n_t}R^{(k)}_{(i)\tau_{l+1},\tau_l}
\right)^2\right\}\right)^{1/2}
\end{equation}

\noindent
for $i=7,6,4,3,2$ (the estimates (\ref{one3005}), (\ref{one3000}) 
($i=5,$ $i=1$) will be used further without 
modification).

Using (\ref{setka600}), we have
\begin{equation}
\label{setka910}
\left({\sf M}\left\{\left(\sum_{l=0}^{n_t}R^{(k)}_{(7)\tau_{l+1},\tau_l}
\right)^2\right\}\right)^{1/2}\le
C\sum\limits_{i,j=1}^m\Biggl(\sum_{l=0}^{n_t}
\Biggl(\int\limits_{\tau_l}^{\tau_{l+1}}
s^{-2\max\limits_{1 \le j \le m} \left|\beta_j\right|}ds\Biggr)^3
\Biggr)^{1/2}.
\end{equation}

Similarly, we obtain (see (\ref{setka800}), (\ref{one101}), (\ref{setka700}), and (\ref{setka900}))
\begin{equation}
\label{setka911}
\left({\sf M}\left\{\left(\sum_{l=0}^{n_t}R^{(k)}_{(i)\tau_{l+1},\tau_l}
\right)^2\right\}\right)^{\hspace{-2mm}1/2}
\hspace{-4mm}\le
C\sum\limits_{j=1}^m\Biggl(\sum_{l=0}^{n_t}
\int\limits_{\tau_l}^{\tau_{l+1}}
s^{-2\max\limits_{1 \le j \le m} \left|\beta_j\right|}ds
\Biggl(\int\limits_{\tau_l}^{\tau_{l+1}}
s^{-\left|\alpha\right|}ds\Biggr)^2
\Biggr)^{1/2},
\end{equation}
\begin{equation}
\label{setka913}
\left({\sf M}\left\{\left(\sum_{l=0}^{n_t}R^{(k)}_{(3)\tau_{l+1},\tau_l}
\right)^2\right\}\right)^{\hspace{-1mm}1/2}
\le C\sum\limits_{j=1}^m\sum_{l=0}^{n_t}
\int\limits_{\tau_l}^{\tau_{l+1}} s^{-\left|\alpha\right|}ds
\int\limits_{\tau_l}^{\tau_{l+1}}
\theta^{-2\max\limits_{1 \le j \le m} \left|\beta_j\right|} d\theta,
\end{equation}
\begin{equation}
\label{setka914}
\left({\sf M}\left\{\left(\sum_{l=0}^{n_t}R^{(k)}_{(2)\tau_{l+1},\tau_l}
\right)^2\right\}\right)^{\hspace{-1mm}1/2}\le
C\sum_{l=0}^{n_t}
\Biggl(\int\limits_{\tau_l}^{\tau_{l+1}}
s^{-\left|\alpha\right|}ds\Biggr)^2,
\end{equation}
where $i=4$ and $i=6$ in (\ref{setka911}).

Splitting the cases $l = 0$ and $l \ge 1$ on the right-hand sides
of (\ref{setka910})--(\ref{setka914}), using the inequality in (\ref{setka510})
and its analogue with $\left|\alpha\right|$ instead of 
$2\max\limits_{1 \le j \le m} \left|\beta_j\right|$,
and applying Lemma~3 for $\Delta_N < 1$ (as before),
we obtain the following estimates:
\begin{equation}
\label{setka930}
\left({\sf M}\left\{\left(\sum_{l=0}^{n_t}R^{(k)}_{(7)\tau_{l+1},\tau_l}
\right)^2\right\}\right)^{1/2}\le
C (\Delta_N)^{\min\{1 - \varepsilon/2,\ 3/2 - 3\max\limits_{1 \le j \le m} \left|\beta_j\right|\}},
\end{equation}
\begin{equation}
\label{setka931}
\left({\sf M}\left\{\left(\sum_{l=0}^{n_t}R^{(k)}_{(i)\tau_{l+1},\tau_l}
\right)^2\right\}\right)^{1/2}\le
C (\Delta_N)^{\min\{1 - \varepsilon/2,\ 
3/2 - \left|\alpha\right| - \max\limits_{1 \le j \le m} \left|\beta_j\right|\}} \ \ \ (i = 4,\ 6),
\end{equation}
\begin{equation}
\label{setka932}
\left({\sf M}\left\{\left(\sum_{l=0}^{n_t}R^{(k)}_{(3)\tau_{l+1},\tau_l}
\right)^2\right\}\right)^{1/2}\le
C (\Delta_N)^{\min\{1 - \varepsilon,\ 2 - \left|\alpha\right| - 
2\max\limits_{1 \le j \le m} \left|\beta_j\right|\}},
\end{equation}
\begin{equation}
\label{setka933}
\left({\sf M}\left\{\left(\sum_{l=0}^{n_t}R^{(k)}_{(2)\tau_{l+1},\tau_l}
\right)^2\right\}\right)^{1/2}\le
C (\Delta_N)^{\min\{1 - \varepsilon,\ 2 - 2\left|\alpha\right|\}}.
\end{equation}

Moreover (see (\ref{one3005}), (\ref{one3000})), 
\begin{equation}
\label{setka934}
\left({\sf M}\left\{\left(\sum_{l=0}^{n_t}R^{(k)}_{(i)\tau_{l+1},\tau_l}
\right)^2\right\}\right)^{1/2}\le
C \Delta_N\ \ \ (i=1,\ 5).
\end{equation}

Combining (\ref{setka530}), (\ref{setka930})--(\ref{setka934}), we get
\begin{equation}
\label{setka2002}
\left({\sf M}\left\{\bigl|D_t^{(k)}\bigr|^2\right\}\right)^{1/2}\le
C (\Delta_N)^{\min\{1 - \varepsilon,\ 2 - 2\left|\alpha\right|,\ 
3/2 - 3\max\limits_{1 \le j \le m} \left|\beta_j\right|,\ 
3/2 - \left|\alpha\right| - \max\limits_{1 \le j \le m} \left|\beta_j\right|\}},
\end{equation}

\noindent
where we used the fact that
$$
3/2 - \left|\alpha\right| - \max\limits_{1 \le j \le m}\left|\beta_j\right|\le
(3/2 - \left|\alpha\right| - \max\limits_{1 \le j \le m}\left|\beta_j\right|)
+ (1/2-\max\limits_{1 \le j \le m}\left|\beta_j\right|)=
$$
$$
=2 - \left|\alpha\right| - 2\max\limits_{1 \le j \le m}\left|\beta_j\right|.
$$

Combining (\ref{comb30}), (\ref{comb40}), (\ref{comb45}), 
and (\ref{setka2000}),  (\ref{setka2001}), 
(\ref{setka2002}), 
and taking into account condition (\ref{Mil12}), 
we arrive at the inequality
$$
Z_u \le C (\Delta_N)^{2\min\{1 - \varepsilon,\ 2 - 2\left|\alpha\right|,\ 
3/2 - 3\max\limits_{1 \le j \le m} \left|\beta_j\right|,\
3/2 - \left|\alpha\right| - \max\limits_{1 \le j \le m} \left|\beta_j\right|\}} +
$$
$$
+ C \int\limits_0^u \left(|\sigma_0(s)| + \sum_{j=1}^m \sigma_j^2(s)\right) Z_s ds, \ \ \ u \in [0, T],
$$

\noindent
from which, using Gronwall's inequality \cite{webb}
and Lyapunov's 
inequality, we find
$$
{\sf M}\bigl\{\left\Vert{\bf x}_T - {\bf y}_N\right\Vert\bigr\} 
\le C (\Delta_N)^{\min\{1 - \varepsilon,\ 2 - 2\left|\alpha\right|,\ 
3/2 - 3\max\limits_{1 \le j \le m} \left|\beta_j\right|,\
3/2 - \left|\alpha\right| - \max\limits_{1 \le j \le m} \left|\beta_j\right|\}}.
$$

Theorem~8 is proved. $\Box$

\section{Numerical Experiments}

\subsection{The Case of It\^{o} SDE with Non-Differentiable 
Diffusion Coefficient}

{\bf Example~1.}\ Consider the following scalar It\^{o} SDE
\begin{equation}
\label{ex001}
x_t=x_0+\int\limits_0^t ax_s ds+\int\limits_0^t b\left|\sin(As)\right|x_s dw_s,\ \ \ x_0=x(0),\ \ \ t\in[0, T],
\end{equation}
where $w_t$ is a scalar standard Wiener process, $a, b, A\in {\bf R}.$

The exact solution to equation (\ref{ex001}) has the form
\begin{equation}
\label{ww1}
x_t=x_0\exp\left(at-\frac{1}{2}b^2\int\limits_0^t \left|\sin(As)\right|^2 ds+b\int\limits_0^t
\left|\sin(As)\right|dw_s\right).
\end{equation}

Note that 
\begin{equation}
\label{ww2}
\int\limits_{\tau}^t \left|\sin(As)\right|^2 ds=\frac{1}{2} (t-\tau)-\frac{1}{4A}
\left(\sin(2At)-\sin(2A\tau)\right)\stackrel{\sf def}{=}\sigma_{t,\tau}^2(A).
\end{equation}

Let $\tau_p=p\Delta$ $(p=0,1,\ldots,N),$
$\Delta=T /N$ $(N>1)$ is a constant step of integration.
Then, from (\ref{ww1}) and (\ref{ww2}) we obtain
\begin{equation}
\label{ww3}
x_{T}=x_0\exp\left(aT-\frac{1}{2}b^2 \sigma_{T,0}^2(A)
+b\sum\limits_{p=0}^{N-1}\sigma_{\tau_{p+1},\tau_p}(A)\zeta_{p+1,p}\right),
\end{equation}
where
\begin{equation}
\label{ww4}
\zeta_{p+1,p}=\frac{1}{\sigma_{\tau_{p+1},\tau_p}(A)}
\int\limits_{\tau_p}^{\tau_{p+1}}
\left|\sin(As)\right|dw_s\ \ \ (p=0, 1,\ldots,N-1)
\end{equation}
are independent standard Gaussian random variables.

It is obvious that the commutativity condition (\ref{2000}) holds for the
It\^{o} SDE (\ref{ex001}). Then, the new modification (\ref{200}) of Milstein method 
for the It\^{o} SDE (\ref{ex001}) takes the form
\begin{equation}
\label{ww4z}
y_{p+1}=y_p\biggl(1+a\Delta+b\sigma_{\tau_{p+1},\tau_p}(A)\zeta_{p+1,p}
+\frac{1}{2}b^2 \sigma_{\tau_{p+1},\tau_p}^2(A)\left(\left(\zeta_{p+1,p}\right)^2-1\right)\biggr),
\end{equation}
where $p=0,1,\ldots,N-1,$ $y_{\tau_p}\stackrel{\sf def}{=}y_p$ 
$(p=0,1,\ldots,N)$,  $\tau_N=T.$

{\bf Numerical Experiment 1.}\ {\it For the It\^{o}
SDE {\rm (\ref{ex001})}
simulate {\rm 2000} {\rm (}$M=20$ groups of $L=100$ realizations each{\rm )} 
independent realizations $x_T^{(k,j)}$ of the random variable $x_T$ at $T=1,$ $x_0=1$  and $a=-0.01,$ $b=1,$
$A=2^6\pi$ 
by the formula {\rm (\ref{ww3})}
and the corresponding independent realizations $y_T^{(k,j)}$ 
of the random variable $y_T$ by the formula {\rm (\ref{ww4z})} 
{\rm ($k=1,\ldots,L,\ j=1,\ldots,M$)}.
Here $N=T/\Delta,$ $\Delta=2^{-15}.$ 
Calculate the estimate $\hat\varepsilon_{100,20}$ of the error 
$\varepsilon={\sf M}\left\{|x_{T}-y_T|\right\}$ 
by the formula 
$$
\hat\varepsilon_{L,M}=\frac{1}{LM}\sum_{k=1}^L\sum_{j=1}^M
\left|x_T^{(k,j)}-y_T^{(k,j)}\right|
$$
and a $90\%$ 
confidence interval for $\varepsilon,$
where $x_T^{(k,j)}$ and $y_T^{(k,j)}$
for a fixed pair $(k,j)$ correspond to
the same realization 
of the Wiener process $w_t,$ $t\in[0,T].$
Repeat the calculations for 
$\Delta=2^{-j},$\ $j=14, 13,\ldots, 3$ and plot 
the dependence $\hat\varepsilon_{100,20}(\Delta)$ with $90\%$ 
confidence intervals for $\varepsilon$ at the specified values 
of ${\Delta.}^1$ At the same time$,$ when doubling
the integration step$,$ use the formula
$$
\sigma_{\tau_{p+2},\tau_p}(A)\zeta_{p+2,p}=
\sigma_{\tau_{p+1},\tau_p}(A)\zeta_{p+1,p}+\sigma_{\tau_{p+2},\tau_{p+1}}(A)\zeta_{p+2,p+1}\ \ \ 
w.~p.~{\rm 1},
$$
i.e. w.~p.~{\rm 1}
$$
\int\limits_{\tau_p}^{\tau_{p+2}}
\left|\sin(As)\right|dw_s=\int\limits_{\tau_p}^{\tau_{p+1}}
\left|\sin(As)\right|dw_s+\int\limits_{\tau_{p+1}}^{\tau_{p+2}}
\left|\sin(As)\right|dw_s.
$$

\footnotetext[1]{Here, we apply
the same methodology to construct confidence intervals
as in \cite{2} (Sect.~9.3).}
}

The result of Numerical Experiment~1 is presented in Fig.~1. 
In this figure, we can see that the new modification (\ref{200})
of the Milstein method has order approximately 1.0 of strong convergence
(see Remark~1).

Next, we compare the classical Milstein method (\ref{11vvv}) and the new modification
(\ref{200}) of the Milstein method using the numerical integration
of the It\^{o} SDE (\ref{ex001}) as an example.

Application of the numerical scheme (\ref{11vvv}) to the It\^{o}
SDE (\ref{ex001}) yields the following result
\begin{equation}
\label{ww5}
\tilde y_{p+1}=\tilde y_p\biggl(1+a\Delta+b\left|sin(A\tau_p)\right|\Delta w_p
+\frac{1}{2}b^2\left|sin(A\tau_p)\right|^2 \left(\left(\Delta w_p\right)^2
-\Delta\right)\biggr),
\end{equation}
where $\Delta w_p=w_{\tau_{p+1}}-w_{\tau_p}$.

Denote
$$
\int\limits_{\tau_p}^{\tau_{p+1}}
\left|\sin(As)\right|dw_s\stackrel{\sf def}{=}J_{\tau_{p+1},\tau_p}.
$$

Moreover, the equalities (\ref{ww3}), (\ref{ww4z}) can be written in the form
\begin{equation}
\label{ww6}
x_{T}=x_0\exp\left(aT-\frac{1}{2}b^2 \sigma_{T,0}^2(A)
+b\sum\limits_{p=0}^{N-1}J_{\tau_{p+1},\tau_p}\right),
\end{equation}
\begin{equation}
\label{ww7}
y_{p+1}=y_p\biggl(1+a\Delta+b J_{\tau_{p+1},\tau_p}
+\frac{1}{2}b^2 \left(\left(J_{\tau_{p+1},\tau_p}\right)^2-\sigma_{\tau_{p+1},\tau_p}^2(A)
\right)\biggr),
\end{equation}
where $\sigma_{\tau_{p+1},\tau_p}^2(A)$ is defined by (\ref{ww2}).

Recall that $\tau_p=p\Delta,$ i.e. $\tau_{p+1}-\tau_p=\Delta.$ It is not difficult to show that
\begin{equation}
\label{701a}
\Delta w_p=\sqrt{\Delta}\cdot \xi_{p+1,p},
\end{equation}
\begin{equation}
\label{702}
J_{\tau_{p+1},\tau_p}=C^{(1)}_{p+1,p}\cdot \eta_{p+1,p}+ C^{(2)}_{p+1,p}\cdot \xi_{p+1,p}, 
\end{equation}
where 
\begin{equation}
\label{9000}
C^{(1)}_{p+1,p}=\sqrt{\int\limits_{\tau_p}^{\tau_{p+1}} \left|\sin(As)\right|^2 ds-
\frac{1}{\Delta}\left(
\int\limits_{\tau_p}^{\tau_{p+1}} \left|\sin(As)\right| ds\right)^2},
\end{equation}
\begin{equation}
\label{700a}
C^{(2)}_{p+1,p}=\frac{1}{\sqrt{\Delta}}\int\limits_{\tau_p}^{\tau_{p+1}} \left|\sin(As)\right| ds,
\end{equation}
and $\xi_{p+1,p},$ $\eta_{p+1,p}$ $(p=0,1,\ldots,N-1)$ are mutually independent
standard Gaussian random variables.

It should be noted that the integral in (\ref{700a}) can be computed
using standard computational procedures.

Let us substitute (\ref{701a}) and (\ref{702}) into 
(\ref{ww5}), (\ref{ww6}), and (\ref{ww7})
\begin{equation}
\label{ww5a}
\tilde y_{p+1}=\tilde y_p\biggl(1+a\Delta+b\sqrt{\Delta}\left|sin(A\tau_p)\right|\xi_{p+1,p},
+\frac{1}{2}b^2\Delta \left|sin(A\tau_p)\right|^2 \left(\left(\xi_{p+1,p}\right)^2
-1\right)\biggr),
\end{equation}
\begin{equation}
\label{ww6a}
x_{T}=x_0\exp\left(aT-\frac{1}{2}b^2 \sigma_{T,0}^2(A)
+b\sum\limits_{p=0}^{N-1}\left(C^{(1)}_{p+1,p}\cdot \eta_{p+1,p}+ C^{(2)}_{p+1,p}\cdot \xi_{p+1,p}\right)
\right),
\end{equation}
$$
y_{p+1}=y_p\biggl(1+a\Delta+b \left(C^{(1)}_{p+1,p}\cdot \eta_{p+1,p}+ 
C^{(2)}_{p+1,p}\cdot \xi_{p+1,p}\right)+\biggr.
$$
\begin{equation}
\label{ww7a}
\biggl.+\frac{1}{2}b^2 \left(
\left(C^{(1)}_{p+1,p}\cdot \eta_{p+1,p}+ C^{(2)}_{p+1,p}\cdot \xi_{p+1,p}\right)^2-
\sigma_{\tau_{p+1},\tau_p}^2(A)\right)
\biggr),
\end{equation}

\vspace{1mm}
\noindent
where $\sigma_{\tau_{p+1},\tau_p}^2(A),$ 
$C^{(1)}_{p+1,p},$ $C^{(2)}_{p+1,p}$ are defined by 
(\ref{ww2}), (\ref{9000}), (\ref{700a}), 
and $\xi_{p+1,p},$ $\eta_{p+1,p}$ $(p=0,1,\ldots,N-1)$ are mutually independent
standard Gaussian random variables.

{\bf Numerical Experiment 2.}\ {\it For the It\^{o}
SDE {\rm (\ref{ex001})}
simulate {\rm 2000} {\rm (}$M=20$ groups of $L=100$ realizations each{\rm )} 
independent realizations $x_T^{(k,j)}$ of the random variable $x_T$ at $T=1,$ $x_0=1$ and $a=-0.01,$ $b=1,$
$A=2^6\pi$ 
by the formula {\rm (\ref{ww6a})}
and the corresponding independent realizations $\tilde y_T^{(k,j)}$ and $y_T^{(k,j)}$ 
of the random variables $\tilde y_T$ and $y_T$ by the formulas {\rm (\ref{ww5a})}
and {\rm (\ref{ww7a})}, correspondently
{\rm ($k=1,\ldots,L,\ j=1,\ldots,M$)}.
Here $N=T/\Delta,$ $\Delta=2^{-15}.$ 
Calculate the estimates $\hat\varepsilon_{100,20}^{(1)},$ $\hat\varepsilon_{100,20}^{(2)}$ of the errors 
$\varepsilon^{(1)}={\sf M}\left\{|x_{T}-\tilde y_T|\right\},$ 
$\varepsilon^{(2)}={\sf M}\left\{|x_{T}-y_T|\right\}$ 
by the formulas 
\begin{equation}
\label{a5000}
\hat\varepsilon_{L,M}^{(1)}=\frac{1}{LM}\sum_{k=1}^L\sum_{j=1}^M
\left|x_T^{(k,j)}-\tilde y_T^{(k,j)}\right|,\ \ \ 
\hat\varepsilon_{L,M}^{(2)}=\frac{1}{LM}\sum_{k=1}^L\sum_{j=1}^M
\left|x_T^{(k,j)}-y_T^{(k,j)}\right|,
\end{equation}

\noindent
and a $90\%$ 
confidence intervals for $\varepsilon^{(1)},$ $\varepsilon^{(2)},$
where $x_T^{(k,j)},$ $\tilde y_T^{(k,j)},$ and $y_T^{(k,j)}$
for a fixed pair $(k,j)$ correspond to
the same realization 
of the Wiener process $w_t,$ $t\in[0,T].$
Repeat the calculations for 
$\Delta=2^{-j},$\ $j=14, 13,\ldots, 3$ and plot 
the dependences $\hat\varepsilon_{100,20}^{(1)}(\Delta),$
$\hat\varepsilon_{100,20}^{(2)}(\Delta)$ with $90\%$ 
confidence intervals for $\varepsilon^{(1)},$ $\varepsilon^{(2)}$ at the specified values 
of ${\Delta.}$ At the same time$,$ when doubling
the integration step$,$ use the formulas
$$
J_{\tau_{p+2},\tau_p}=J_{\tau_{p+2},\tau_{p+1}}+J_{\tau_{p+1},\tau_p},
$$
$$
w_{\tau_{p+2}}-w_{\tau_p}=\Delta w_{p+1}+\Delta w_{p}
$$
w.~p.~{\rm 1}.
}

\vspace{2mm}

Fig.~2 shows that the new modification (\ref{200}) of the Milstein method 
has order approximately 1.0 of strong convergence. At the same time, the Milstein 
method (\ref{11vvv}) fails to converge in the range of large 
integration step sizes $\Delta\in [2^{-7}, 2^{-3}]$, converges with strong order approximately
1.0 for small step sizes $\Delta <2^{-7}$, yet its error significantly exceeds that of method (\ref{200}).
This is due to the fact that the sufficient conditions for 
convergence are not satisfied for the Milstein method (\ref{11vvv}), and it fails 
to converge for large integration step sizes due to the non-differentiable 
function 
$b\left|\sin(As)\right|x.$

\subsection{The Case of It\^{o} SDEs with Singular Drift and Diffusion Coefficients} 

The numerical experiments in this section demonstrate good agreement
between the experimental 
orders of convergence of numerical methods (\ref{11}) and (\ref{euler1})
and the theoretical orders of convergence 
from Theorems 9 and 7.

{\bf Example~2.}\  Consider the following scalar commutative It\^{o} SDE
\begin{equation}
\label{1000e}
x_t=x_0+\int\limits_0^t a(s)x_s ds+
\sum\limits_{j=1}^3
\int\limits_0^t b_j(s)x_s d{\bf w}_{s}^{(j)},\ \ \
x_0=x(0),\ \ \ t\in[0, T],
\end{equation}
where $a(t)=A t^{\alpha},$
$b_j(t)= C_j t^{\beta_j},$ $\alpha\in (-1,0)$, $\beta_j\in (-1/2,0)$ $(j=1,2,3)$,
$A,C_j\in {\bf R}$ $(j=1,2,3),$
${\bf w}_t^{(j)}$ $(j=1,2,3)$ are independent standard Wiener processes.

The exact solution to equation (\ref{1000e}) has the following form
\begin{equation}
\label{1001e}
x_t=x_0 \exp\left(\int\limits_0^t \left(a(s)-\frac{1}{2}\sum\limits_{j=1}^3 b_j^2(s)\right)ds+\sum\limits_{j=1}^3
\int\limits_0^t b_j(s)d{\bf w}_s^{(j)}\right).
\end{equation}

Let $\tau_p=p\Delta$ $(p=0,1,\ldots,N),$
$\Delta=T /N$ $(N>1)$ is a constant step of integration.

Then, from (\ref{1001e}) we obtain
\begin{equation}
\label{1001ee}
x_T=x_0 \exp\left(\int\limits_0^T \left(a(s)-\frac{1}{2}\sum\limits_{j=1}^3 b_j^2(s)\right)ds+\sum\limits_{j=1}^3
\sum\limits_{p=0}^{N-1}\int\limits_{\tau_p}^{\tau_{p+1}} b_j(s)d{\bf w}_s^{(j)}\right).
\end{equation}

Application of the classical Milstein method (\ref{11vvv}) to the It\^{o}
SDE (\ref{1000e}) yields the following result
$$
\tilde y_{p+1}=
\tilde y_p\Biggl(1+a(\tau_p)\Delta+\sum\limits_{j=1}^3 b_j(\tau_p)\Delta {\bf w}_p^{(j)}+\Biggr.
$$
\begin{equation}
\label{a2000}
\Biggl.+
\frac{1}{2}\sum\limits_{l,j=1}^3 b_l(\tau_p)b_j(\tau_p)\left(
\Delta {\bf w}_p^{(j)}\Delta {\bf w}_p^{(l)}-{\bf 1}_{\{l=j\}}\Delta\right)\Biggr),
\end{equation}

\noindent
where $\Delta {\bf w}_p^{(j)}={\bf w}_{\tau_{p+1}}^{(j)}-{\bf w}_{\tau_p}^{(j)}$
$(j=1,2,3).$

At the same time, the new modification (\ref{200}) of Milstein method 
for the It\^{o} SDE (\ref{1000e}) takes the form
\begin{equation}
\label{a2001}
y_{p+1}=y_{p}\Biggl(1+
I^{a}_{\tau_{p+1},\tau_p}+
\sum\limits_{j=1}^3
I^{b_j(j)}_{\tau_{p+1},\tau_p}
+\frac{1}{2}\sum\limits_{l,j=1}^3
\Biggl(I^{b_j(j)}_{\tau_{p+1},\tau_p}\cdot I^{b_l(l)}_{\tau_{p+1},\tau_p}-
{\bf 1}_{\{l=j\}}
I^{b_l^2}_{\tau_{p+1},\tau_p}\Biggr)\Biggr),
\end{equation}

\noindent
where
\begin{equation}
\label{a4000}
I^{b_l^2}_{\tau,t}=\int\limits_{\tau}^{t}b_l^2(s)ds,\ \ \ 
I^{a}_{\tau,t}=\int\limits_{\tau}^{t} a(s)ds,\ \ \ 
I^{b_j(j)}_{\tau,t}=
\int\limits_{\tau}^{t} b_j(s)d{\bf w}_{s}^{(j)}\ \ \ (j=1,2,3).
\end{equation}

Furthermore, the modified Euler's method (\ref{euler1}) for the It\^{o} SDE (\ref{1000e})
can be written as
\begin{equation}
\label{a2001xxx}
\hat y_{p+1}=\hat y_{p}\left(1+
I^{a}_{\tau_{p+1},\tau_p}+
\sum\limits_{j=1}^3
I^{b_j(j)}_{\tau_{p+1},\tau_p}\right).
\end{equation}

By analogy with (\ref{701a}), (\ref{702}) we have
\begin{equation}
\label{a1000}
\Delta {\bf w}_p^{(j)}=\sqrt{\Delta}\cdot \xi_{p+1,p}^{(j)},
\end{equation}
\begin{equation}
\label{a1001}
I^{b_j(j)}_{\tau_{p+1},\tau_p}=D^{(1,j)}_{p+1,p}\cdot \eta_{p+1,p}^{(j)}+ 
D^{(2,j)}_{p+1,p}\cdot \xi_{p+1,p}^{(j)}, 
\end{equation}
where 
\begin{equation}
\label{a1002}
D^{(1,j)}_{p+1,p}=\sqrt{\int\limits_{\tau_p}^{\tau_{p+1}} b_j^2(s) ds-
\frac{1}{\Delta}\left(
\int\limits_{\tau_p}^{\tau_{p+1}} b_j(s) ds\right)^2},
\end{equation}
\begin{equation}
\label{a1003}
D^{(2,j)}_{p+1,p}=\frac{1}{\sqrt{\Delta}}\int\limits_{\tau_p}^{\tau_{p+1}} b_j(s) ds,
\end{equation}
and $\xi_{p+1,p}^{(j)},$ $\eta_{p+1,p}^{(j)}$ $(j=1,2,3,\ p=0,1,\ldots,N-1)$ are mutually independent
standard Gaussian random variables.

Note that for $b_j(s)=C_j s^{\beta_j},$ $\beta_j\in (-1/2,0)$ $(j=1,2,3)$, $\tau_p=p\Delta$ 
$(p=0,1,\ldots,N)$
from (\ref{a1002}) and (\ref{a1003})
we obtain
\begin{equation}
\label{a1003hh}
D^{(1,j)}_{p+1,p}=\Delta^{\beta_j+1/2}\left|C_j\right|
\sqrt{\frac{(p+1)^{2\beta_j+1}-p^{2\beta_j+1}}{2\beta_j+1}-\left(\frac{(p+1)^{\beta_j+1}-p^{\beta_j+1}}
{\beta_j+1}\right)^2},
\end{equation}
\begin{equation}
\label{a1004hh}
D^{(2,j)}_{p+1,p}=\Delta^{\beta_j+1/2}\ C_j\
\frac{(p+1)^{\beta_j+1}-p^{\beta_j+1}}{\beta_j+1},
\end{equation}
where $p=0,1,\ldots,N-1.$

Let us substitute (\ref{a1000}) and (\ref{a1001}) into 
(\ref{1001ee}), (\ref{a2000}), (\ref{a2001}), and (\ref{a2001xxx})
\begin{equation}
\label{1001eee}
x_T=x_0 \exp\left(I^{a}_{T,0}-\frac{1}{2}\sum\limits_{j=1}^3 I^{b_j^2}_{T,0}
+\sum\limits_{j=1}^3
\sum\limits_{p=0}^{N-1}\left(D^{(1,j)}_{p+1,p}\cdot \eta_{p+1,p}^{(j)}+ 
D^{(2,j)}_{p+1,p}\cdot \xi_{p+1,p}^{(j)}\right)\right),
\end{equation}
$$
\tilde y_{p+1}=
\tilde y_p\Biggl(1+a(\tau_p)\Delta+\sqrt{\Delta}\sum\limits_{j=1}^3 b_j(\tau_p)\xi_{p+1,p}^{(j)}
+\Biggr.
$$
\begin{equation}
\label{a2000s}
\Biggl.+
\frac{\Delta}{2}\sum\limits_{l,j=1}^3 b_l(\tau_p)b_j(\tau_p)\left(
\xi_{p+1,p}^{(j)}\cdot \xi_{p+1,p}^{(l)},-{\bf 1}_{\{l=j\}}\right)\Biggr),
\end{equation}

\vspace{1mm}
$$
y_{p+1}=y_{p}\Biggl(1+
I^{a}_{\tau_{p+1},\tau_p}+
\sum\limits_{j=1}^3
\left(D^{(1,j)}_{p+1,p}\cdot \eta_{p+1,p}^{(j)}+ 
D^{(2,j)}_{p+1,p}\cdot \xi_{p+1,p}^{(j)}\right)
+\Biggr.
$$
$$
+\frac{1}{2}\sum\limits_{l,j=1}^3
\Biggl(\left(D^{(1,j)}_{p+1,p}\cdot \eta_{p+1,p}^{(j)}+ 
D^{(2,j)}_{p+1,p}\cdot \xi_{p+1,p}^{(j)}\right)
\left(D^{(1,l)}_{p+1,p}\cdot \eta_{p+1,p}^{(l)}+ 
D^{(2,l)}_{p+1,p}\cdot \xi_{p+1,p}^{(l)}\right)-
$$
\begin{equation}
\label{a2001s}
\Biggl.-
{\bf 1}_{\{l=j\}}
I^{b_l^2}_{\tau_{p+1},\tau_p}\Biggr)\Biggr),
\end{equation}

\begin{equation}
\label{a2001sxxx}
\hat y_{p+1}=\hat y_{p}\left(1+
I^{a}_{\tau_{p+1},\tau_p}+
\sum\limits_{j=1}^3
\left(D^{(1,j)}_{p+1,p}\cdot \eta_{p+1,p}^{(j)}+ 
D^{(2,j)}_{p+1,p}\cdot \xi_{p+1,p}^{(j)}\right)\right),
\end{equation}

\noindent
where $I^{a}_{T,0},$ $I^{a}_{\tau_{p+1},\tau_p},$ $I^{b_j^2}_{T,0},$ $I^{b_l^2}_{\tau_{p+1},\tau_p},$
$D^{(1,j)}_{p+1,p},$ $D^{(2,j)}_{p+1,p}$ are defined by (\ref{a4000}), 
(\ref{a1002}), (\ref{a1003}), (\ref{a1003hh}), (\ref{a1004hh}), 
and $\xi_{p+1,p}^{(j)},$ $\eta_{p+1,p}^{(j)}$ 
$(j=1,2,3,\ p=0,1,\ldots,N-1)$ are mutually independent
standard Gaussian random variables.

{\bf Numerical Experiment 3.}\ {\it For the It\^{o}
SDE {\rm (\ref{1000e})}
simulate {\rm 10000} {\rm (}$M=100$ groups of $L=100$ realizations each{\rm )} 
independent realizations $x_T^{(k,j)}$ of the random variable $x_T$ at $T=1,$ $x_0=1$ 
and $a(s)=-s^{-1/3},$ 
$b_1(s)=s^{-1/3},$ $b_2(s)=s^{-1/4},$ $b_3(s)=s^{-1/5}$
by the formula {\rm (\ref{1001eee})}
and the corresponding independent realizations $\tilde y_T^{(k,j)},$ $y_T^{(k,j)},$ and $\hat y_T^{(k,j)}$
of the random variables $\tilde y_T,$ $y_T,$ and $\hat y_T$ by the formulas {\rm (\ref{a2000s}),}
{\rm (\ref{a2001s}),} and {\rm (\ref{a2001sxxx}),} correspondently
{\rm ($k=1,\ldots,L,\ j=1,\ldots,M$)}.
Here $N=T/\Delta,$ $\Delta=2^{-15}.$ When performing
the first step by the Milstein method {\rm (\ref{a2000s}),} set $\tau_0=10^{-15}$ instead of $\tau_0=0.$
Calculate the estimates $\hat\varepsilon_{100,100}^{(1)},$ $\hat\varepsilon_{100,100}^{(2)},$ 
$\hat\varepsilon_{100,100}^{(3)}$ of the errors 
$\varepsilon^{(1)}={\sf M}\left\{|x_{T}-\tilde y_T|\right\},$ 
$\varepsilon^{(2)}={\sf M}\left\{|x_{T}-y_T|\right\},$ 
$\varepsilon^{(3)}={\sf M}\left\{|x_{T}-\hat y_T|\right\}$ 
by the formulas 
$$
\hat\varepsilon_{L,M}^{(1)}=\frac{1}{LM}\sum_{k=1}^L\sum_{j=1}^M
\left|x_T^{(k,j)}-\tilde y_T^{(k,j)}\right|,
$$
$$
\hat\varepsilon_{L,M}^{(2)}=\frac{1}{LM}\sum_{k=1}^L\sum_{j=1}^M
\left|x_T^{(k,j)}-y_T^{(k,j)}\right|,
$$
$$
\hat\varepsilon_{L,M}^{(3)}=\frac{1}{LM}\sum_{k=1}^L\sum_{j=1}^M
\left|x_T^{(k,j)}-\hat y_T^{(k,j)}\right|,
$$

\noindent
and a $90\%$ 
confidence intervals for $\varepsilon^{(1)},$ $\varepsilon^{(2)},$ $\varepsilon^{(3)},$
where $x_T^{(k,j)},$ $\tilde y_T^{(k,j)},$ $y_T^{(k,j)},$  and $\hat y_T^{(k,j)}$
for a fixed pair $(k,j)$ correspond to
the same realization 
of the Wiener processes ${\bf w}_t^{(l)},$ $t\in[0,T]$ $(l=1,2,3).$
Repeat the calculations for 
$\Delta=2^{-j},$\ $j=14, 13,\ldots, 3$ and plot 
the dependences $\hat\varepsilon_{100,100}^{(1)}(\Delta),$
$\hat\varepsilon_{100,100}^{(2)}(\Delta),$ $\hat\varepsilon_{100,100}^{(3)}(\Delta)$ with $90\%$ 
confidence intervals for $\varepsilon^{(1)},$ $\varepsilon^{(2)},$ $\varepsilon^{(3)}$
at the specified values 
of ${\Delta.}$ At the same time$,$ when doubling
the integration step$,$ use the formulas
$$
I^{b_j(j)}_{\tau_{p+2},\tau_p}=
I^{b_j(j)}_{\tau_{p+2},\tau_{p+1}}+I^{b_j(j)}_{\tau_{p+1},\tau_p},\ \ \ 
{\bf w}_{\tau_{p+2}}^{(j)}-{\bf w}_{\tau_p}^{(j)}=
\Delta {\bf w}_{p+1}^{(j)}+\Delta {\bf w}_{p}^{(j)}\ \ \ (j=1,2,3)
$$
w.~p.~{\rm 1}.
}

Shifting the initial point $\tau_0=0$ to $\tau_0 = 10^{-15}$ 
formally eliminates the singularity, 
transforming the It\^{o} SDE (\ref{1000e}) into a regular one on the interval 
$[10^{-15}, 1].$ However, the extreme dynamics of the coefficients in the vicinity 
of the starting point yield a massive Lipschitz constant. Consequently, the classical 
Milstein method (\ref{11vvv}) exhibits pre-asymptotic behavior with a gigantic error constant 
(resulting in large errors), 
which practically equates to a lack of convergence to the true solution of the 
original problem. Conversely, the proposed modified Milstein and Euler methods (\ref{200}) and 
(\ref{a2001xxx})
analytically 
smooth these dynamics within the stochastic integrals, ensuring stable and 
robust convergence under standard integration step sizes (see Fig.~3).
Note that the new modification (\ref{200}) of the Milstein method has an experimental 
strong order of convergence of approximately 0.60 for $\Delta\in [2^{-15}, 2^{-6}]$ (see Fig.~3), 
while the modified Euler method (\ref{a2001xxx}) has an experimental strong order of convergence 
of approximately 0.30 for $\Delta\in [2^{-15}, 2^{-6}]$ (see Fig.~3).
Note that the theoretical orders of strong convergence of methods
(\ref{200}) and (\ref{a2001xxx}) are $1/2$ and $1/3,$ respectively
(see Theorems 9, 7).

{\bf Numerical Experiment 3.1.}\ {\it Repeat Numerical Experiment {\rm 3}
for $a(s)=-s^{-1/5},$ 
$b_1(s)=s^{-1/5},$ $b_2(s)=s^{-1/6},$ $b_3(s)=s^{-1/7}.$}

The results of Numerical Experiment~3.1 are presented in Fig.~4. 
In this figure, we see that the new modification (\ref{200})
of the Milstein method and 
the modified Euler method (\ref{a2001xxx}) have strong orders of convergence 
of approximately  1.0  and 0.5 for $\Delta\in [2^{-15}, 2^{-6}]$, respectively (see Fig.~4).
At the same time, the classical Milstein method (\ref{11vvv}) 
exhibits pre-asymptotic behavior
behavior with large errors.
The theoretical orders of strong convergence of methods
(\ref{200}) and (\ref{a2001xxx}) are $9/10$ and $1/2-\varepsilon,$ respectively
(see Theorems 9, 7), where $\varepsilon$ is an arbitrary
small positive real number.

{\bf Numerical Experiment 4.}\ {\it Repeat Numerical Experiment {\rm 3}. 
In doing so, the first step of numerical scheme 
{\rm (\ref{a2000s})} is executed using formula {\rm (\ref{a2001s})}
$(\tau_0=0).$
}

The result of Numerical Experiment~4 is shown in Fig.~5. As expected, doing 
the first step of the Milstein method (\ref{11vvv}) using the new modification 
(\ref{200}) of the Milstein method greatly improved the behavior of the Milstein 
method (\ref{11vvv}) compared to Numerical Experiment~3.
However, the Milstein method (\ref{11vvv}) almost does not converge in the region 
of small integration step sizes $\Delta\in [2^{-14}, 2^{-15}]$ (see Fig.~5).

{\bf Numerical Experiment 4.1.}\ {\it Repeat Numerical Experiment {\rm 3.1}. 
In doing so, the first step of numerical scheme 
{\rm (\ref{a2000s})} is executed using formula {\rm (\ref{a2001s})}
$(\tau_0=0).$
}

The result of Numerical Experiment~4.1 is shown in Fig.~6.
In this figure, we can observe similar behavior of the Milstein method (\ref{11vvv}) 
to that in Fig.~5.

{\bf Numerical Experiment 4.2.}\ {\it Repeat Numerical Experiment {\rm 3.1}. 
In doing so, the first $10$ steps of numerical scheme 
{\rm (\ref{a2000s})} is executed using formula {\rm (\ref{a2001s})}
$(\tau_0=0).$
}

The result of Numerical Experiment~4.2 is shown in Fig.~7. In this case, 
doing the first 10 integration steps using the new modification (\ref{200})
of the Milstein method leads to an even better behavior of the Milstein 
method (\ref{11vvv}) (than in Fig.~6).
However, one should not have illusions about the capabilities of the Milstein 
method (\ref{11vvv}) when applied to singular It\^{o} SDEs, because, as noted in Remark~5, one 
can consider It\^{o} SDEs with a finite number of singularity points. 
This situation makes the Milstein method (\ref{11vvv}) practically inapplicable 
to the indicated class of It\^{o} SDEs, while the new modification (\ref{200})
of the Milstein method will have the same theoretical order of strong convergence as in Theorem 5.

{\bf Numerical Experiment 5.}\ {\it Repeat Numerical Experiment {\rm 3} with
$M=20,$ $L=100$ and
$a(s)=-s^{1/2},$ 
$b_1(s)=s^{1/3},$ $b_2(s)=s^{1/4},$ $b_3(s)=s^{1/5}.$
In doing so, the first step of numerical scheme 
{\rm (\ref{a2000s})} is executed with $\tau_0=0.$
}

The results of Numerical Experiment~5 are shown in Fig.~8. 
In this case, we consider a non-singular It\^{o} SDE with non-differentiable 
coefficients. Therefore, according to Remark~1, the new modification (\ref{200})
of the Milstein method has the order of strong convergence $1.0$, while the 
modified Euler method (\ref{a2001xxx}) has the order of strong convergence $0.5$ 
(see Fig.~8). Formally, the sufficient conditions ensuring the order 
of strong convergence $1.0$ are not satisfied for the Milstein method (\ref{11vvv}). 
It can be seen in Fig.~8 that the experimental order of strong convergence 
of the Milstein method (\ref{11vvv}) is between $0.5$ and $1.0$, that is, less than $1.0$, 
as expected.

{\bf Example~3.}\ Consider the following scalar non-commutative It\^{o} SDE
\begin{equation}
\label{15000}
x_t=x_0+\int\limits_0^t
\sigma_0(s)x_s ds+\int\limits_0^t
\sigma_1(s)x_s d{\bf w}_s^{(1)}+\int\limits_0^t\sigma_2(s)d{\bf w}_s^{(2)},\ \ \ x_0=x(0),
\end{equation}

\noindent
where $t\in [0, T],$
$\sigma_0(t)=a t^{\alpha},$
$\sigma_j(t)= c_j t^{\beta_j},$ $\alpha\in (-1,0)$, $\beta_j\in (-1/2,0)$ $(j=1,2)$,
$a,c_j\in {\bf R}$ $(j=1,2),$
${\bf w}_t^{(j)}$ $(j=1,2)$ are independent standard Wiener processes.

Let $\tau_p=p\Delta$ $(p=0,1,\ldots,N),$
$\Delta=T /N$ $(N>1)$ is a constant step of integration.

Application of the classical Milstein method (\ref{11vv}) to the It\^{o}
SDE (\ref{15000}) yields the following result
$$
\tilde y_{p+1}=\tilde y_p\left(1+\sigma_0(\tau_p)\Delta +\sigma_1(\tau_p)I_{\tau_{p+1},\tau_p}^{(1)}+
\sigma_1^2(\tau_p) I_{\tau_{p+1},\tau_p}^{(11)}\right)+
$$
\begin{equation}
\label{20000}
+
\sigma_2(\tau_p)I_{\tau_{p+1},\tau_p}^{(2)}+\sigma_1(\tau_p)\sigma_2(\tau_p)I_{\tau_{p+1},\tau_p}^{(21)},
\end{equation}

\vspace{1mm}
\noindent
where iterated It\^{o} stochastic integrals 
$I_{\tau_{p+1},\tau_p}^{(1)},$ $I_{\tau_{p+1},\tau_p}^{(2)},$ $I_{\tau_{p+1},\tau_p}^{(11)},$
$I_{\tau_{p+1},\tau_p}^{(21)}$ are defined by formula (\ref{int1})
and for their mean-square approximation formulas (\ref{15000a})--(\ref{15000d})
can be used.

At the same time, the new modification (\ref{11}) of Milstein method 
for the It\^{o} SDE (\ref{15000}) takes the form
\begin{equation}
\label{20001}
y_{p+1}=y_p\left(1+I^{\sigma_0}_{\tau_{p+1},\tau_p}+ 
I^{\sigma_1(1)}_{\tau_{p+1},\tau_p}+
I_{\tau_{p+1},\tau_p}^{\sigma_1\sigma_1(11)}\right)+
I^{\sigma_2(2)}_{\tau_{p+1},\tau_p}+I_{\tau_{p+1},\tau_p}^{\sigma_2\sigma_1(21)},
\end{equation}

\noindent
where iterated It\^{o} stochastic integrals $I^{\sigma_1(1)}_{\tau_{p+1},\tau_p},$
$I^{\sigma_2(2)}_{\tau_{p+1},\tau_p},$ $I_{\tau_{p+1},\tau_p}^{\sigma_1\sigma_1(11)},$
$I_{\tau_{p+1},\tau_p}^{\sigma_2\sigma_1(21)}$ are difined by (\ref{11b}) and (\ref{11c}),
and the integral $I^{\sigma_0}_{\tau_{p+1},\tau_p}$ has the form (\ref{11a}).

Furthermore, the modified Euler's method (\ref{euler1}) 
for the It\^{o} SDE (\ref{15000}) can be written as 
\begin{equation}
\label{20001xxx}
\hat y_{p+1}=\hat y_p\left(1+I^{\sigma_0}_{\tau_{p+1},\tau_p}+ 
I^{\sigma_1(1)}_{\tau_{p+1},\tau_p}\right)+
I^{\sigma_2(2)}_{\tau_{p+1},\tau_p}.
\end{equation}

To obtain the mean-square approximations of iterated It\^{o}
stochastic integrals 
$I^{\sigma_1(1)}_{\tau_{p+1},\tau_p},$
$I^{\sigma_2(2)}_{\tau_{p+1},\tau_p},$ $I_{\tau_{p+1},\tau_p}^{\sigma_1\sigma_1(11)},$
$I_{\tau_{p+1},\tau_p}^{\sigma_2\sigma_1(21)},$
one can apply formulas (\ref{100})--(\ref{106}), 
(\ref{1000}), (\ref{1001}). Moreover,
\begin{equation}
\label{xx1}
I_{\tau_{p+1},\tau_p}^{\sigma_1\sigma_1(11)q_1}=\frac{1}{2}\left(\left(
I^{\sigma_1(1)q_1}_{\tau_{p+1},\tau_p}\right)^2-
I^{\sigma_1^2}_{\tau_{p+1},\tau_p}\right),
\end{equation}
\begin{equation}
\label{xx2}
I^{\sigma_2\sigma_1(21)q_2}_{\tau_{p+1},\tau_p}
=\sum_{k,r=0}^{q_2}
C_{r k}^{\sigma_l\sigma_j}\zeta_{k}^{(l)}\zeta_{r}^{(j)},
\end{equation}
where $q_1, q_2\in {\bf N};$ another notations in (\ref{xx1}), (\ref{xx2})
are the same as in (\ref{100})--(\ref{103}), (\ref{1000}), (\ref{1001}).

The exact solution to It\^{o} SDE (\ref{15000}) is known and 
has the following form
\begin{equation}
\label{20002}
x_t=y_t\left(x_0+\int\limits_0^t y_s^{-1}\sigma_2(s)d{\bf w}_s^{(2)}\right),
\end{equation}
where
\begin{equation}
\label{20003}
y_t=\exp\left(\int\limits_0^t\left(\sigma_0(s)-\frac{1}{2}\sigma_1^2(s)\right)ds +
\int\limits_0^t \sigma_1(s)d{\bf w}_s^{(1)}\right).
\end{equation}

From (\ref{20002}) and (\ref{20003}) it follows that the problem of numerical simulation 
of the exact solution to It\^{o} SDE (\ref{15000}) is even more complex 
than the numerical implementation of formula (\ref{20001}).
Therefore, as a reference solution, we choose a random sequence obtained 
using the numerical implementation of formula (\ref{20001}) with a sufficiently small integration 
step $\Delta=2^{-15}T.$ Further, the reference solution $x^{ref}_T$ at time moment
$T$ will be compared with $\tilde y_T,$ $y_T$ obtained using 
classical Milstein method (\ref{20000}) and its new modification (\ref{20001})
with integration steps $\Delta=2^{-k}T,$ where $k=2,3,\ldots,14.$
In this case, the iterated It\^{o} stochastic integrals included in (\ref{20000}) and (\ref{20001}) 
are first modeled at the smallest integration step $\Delta=2^{-15}T.$ Then, when the 
integration step is doubled, the following transition formulas for the 
iterated It\^{o} stochastic integrals are used
\begin{equation}
\label{20004}
I_{\tau_{p+2},\tau_p}^{\sigma_2\sigma_1(21)}=
I_{\tau_{p+2},\tau_{p+1}}^{\sigma_2\sigma_1(21)}+
I_{\tau_{p+1},\tau_p}^{\sigma_2\sigma_1(21)}+
I_{\tau_{p+2},\tau_{p+1}}^{\sigma_1(1)}\cdot
I_{\tau_{p+1},\tau_p}^{\sigma_2(2)},
\end{equation}
\begin{equation}
\label{20005}
I_{\tau_{p+2},\tau_p}^{\sigma_j(j)}=I_{\tau_{p+2},\tau_{p+1}}^{\sigma_j(j)}+
I_{\tau_{p+1},\tau_p}^{\sigma_j(j)}\ \ \ (j=1,2),
\end{equation}
\begin{equation}
\label{20006}
I_{\tau_{p+2},\tau_p}^{(21)}=
I_{\tau_{p+2},\tau_{p+1}}^{(21)}+
I_{\tau_{p+1},\tau_p}^{(21)}+
I_{\tau_{p+2},\tau_{p+1}}^{(1)}\cdot
I_{\tau_{p+1},\tau_p}^{(2)},
\end{equation}
\begin{equation}
\label{20007}
I_{\tau_{p+2},\tau_p}^{(j)}=I_{\tau_{p+2},\tau_{p+1}}^{(j)}+
I_{\tau_{p+1},\tau_p}^{(j)}\ \ \ (j=1,2)
\end{equation}
w.~p.~1.

\vspace{2mm}

{\bf Numerical Experiment 6.}\ {\it For the It\^{o}
SDE {\rm (\ref{15000})}
simulate {\rm 2000} {\rm (}$M=20$ groups of $L=100$ realizations each{\rm )} 
independent realizations $x_T^{ref(k,j)}$ of the random variable $x_T^{ref}$ at $T=1, x_0=1$ 
and $\sigma_0(s)=-s^{-1/2},
\sigma_1(s)=s^{-1/3}, \sigma_2(s)=s^{-1/4}$
by the formula {\rm (\ref{20001})} with $\Delta=2^{-15}T$
and the corresponding independent realizations $\tilde y_T^{(k,j)},$ $y_T^{(k,j)},$ and $\hat y_T^{(k,j)},$
of the random variables $\tilde y_T,$ $y_T,$ and $\hat y_T,$ by the formulas {\rm (\ref{20000}),}
{\rm (\ref{20001}),} and {\rm (\ref{20001xxx})} with $\Delta=2^{-14}T,$ correspondently
{\rm ($k=1,\ldots,L,\ j=1,\ldots,M$)}.
When performing
the first step by scheme {\rm (\ref{20000}),} set $\tau_0=10^{-9}$ instead of $\tau_0=0.$
Calculate the estimates $\hat\varepsilon_{100,20}^{(1)},$ $\hat\varepsilon_{100,20}^{(2)},$
$\hat\varepsilon_{100,20}^{(3)}$ of the errors 
$$
\varepsilon^{(1)}={\sf M}\bigl\{|x_{T}^{ref}-\tilde y_T|\bigr\},\ \ \
\varepsilon^{(2)}={\sf M}\bigl\{|x_{T}^{ref}-y_T|\bigr\},\ \ \
\varepsilon^{(2)}={\sf M}\bigl\{|x_{T}^{ref}-\hat y_T|\bigr\}
$$

\noindent 
by the formulas 
$$
\hat\varepsilon_{L,M}^{(1)}=\frac{1}{LM}\sum_{k=1}^L\sum_{j=1}^M
\left|x_T^{ref(k,j)}-\tilde y_T^{(k,j)}\right|,
$$
$$
\hat\varepsilon_{L,M}^{(2)}=\frac{1}{LM}\sum_{k=1}^L\sum_{j=1}^M
\left|x_T^{ref(k,j)}-y_T^{(k,j)}\right|
$$
$$
\hat\varepsilon_{L,M}^{(3)}=\frac{1}{LM}\sum_{k=1}^L\sum_{j=1}^M
\left|x_T^{ref(k,j)}-\hat y_T^{(k,j)}\right|
$$

\noindent
and a $90\%$ 
confidence intervals for $\varepsilon^{(1)},$ $\varepsilon^{(2)},$ $\varepsilon^{(2)}.$
For numerical modeling of the iterated It\^{o} stochastic integrals
$I^{\sigma_1(1)}_{\tau_{p+1},\tau_p},$
$I^{\sigma_2(2)}_{\tau_{p+1},\tau_p},$ $I_{\tau_{p+1},\tau_p}^{\sigma_1\sigma_1(11)},$
$I_{\tau_{p+1},\tau_p}^{\sigma_2\sigma_1(21)},$
$I^{(j)}_{\tau_{p+1},\tau_p}$ $(j=1,2),$
$I_{\tau_{p+1},\tau_p}^{(11)},$
$I_{\tau_{p+1},\tau_p}^{(21)}$
use formulas {\rm (\ref{100})--(\ref{103}),} 
{\rm (\ref{1000})--(\ref{15000b}),} {\rm (\ref{15000d}),}
{\rm (\ref{xx1})}, {\rm (\ref{xx2}),} 
where $q,q_1,q_2=10.$
To recalculate the above iterated It\^{o} stochastic integrals
when doubling the integration step, use formulas 
{\rm (\ref{20004})}--{\rm (\ref{20007}).} 
Repeat the calculations for 
$\Delta=2^{-j},$\ $j=13, 12, 11,\ldots,2$ and plot 
the dependences $\hat\varepsilon_{100,20}^{(1)}(\Delta),$
$\hat\varepsilon_{100,20}^{(2)}(\Delta),$ $\hat\varepsilon_{100,20}^{(3)}(\Delta)$ with $90\%$ 
confidence intervals for $\varepsilon^{(1)},$ $\varepsilon^{(2)},$ 
$\varepsilon^{(3)}$ at the specified values 
of ${\Delta.}$ 
}

\vspace{2mm}
                     
The results of Numerical Experiment~6 are shown in Fig.~9.
Just as in Numerical Experiment~3,
the classical 
Milstein method (\ref{11vvv}) exhibits pre-asymptotic behavior with large errors,
because of the extreme dynamics of the coefficients of the It\^{o} SDE in the vicinity 
of the starting point.
From the other hand, the proposed modified Milstein and Euler methods (\ref{200}) and 
(\ref{a2001xxx})
analytically smooth these dynamics within the stochastic integrals, ensuring 
convergence under standard integration step sizes (see Fig.~9).
Note that the new modification (\ref{200}) of the Milstein method has an experimental 
strong order of convergence of approximately 0.5 (see Fig.~9), 
while the modified Euler method (\ref{a2001xxx}) has an experimental strong order of convergence 
of approximately 0.27 (see Fig.~9).
The theoretical orders of strong convergence of methods
(\ref{200}) and (\ref{a2001xxx}) are $1/2$ and $1/3,$ respectively
(see Theorems 9, 7).

\vspace{2mm}

{\bf Numerical Experiment 6.1.}\ {\it Repeat Numerical Experiment {\rm 6}
for $a(s)=-s^{-1/8},$ 
$b_1(s)=s^{-1/6},$ $b_2(s)=s^{-1/7}.$}

\vspace{2mm}

The results of Numerical Experiment~6.1 are presented in Fig.~10. 
In this figure, we see that the new modification (\ref{200})
of the Milstein method and 
the modified Euler method (\ref{a2001xxx}) have experimental strong orders of convergence 
of approximately 1.0 and 0.5, respectively (see Fig.~10).
At the same time, the classical Milstein method (\ref{11vvv}) 
exhibits pre-asymptotic behavior
with large errors.
The theoretical orders of strong convergence of methods
(\ref{200}) and (\ref{a2001xxx}) are $1-\varepsilon$ and $1/2-\varepsilon,$ respectively
(see Theorems 9, 7), where $\varepsilon$ is an arbitrary
small positive real number.

\vspace{2mm}

{\bf Numerical Experiment 7.}\ {\it Repeat Numerical Experiment {\rm 6}. 
In doing so, the first step of numerical scheme 
{\rm (\ref{20000})} is executed using formula {\rm (\ref{20001})}
$(\tau_0=0).$
}

\vspace{2mm}

The result of Numerical Experiment~7 is shown in Fig.~11. As expected, doing 
the first step of the Milstein method (\ref{11vvv}) using the new modification 
(\ref{200}) of the Milstein method greatly improved the behavior of the Milstein 
method (\ref{11vvv}) compared to Numerical Experiment~6.

\vspace{2mm}

{\bf Numerical Experiment 7.1.}\ {\it Repeat Numerical Experiment {\rm 7}
for $a(s)=-s^{-1/8},$ 
$b_1(s)=s^{-1/6},$ $b_2(s)=s^{-1/7}.$}

\vspace{2mm}

The result of Numerical Experiment~7.1 is shown in Fig.~12.
In this figure, we can observe similar behavior of the Milstein method (\ref{11vvv}) 
to that in Fig.~11.

\vspace{2mm}

{\bf Numerical Experiment 7.2.}\ {\it Repeat Numerical Experiment {\rm 7.1}.
In doing so, the first $10$ steps of numerical scheme 
{\rm (\ref{20000})} is executed using formula {\rm (\ref{20001})}
$(\tau_0=0).$
}

\vspace{2mm}

The result of Numerical Experiment~7.2 is shown in Fig.~13. In this case, 
doing the first 10 integration steps using the new modification (\ref{200})
of the Milstein method leads to an even better behavior of the Milstein 
method (\ref{11vvv}) (than in Fig.~12).
As noted in Remark~5, one 
can consider It\^{o} SDEs with a finite number of singularity points.
Therefore, one should not have illusions about the potential of 
the Milstein method (\ref{11vvv}) when applied to singular It\^{o} SDEs.

\vspace{2mm}

{\bf Numerical Experiment 8.}\ {\it Repeat Numerical Experiment {\rm 6} with
$\sigma_0(s)=-s^{1/2},
\sigma_1(s)=s^{1/3}, \sigma_2(s)=s^{1/4}.$
In doing so, the first step of numerical scheme 
{\rm (\ref{20000})} is executed with $\tau_0=0.$
}

\vspace{2mm}

The results of Numerical Experiment~8 are shown in Fig.~14. 
In this case, we have a non-singular It\^{o} SDE with non-differentiable 
coefficients. Thus, according to Remark~1, the new modification (\ref{200})
of the Milstein method has the order of strong convergence $1.0$, while the 
modified Euler method (\ref{a2001xxx}) has the order of strong convergence $0.5$ 
(see Fig.~14). In the case under consideration, sufficient conditions ensuring the order 
of strong convergence $1.0$ are not satisfied for the Milstein method (\ref{11vvv}). 
It can be seen in Fig.~14 that the experimental order of strong convergence 
of the Milstein method (\ref{11vvv}) is between $0.5$ and $1.0$, that is, less than $1.0$, 
as expected.

\vspace{20mm}

\section{Figures}

This section presents the figures illustrating the results of the numerical 
experiments throughout this paper.

\newcommand{\experimentpagetop}{\vspace*{-7mm}}
\newcommand{\experimentgap}{\vspace{16pt}}
\newcommand{\experimentfigure}[3]{%
  \noindent\begin{minipage}{\textwidth}
  \centering
  \includegraphics[width=\textwidth,height=0.45\textheight,
    keepaspectratio]{#1}\par
  \captionsetup{labelsep=period,skip=10pt}%
  \captionof{figure}{Results of Experiment #2}\label{#3}
  \end{minipage}\par
}
\clearpage
\begingroup
\centering
\setlength{\parskip}{0pt}
\experimentpagetop
\experimentfigure{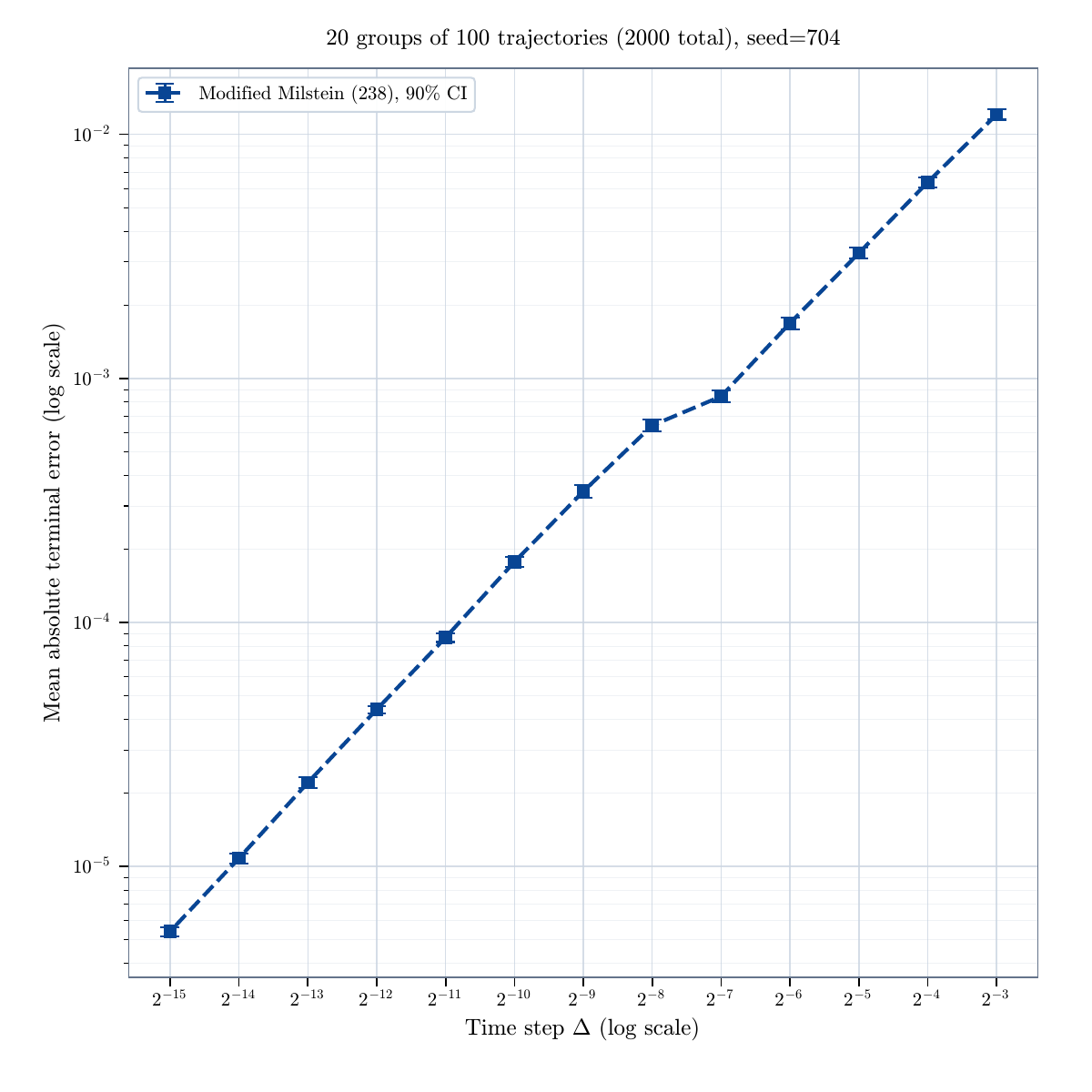}{1}{fig:exp-1}
\experimentgap
\experimentfigure{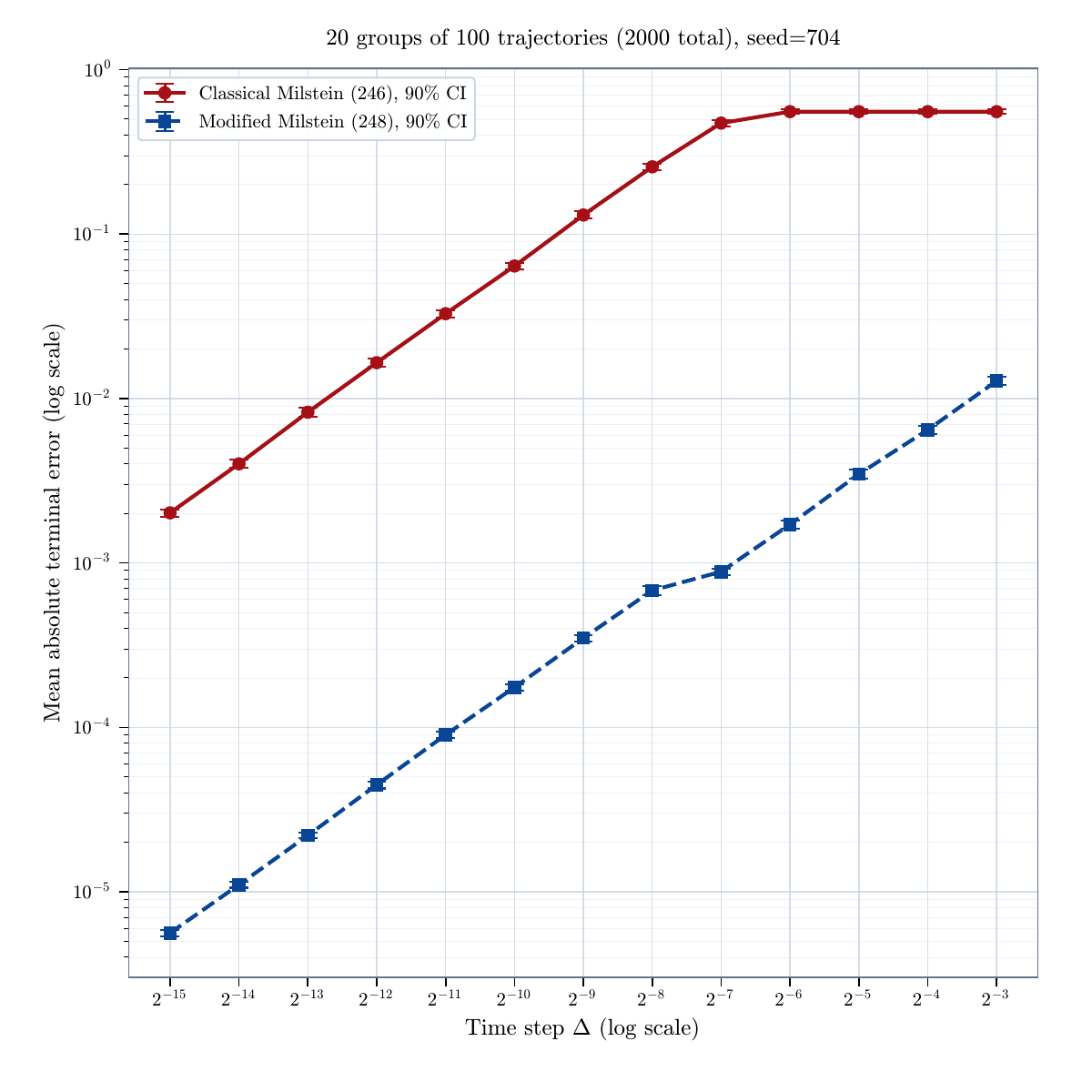}{2}{fig:exp-2}
\clearpage
\experimentpagetop
\experimentfigure{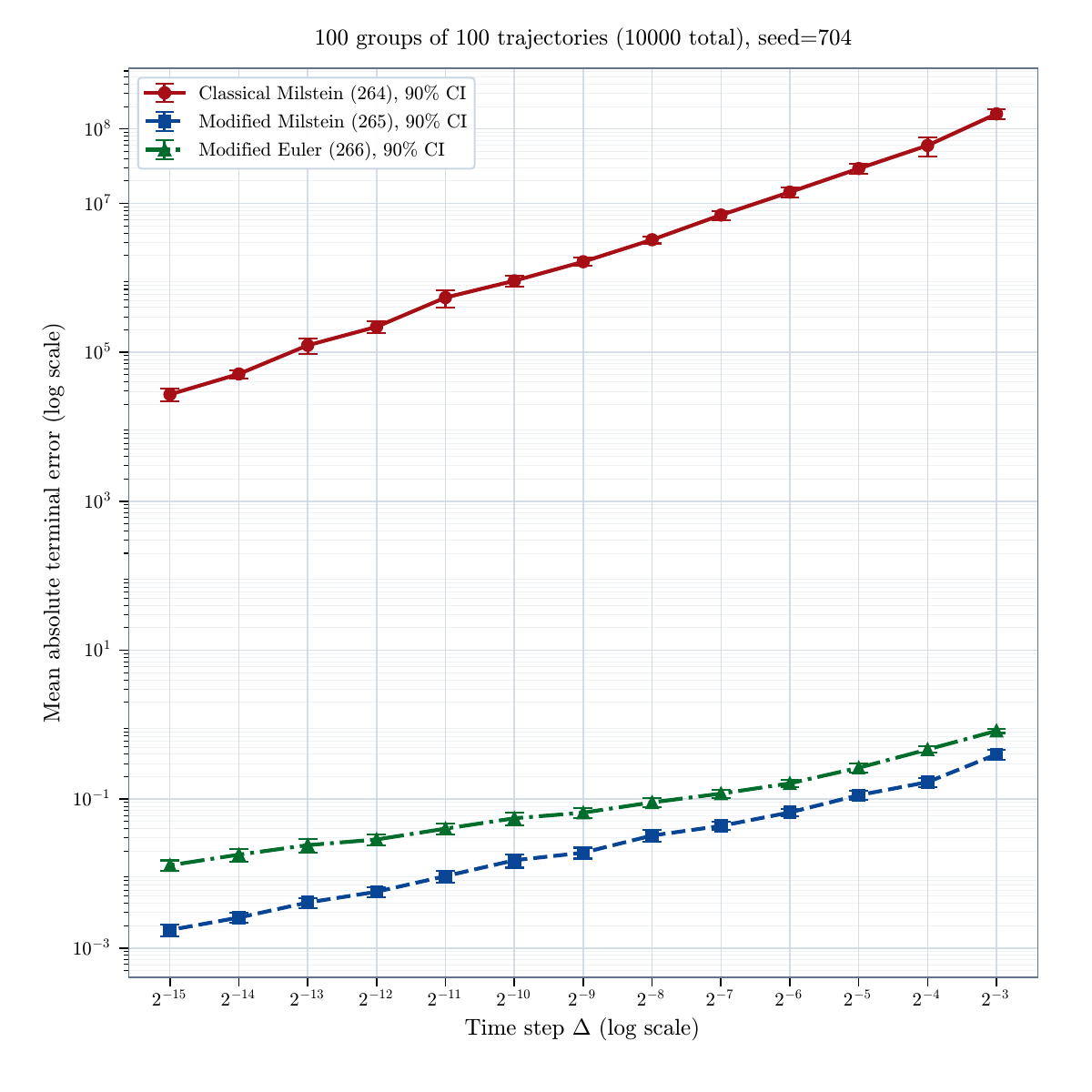}{3}{fig:exp-3-9}
\experimentgap
\experimentfigure{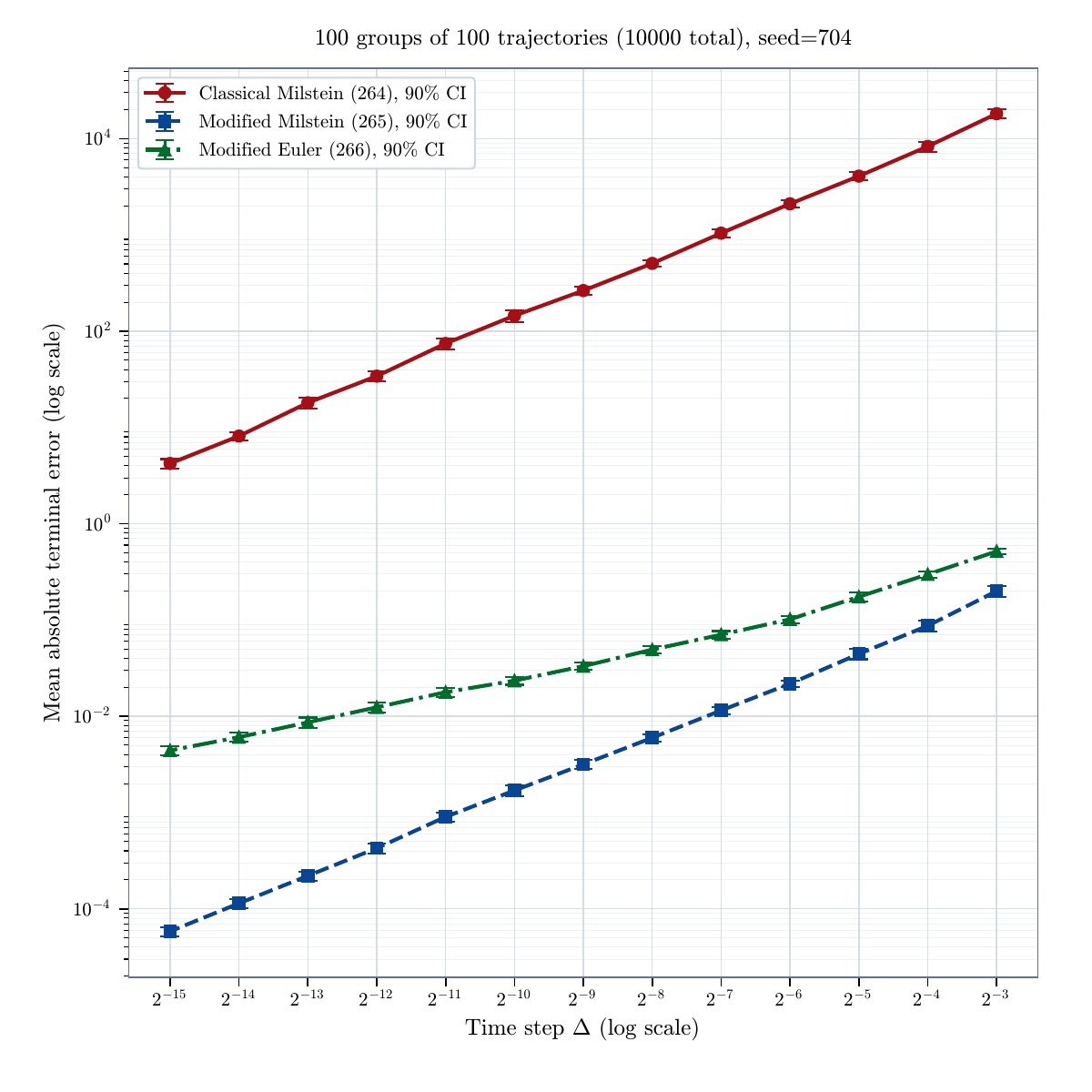}{3.1}{fig:exp-3-1-9}
\clearpage
\experimentpagetop
\experimentfigure{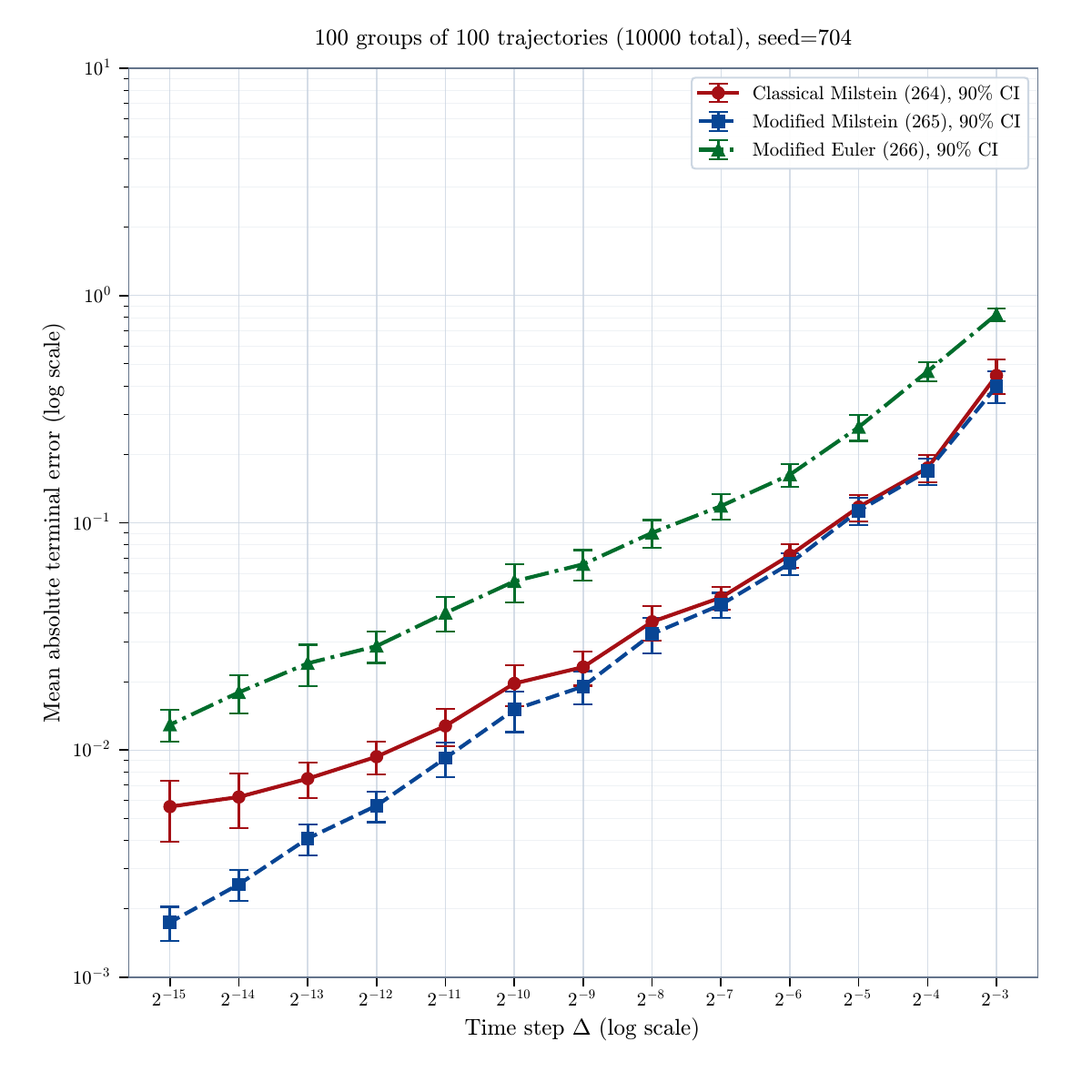}{4}{fig:exp-4}
\experimentgap
\experimentfigure{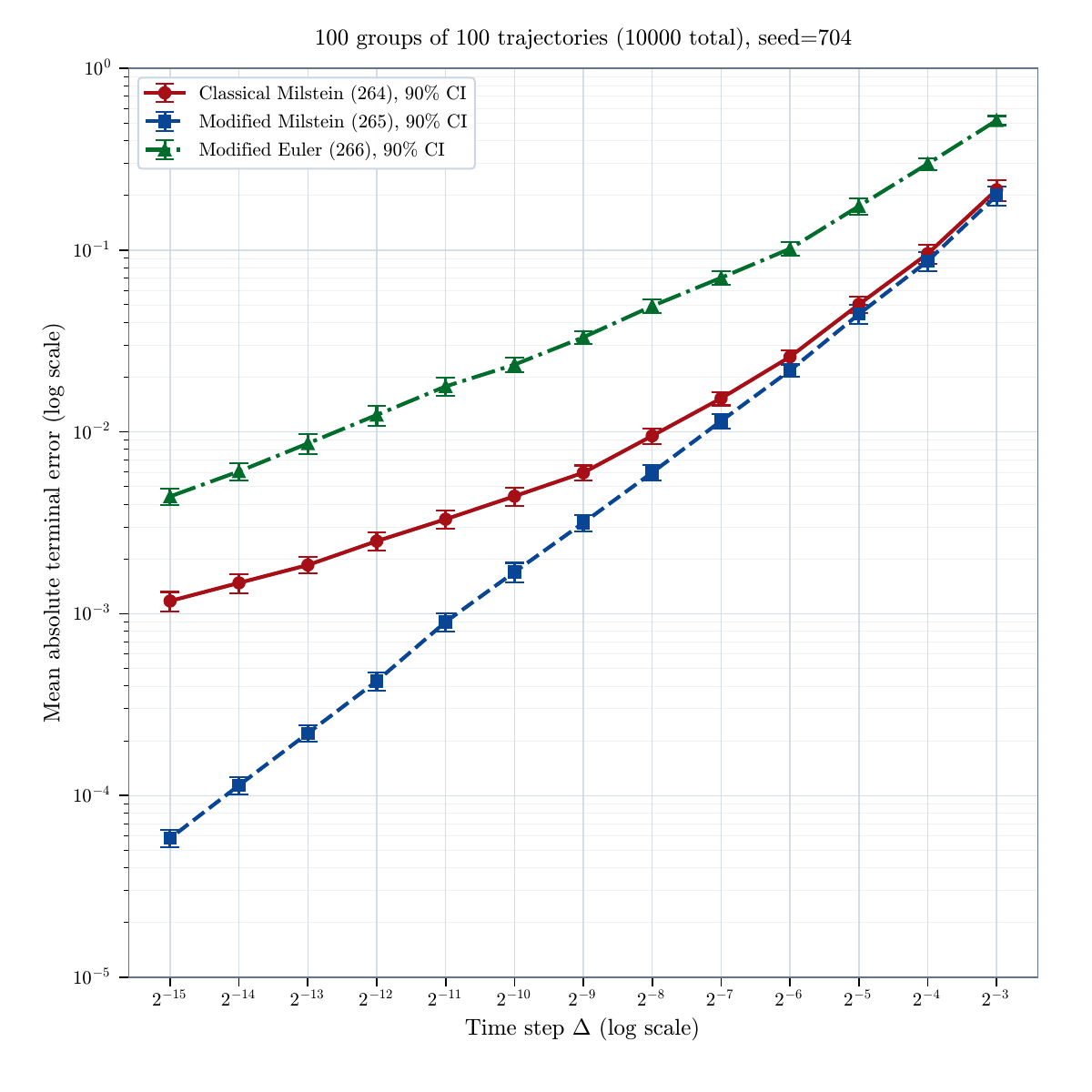}{4.1}{fig:exp-4-1}
\clearpage
\experimentpagetop
\experimentfigure{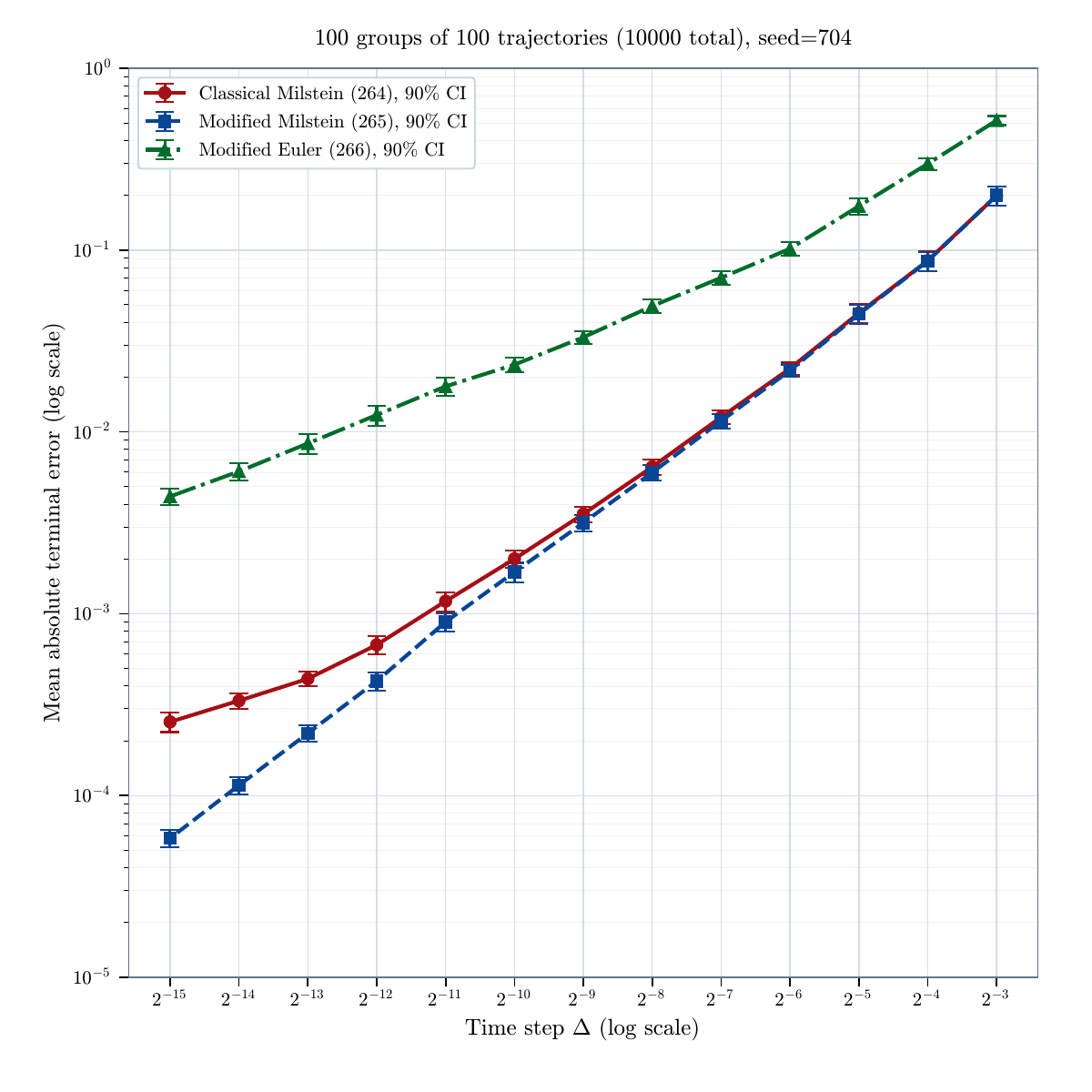}{4.2}{fig:exp-4-2}
\experimentgap
\experimentfigure{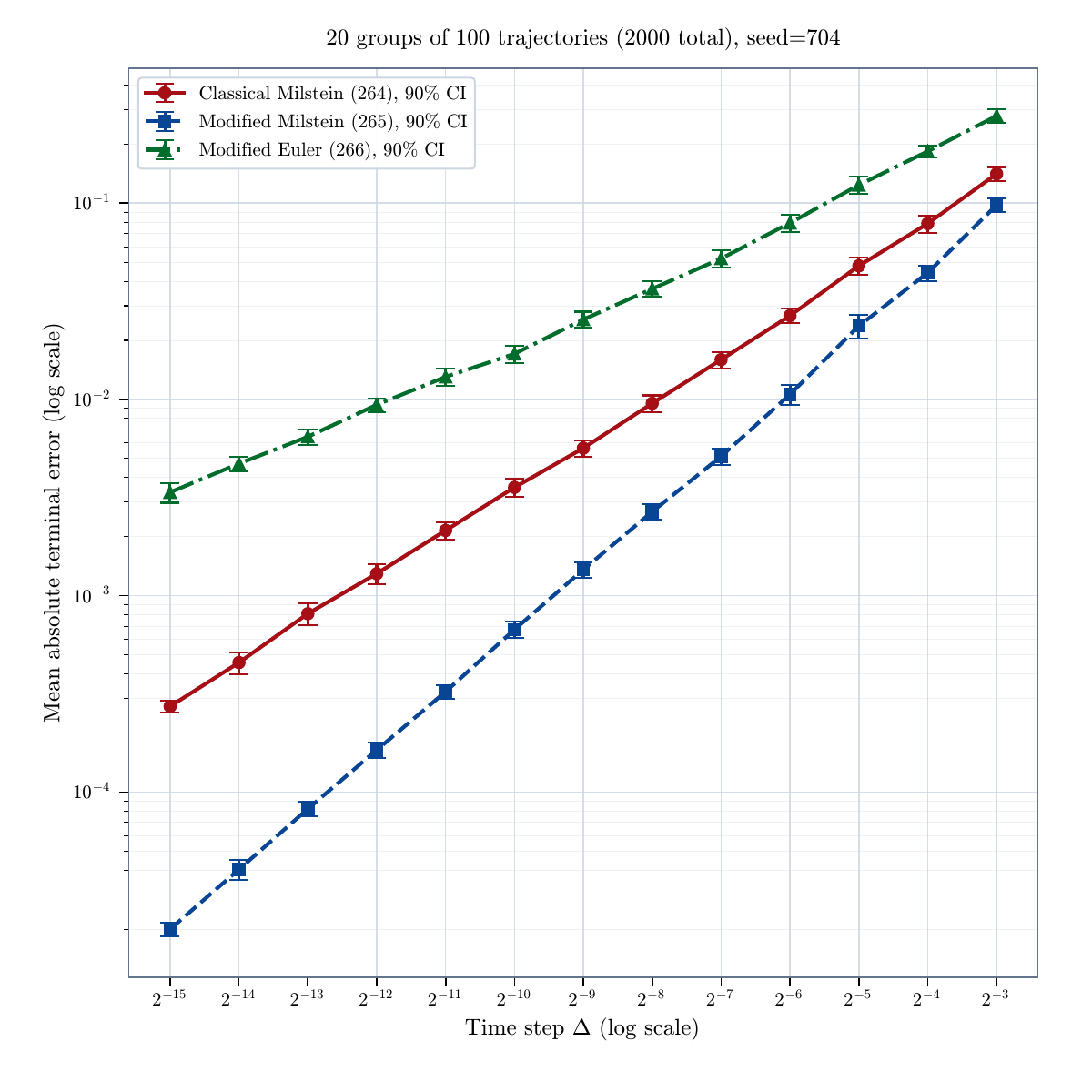}{5}{fig:exp-5}
\clearpage
\experimentpagetop
\experimentfigure{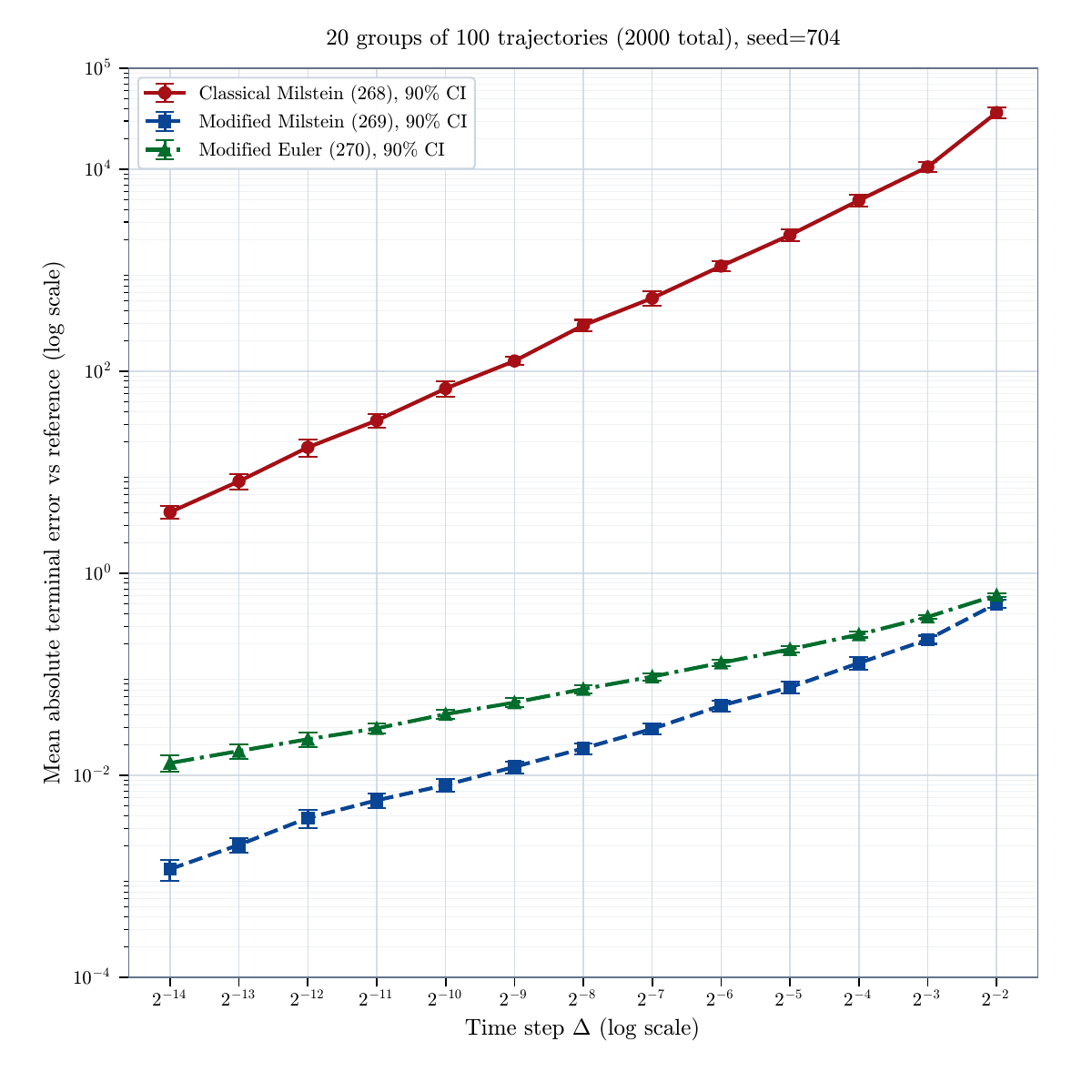}{6}{fig:exp-6}
\experimentgap
\experimentfigure{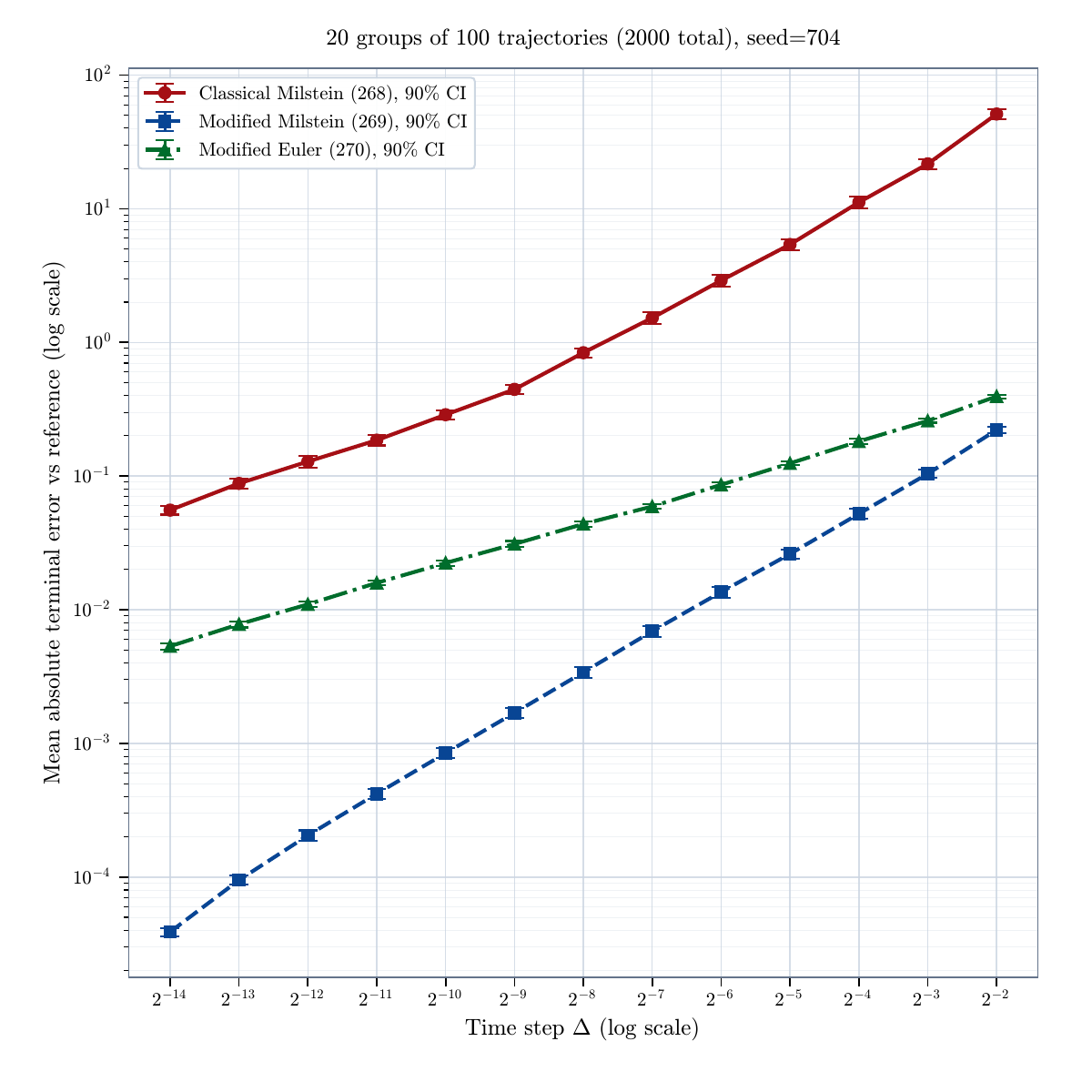}{6.1}{fig:exp-6-1}
\clearpage
\experimentpagetop
\experimentfigure{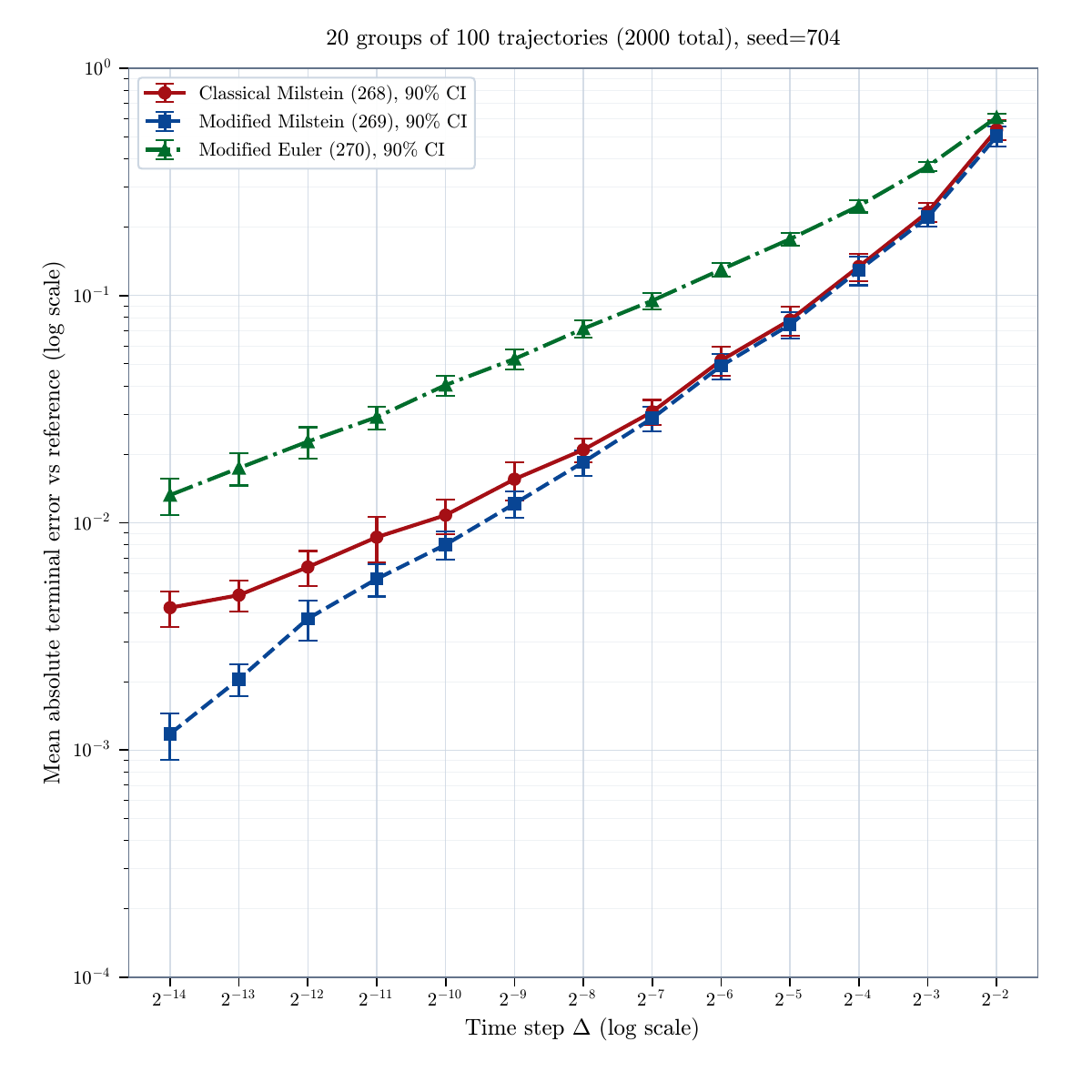}{7}{fig:exp-7}
\experimentgap
\experimentfigure{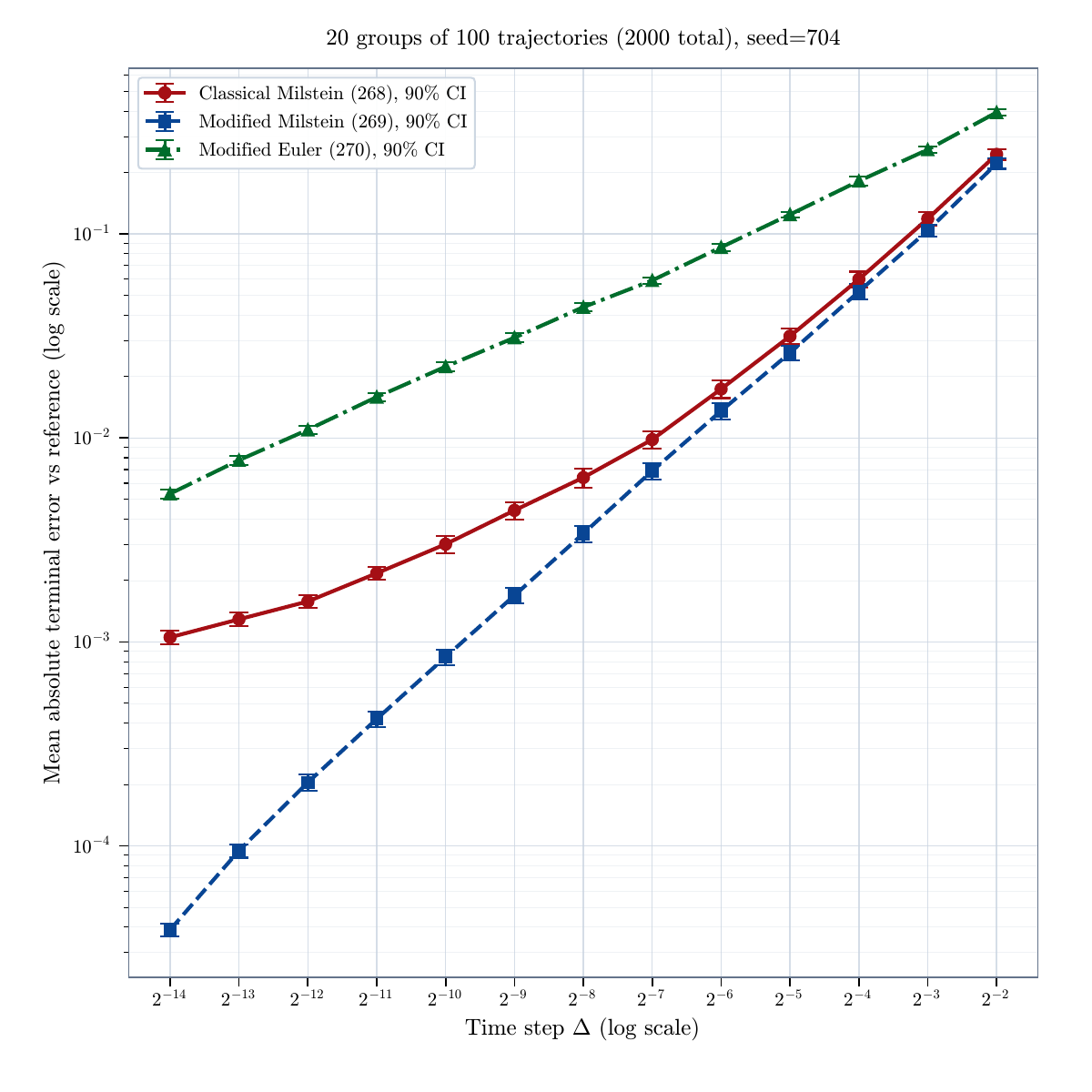}{7.1}{fig:exp-7-1}
\clearpage
\experimentpagetop
\experimentfigure{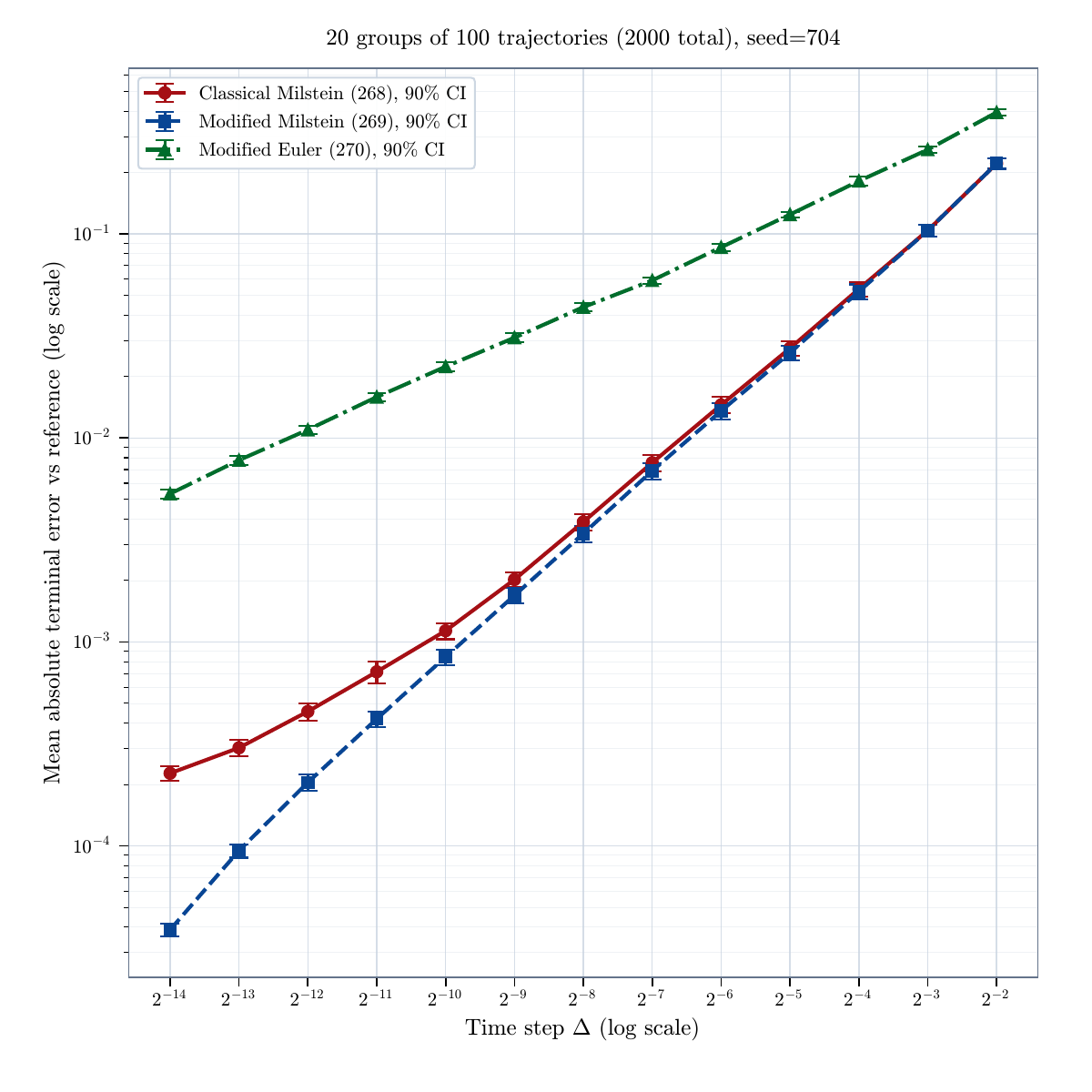}{7.2}{fig:exp-7-2}
\experimentgap
\experimentfigure{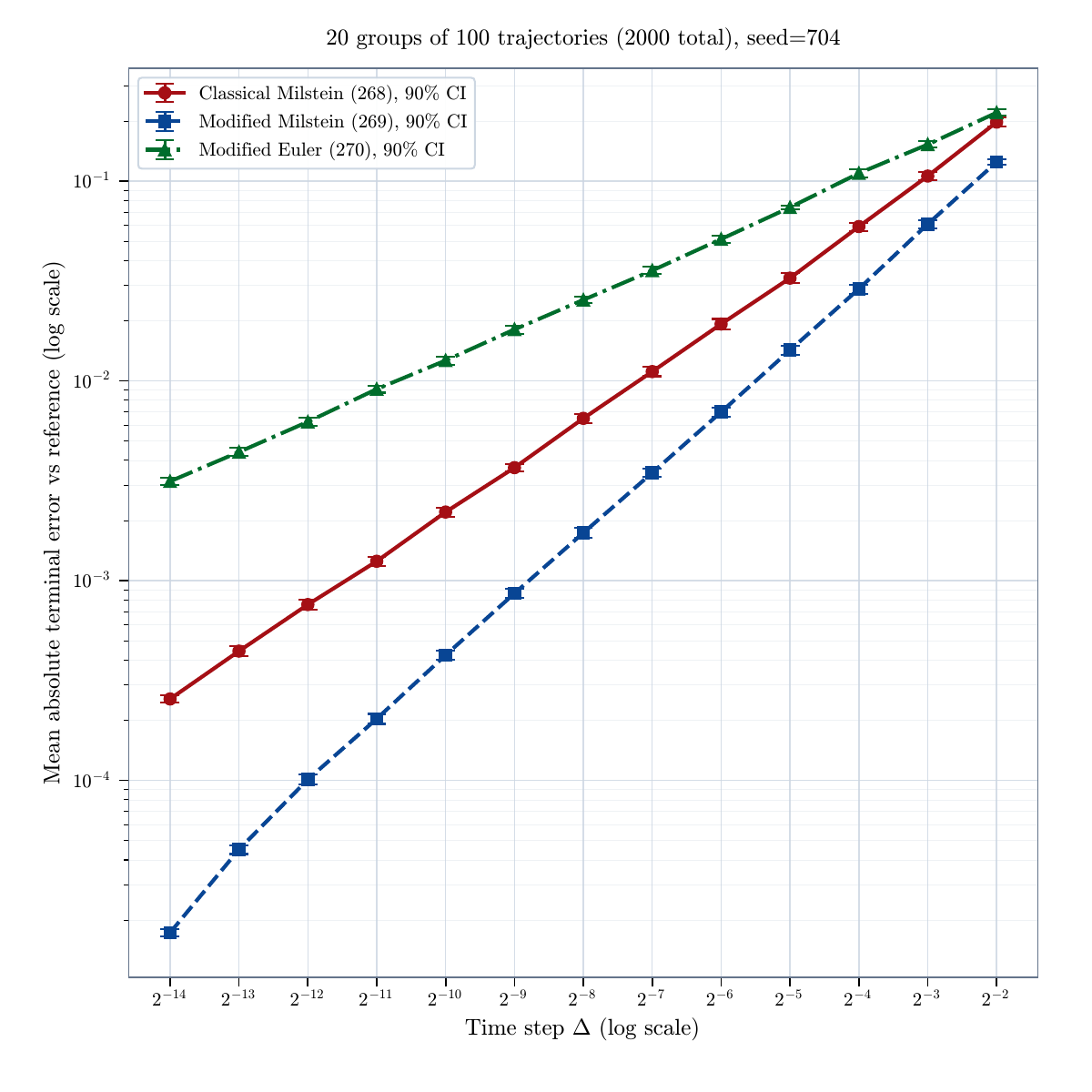}{8}{fig:exp-8}
\clearpage
\endgroup

\linespread{0.7}

\end{document}